\documentclass[final]{amsart}
\usepackage[utf8]{inputenc}
\usepackage{amsmath, amsthm} 
\usepackage{amsfonts,amssymb}
\usepackage{graphicx,epstopdf}
\usepackage{hyperref}
\usepackage{cite}

\newcommand{\RR}{\mathbb{R}}

\renewcommand{\k}{{\bf k}}
\renewcommand{\u}{{\bf u}}
\renewcommand{\v}{{\bf v}}
\newcommand{\F}{{\sf F}}

\DeclareMathOperator{\dd}{d\!}
\DeclareMathOperator{\sech}{sech}

\newtheorem{thm}{Theorem}[section]
\newtheorem{conjecture}[thm]{Conjecture}
\newtheorem{proposition}[thm]{Proposition}
\newtheorem{remark}[thm]{Remark}

\begin{document}
\title[Numerical study of the SGN and WGN equations]{Numerical study of the Serre-Green-Naghdi equations and a fully dispersive counterpart}
\author{Vincent Duch\^ene}
\address{Univ Rennes, CNRS, IRMAR - UMR 6625, F-35000 Rennes, France\\
E-mail \url{vincent.duchene@univ-rennes1.fr}}

\author{Christian Klein}
\address{Institut de Math\'ematiques de Bourgogne, UMR 5584\\
Institut Universitaire de France\\
                Universit\'e de Bourgogne-Franche-Comt\'e, 9 avenue Alain Savary, 21078 Dijon
                Cedex, France\\
    E-mail \url{Christian.Klein@u-bourgogne.fr}}

\date{\today}

\begin{abstract}
	We perform numerical experiments on the Serre-Green-Naghdi (SGN) equations
	and a fully dispersive ``Whitham-Green-Naghdi'' (WGN) counterpart in dimension 1. 
	In particular, solitary wave solutions of the WGN equations are constructed and their stability, 
	along with the explicit ones of the SGN equations, is studied. Additionally, 
	the emergence of modulated oscillations and 
	the possibility of a blow-up of solutions in various situations is investigated. 
	
	We argue that a simple numerical scheme based on a Fourier spectral method 
	combined with the Krylov subspace iterative technique GMRES
	to address the elliptic problem and a fourth order explicit Runge-Kutta scheme in time
	allows to address efficiently even computationally challenging problems.
\end{abstract}


\thanks{This work is partially supported by 
the ANR-FWF project ANuI - ANR-17-CE40-0035, the isite BFC project 
NAANoD, the  EIPHI Graduate School (contract ANR-17-EURE-0002) and by the 
European Union Horizon 2020 research and innovation program under the 
Marie Sklodowska-Curie RISE 2017 grant agreement no. 778010 IPaDEGAN}
\maketitle

\section{Introduction}\label{S.introduction}
\subsection{Motivation}

The Serre-Green-Naghdi (SGN) model is a popular model for the propagation of surface gravity waves in coastal oceanography. It is expected to provide a reasonable approximation of the response to gravity forces of a layer of homogeneous incompressible fluid with a free surface (hereafter referred to as the {\em water waves problem}) in the so-called {\em shallow-water regime}, that is for weakly dispersive but possibly strongly nonlinear flows. It has been derived and studied by many authors, including~\cite{Serre53,SuGardner69,GreenNaghdi76,MilesSalmon85,BazdenkovMorozovPogutse87,Seabra-SantosRenouardTemperville87,Salmon88,WeiKirbyGrilliEtAl95,CamassaHolmLevermore96,KimBaiErtekinEtAl01,Barthelemy04,LannesBonneton09,ClamondDutykh12,Ionescu-Kruse12,FedotovaKhakimzyanovDutykh14,GavrilyukKalischKhorsand15}.
Its rigorous justification as an asymptotic model in the shallow-water limit has been obtained in~\cite{Makarenko86,Li06,Israwi11,Alvarez-SamaniegoLannes08a,FujiwaraIguchi15}. In addition to its validity as a model for the water waves problem, the SGN equations have attracted interest as they are natural dispersive generalizations of the equations for isentropic compressible flows  and as such can be studied through the Lagrange formalism~\cite{GavrilyukTeshukov01,Gavrilyuk11,GavrilyukNkongaShyueEtAl20}. In addition, in the irrotational framework they can be obtained through canonical Hamilton's equations~\cite{Salmon88,DucheneIsrawi18}, consistently with Zakharov's Hamiltonian formulation of the water waves problem~\cite{Zakharov68}. Shortly put, the SGN system enjoys strong structural properties.

In this work, we numerically compare the SGN equations with a model introduced by the first author and collaborators in~\cite{DucheneIsrawiTalhouk16}. The model is obtained using Hamilton's equations with a modified Hamiltonian, and hence preserves at least part of the structure of the SGN equations, while having the additional property that the dispersion relation of the system coincides exactly with the one of the water waves problem. Models with such properties are often said to be {\em fully dispersive}, and have been advocated by G.~B.~Whitham~\cite{Whitham67} as a way to reproduce ---at least qualitatively--- in a better way some key properties of the water waves problem, such as wavebreaking and non-smooth travelling waves of extreme height. The price to pay is that the equations include non-local pseudodifferential operators (Fourier multipliers). Whitham's prediction turned out to be valid at least for the unidirectional model which bears his name, as shown by~\cite{Hur17,EhrnstromWahlen19,TruongWahlenWheeler20,SautWang20}. This fact triggered renewed activity on bidirectional models (systems), and we refer to the surveys~\cite{KleinLinaresPilodEtAl18,Carter18,DinvayDutykhKalisch19} for more information. The aforementioned model refines systems studied therein (sometimes called Whitham-Boussinesq systems) so as to offer  improved precision in strongly nonlinear situations. We refer to it in this work as the Whitham-Green-Naghdi (WGN) model. It has been rigorously justified among other fully dispersive models in~\cite{Emerald,MM4WW}.\footnote{More precisely, it is proved to be an approximate model to the water waves system of order $\mathcal{O}(\epsilon\delta^4)$ in the sense of consistency, where $\epsilon$ is the ``nonlinearity'' dimensionless parameter defined as the ratio of the typical amplitude of the wave to the reference layer depth, and $\delta$ is the ``shallowness'' dimensionless parameter defined as the ratio of the reference layer depth to the typical horizontal wavelength. The corresponding precision of the SGN equations is $\mathcal{O}(\delta^4)$, and the one of Whitham-Boussinesq equations is $\mathcal{O}(\epsilon\delta^2)$. The improvement between the WGN predictions and SGN predictions can be witnessed in~\cite[Figures~3 and~4]{DucheneNilssonWahlen18} in the context of small-amplitude solitary waves and in~\cite[\S I.5]{MM4WW} for time-evolving profiles. }

In this work, we numerically investigate properties of the SGN and WGN equations in extreme situations. More precisely, we will investigate features and stability of solitary waves with large height and large velocity, as well as solutions whose evolution produces steep gradients.
It must be emphasized that both the SGN and WGN systems are expected to provide poor approximations to the water waves problem in the above scenarii since we voluntarily depart from the shallow-water regime of validity ($\delta\ll 1$). Our motivation is theoretical as we aim at extracting information on the role of dispersive properties for such fully nonlinear models.
 We choose the SGN and WGN equations as our subject of investigations in order to step out of the world of unidirectional scalar (nonlinear and dispersive)  equations for which similar studies have been realized~\cite{AGK,GK12,GK08,KS,etna}, while retaining strong structural properties. In particular, solitary waves can be identified with critical points of functionals which directly derive from the aforementioned Hamiltonian structure~\cite{DucheneNilssonWahlen18}. Moreover, the two systems of equations can (and will) be numerically approximated using identical numerical strategies, specifically Fourier pseudospectral methods.

\subsection{The equations}

In the (horizontal) one-dimensional setting and flat-bottom framework, the SGN equations read (see {\em e.g.}~\cite{DucheneIsrawi18})
   \begin{equation}\label{SGN}\tag{SGN}
      \left\{ \begin{array}{l}
      \partial_t\zeta+\partial_x(hu) =0,\\ \\
   \partial_t\big(u-\frac{1}{3h}\partial_x(h^3\partial_x u)\big)+g\partial_x\zeta+u\partial_x u=\partial_x \big(\frac{u}{3h}\partial_x(h^3\partial_x u)+\frac12 h^2(\partial_x u)^2\big).
      \end{array}\right.
      \end{equation}
 Here, $d>0$ is the reference layer depth, $g>0$ is the gravitation 
 constant,\footnote{By scaling arguments (specifically setting $\zeta(x,t)=d \tilde\zeta(\frac1d x  ,\frac{\sqrt{gd}}{d} t)$, $u(x,t)=\sqrt{gd} \tilde u(\frac1d x  ,\frac{\sqrt{gd}}{d} t)$) it is always possible to 
 set $g=d=1$, and we shall do so in the following.} and $u(t,x)$ represents the layer-averaged horizontal velocity, $\zeta(t,x)$ (or rather its graph) represents the surface deformation and $h(t,x)=d+\zeta(t,x)$ represents the water depth at time $t$ and horizontal position $x\in\RR$. 
 We refer to~\cite{DucheneIsrawi18} for a description of its canonical Hamiltonian structure. Known conserved quantities of~\eqref{SGN} are $\int_{\RR} f_i\,\dd x$ ($i=1,\dots,4$) with densities
\begin{equation}
f_1=\zeta, \quad f_2=hu , \quad f_3= g\zeta^2+hu^2+\frac13 h^3(\partial_x u)^2, \quad f_4= u-\frac{1}{3h}\partial_x(h^3\partial_x u) 
    \label{SGNcons}
\end{equation}
representing respectively the mass, momentum (or horizontal impulse), total energy and the rescaled tangential fluid velocity at the free interface~\cite{GavrilyukKalischKhorsand15}.

The fully dispersive model introduced in~\cite{DucheneIsrawiTalhouk16} (when restricted to the one-layer case and neglecting surface tension) is
\begin{equation}\label{WGN}\tag{WGN}
      \left\{ \begin{array}{l}
      \partial_t\zeta+\partial_x(hu) =0,\\ \\
   \partial_t\big(u-\frac{1}{3h}\partial_x\F(h^3\partial_x\F u)\big)+g\partial_x\zeta+u\partial_x u=\partial_x \big(\frac{u}{3h}\partial_x\F(h^3\partial_x\F u)+\frac12 h^2(\partial_x\F u)^2\big),
      \end{array}\right.
   \end{equation}
   where $\F$ is the Fourier multiplier defined by
   \[\forall \varphi\in L^2(\RR) , \quad    \widehat{\F \varphi}(\xi)=F(d|\xi|)\widehat{\varphi}(\xi) \qquad \text{ where }  F(k)=\sqrt{\frac{3}{|k|\tanh(|k|)}-\frac3{ |k|^2}}.\]
Known conserved quantities of~\eqref{WGN} are $\int_{\RR} f_i\,\dd x$ ($i=1,\dots,4$) with densities
\begin{equation}
f_1=\zeta, \quad f_2= hu , \quad f_3= g\zeta^2+hu^2+\frac13 h^3(\partial_x\F u)^2, \quad f_4= u-\frac{1}{3h}\partial_x\F(h^3\partial_x\F u) .
      \label{WGNcons}
\end{equation}

Our convention for the Fourier transform is the one for which the following identities hold for sufficiently regular and localized functions $g$:
\begin{align*}
  \forall k\in\RR, \quad  \widehat{g}(k) & :=\frac{1}{(2\pi)^{1/2}}\int_\RR
    e^{-ikx}\,g(x) \,\dd x,
    \\
 \forall x\in\RR, \quad    g(x) & 
    =\frac{1}{(2\pi)^{1/2}}\int_\RR
    e^{+ikx}\,\widehat{g}(k) \,\dd k.
\end{align*}

    \subsection{Outline}
    
    Let us now present the structure of this manuscript. In section 
	\ref{S.solitary} we numerically construct solitary wave solutions 
	to the \ref{WGN} equations. In section~\ref{S.evolution} we present our 
	numerical approach for the numerical solution of the initial-value problem
	for the \ref{SGN} and \ref{WGN} equations, and test its validity on 
	solitary wave solutions. The stability of solitary 
	waves for both \ref{SGN} and \ref{WGN} equations is numerically investigated in section 
	\ref{S.stability}. In section~\ref{S.DSW} we study the emergence 
	of zones of rapid modulated oscillations within solutions to both equations
	 starting from unidirectional, long, smooth and localized initial data. 
	Based on the relationship between the \ref{SGN} equations and the Camassa-Holm 
	equation, we study in section~\ref{S.CH} the behaviour of solutions to the former 
	with initial data leading to finite-time blow-up for the latter. In section~\ref{S.cavitation},
	we study the possibility of a blow-up for initial data near cavitation, that is vanishing depth. 
	We summarize our findings and add some concluding remarks in section~\ref{S.conc}. 

\section{Solitary waves}\label{S.solitary}
In this section we study solitary wave solutions
 to the (fully dispersive) 
Whitham-Green-Naghdi equations~\eqref{WGN}, that is solutions of the form
\begin{equation}\label{solitary}
\zeta(t,x)=\zeta_c(x-ct) \qquad u(t,x)= u_c(x-ct), \qquad \lim_{|x|\to\infty}|\zeta_c|(x)+| u_c|(x)=0
\end{equation}
  where the constant  $c\in\RR$ is the solitary wave velocity. 
  
  It is well-known that for any supercritical velocity $c>1$ (recall we set $g=d=1$), there exists a smooth solitary wave solution to ~\eqref{SGN} with explicit formula
   \footnote{In this work we consider only smooth solitary waves maintaining 
    positive depth; see e.g.~\cite{JiangBi17}.
   }
  \begin{equation}
      \zeta_c(x) = (c^2-1)\sech^{2}\Big(\tfrac{\sqrt{3}}{2}\sqrt{\tfrac{c^2-1}{c^2}}\,x\Big),\quad 
       u_c(x) = \frac{c\zeta_c(x)}{1+\zeta_c(x)}
      \label{zetaSGN}.
  \end{equation}
The functions $\zeta_c$ and $u_c$ are smooth, even and positive and have a unique non-degenerate maximum at the origin. Such an explicit formula is of course unexpected for the fully dispersive system~\eqref{WGN}. However, the following result has been shown in~\cite{DucheneNilssonWahlen18}:
\begin{proposition}
There exists $(\zeta^{(q)},u^{(q)})_{q>0}$ a one-parameter family of smooth square-integrable functions such that for all $q>0$, 
 $(\zeta_{c_q},u_{c_q}):=(\zeta^{(q)},u^{(q)})$ provides a solitary wave solution to~\eqref{WGN} with velocity $c_q>1$, and
\[ c_q \to 1 \quad \text{ and }  \quad \big\Vert (c_q^2-1)^{-1}\zeta^{(q)}((c_q^2-1)^{-1/2}\cdot)-\sech^2 (\tfrac{\sqrt3}2\cdot ) \big\Vert_{H^1} \to 0 \quad (q\to 0).\]
\end{proposition}
We also refer to~\cite{DucheneNilssonWahlen18} for numerical computations of \ref{WGN} solitary waves with small velocities $0<c-1\ll 1$. In the following we numerically investigate the existence and behavior of \ref{WGN} solitary waves for large velocities. Based on these numerical experiments we conjecture the following.
\begin{conjecture}\label{conj1}
For all $c>1$, there exist smooth and rapidly decaying solitary wave solutions to the Whitham-Green-Naghdi system~\eqref{WGN} with velocity $c$ and such that the following holds.
\begin{enumerate}
\item For all $c>1$, the profiles $\zeta_c,u_c$ are positive on the real line.
\item For all $c>1$, the profiles $\zeta_c,u_c$ have a unique critical point corresponding to their maximum, and it is non-degenerate.
\item For all $c>1$, the profiles $\zeta_c,u_c$ are symmetric about their maximum.
\item $\max \zeta_c\sim c^2$ and $\max u_c\sim c$ as $c\to\infty$.
\end{enumerate}
\end{conjecture}
This is in sharp contrast to the celebrated result~\cite{AmickToland81} on the existence of (peaked) solitary waves of extreme height for the water waves problem, and the corresponding result obtained on the Whitham equation~\cite{EhrnstromWahlen19,TruongWahlenWheeler20} (see also 
\cite{KS} and references therein for a numerical investigation), and invalidates the naive thinking that this feature relies only on the dispersion relation of the equations linearized about the rest state.

\subsection{The equations for solitary waves}

  Plugging~\eqref{solitary} into~\eqref{WGN}  yields for the first equation 
    \begin{equation}
        \zeta_c=\frac{h_cu_c}{c}=\frac{u_c}{c-u_c},\quad h_c=1+\zeta_c=\frac{c}{c-u_c}
        \label{htravel},
    \end{equation}
    and for the second
    \begin{equation}
        \frac{u_c-c}{3h_c}\partial_{x}\mathsf{F}(h_c^{3}\partial_{x}\mathsf{F} u_c)+\frac{1}{2}
	h_c^{2}(\partial_{x}\mathsf{F}u_c)^{2}+cu_c-\zeta_c-\frac{u_c^{2}}{2}=0.
        \label{utravel}
    \end{equation} 
From~\eqref{htravel} we infer
     \begin{equation}
           u_c=\frac{c\zeta_c}{h_c}
            \label{hurelation},
        \end{equation}
and plugging~\eqref{hurelation} into~\eqref{htravel} yields a single equation for $\zeta_c$, namely
        \begin{equation}
            \frac{-1}{3h_c^{2}}\partial_{x}\mathsf{F}(h_c^{3}\partial_{x}\mathsf{F} \tfrac{\zeta_c}{h_c})+\frac{1}{2}
    	h_c^{2}(\partial_{x}\mathsf{F}\tfrac{\zeta_c}{h_c})^{2}+\frac{\zeta_c}{h_c}-\frac{\zeta_c}{c^2}-\frac{\zeta_c^2}{2h_c^2}=0.
            \label{htravel-2}
        \end{equation}
        Similarly, we may use~\eqref{hurelation} with~\eqref{utravel} to produce a single equation for $u_c$, namely
                 \begin{equation}
                        -\frac{(c-u_c)^2}{3c}\partial_{x}\mathsf{F}\big((\tfrac{c}{c-u_c})^{3}\partial_{x}\mathsf{F} u_c\big)+\tfrac{1}{2}
                	(\tfrac{c}{c-u_c})^{2}(\partial_{x}\mathsf{F}u_c)^{2}+cu_c-\frac{u_c}{c-u_c}-\tfrac12 u_c^2=0.
                        \label{utravel-2}
                    \end{equation}  
Finally, we find it convenient to solve~\eqref{utravel-2} using the variable $\eta_c=\frac{u_c}{c}=\frac{\zeta_c}{h_c}$:
  \begin{equation}
             \frac{-(1-\eta_c)^2}{3}\partial_{x}\mathsf{F}\Big(\frac{\partial_{x}\mathsf{F}\eta_c}{(1-\eta_c)^3}\Big)+\frac{1}{2(1-\eta_c)^2}
     	\big(\partial_{x}\mathsf{F}\eta_c\big)^{2}+\eta_c-\frac{\eta_c}{c^2(1-\eta_c)}-\tfrac{1}{2}\eta_c^2=0.
             \label{utravel-3}
         \end{equation}         
\begin{remark}\label{R.Lagrange}
Equation~\eqref{htravel-2}  can be written as $\delta_{\zeta_c} \mathcal L=0$  with
\[\mathcal L:= \int_{\RR} \ell (\zeta_c,c\tfrac{\zeta_c}{1+\zeta_c}) \dd x\]
where
\[\ell(\zeta,u) := \frac12 \zeta^2 - \frac12 (1+\zeta) u^2 - \frac16 (1+\zeta)^3 (\partial_{x}\mathsf{F} u)^2\]
is the Lagrangian density naturally associated with the Hamiltonian formulation of the equations, and physically corresponds to the difference between the potential energy and the kinetic energy. In particular, the Jacobian of $\delta\mathcal L$ is, by definition, symmetric for the $L^2$-inner product.
\end{remark}

\subsection{Numerical construction of solitary waves}

We seek numerical approximations of solitary waves for~\eqref{WGN} with fixed velocity $c>1$
 through zeroes of a finite-dimensional vector-field accounting for the left-hand side of~\eqref{utravel-2} or~\eqref{utravel-3}:
\begin{equation}\label{discrete}\mathcal F(\u) = 0.\end{equation}
The vector $\u=(u(x_{1}),\ldots, u(x_{n}))$ represents values of the function $u$ at collocation points 
$x_{n}=-\pi L+ n2\pi L/N$, $n=1,\ldots,N$ for $x\in L[-\pi,\pi]$. Nonlinear operations are naturally computed at collocation points,
while $\partial_{x}\mathsf{F}$ is approximated via a discrete Fourier 
transform computed efficiently with a \emph{Fast Fourier 
transform} (FFT) and multiplication in Fourier space, that is
 \[ \forall \v \in \RR^N, \qquad \widehat{\partial_{x}\mathsf{F} \v}_k := {\rm i} \k_k \, F(\k_k)\, \widehat \v_k\]
 where we denote $\widehat\v=(\widehat \v_{-N/2+1},\ldots,\widehat 
 \v_{N/2})$ the coefficients of the Fast Fourier transform of $\v$ 
 (which we slightly incorrectly refer to as {\em Fourier 
 coefficients} in the following), and $\k_k=(k/L)$, $k=-N/2+1,\ldots,N/2$ the discrete Fourier modes.
Here $L>0$ is a constant chosen such that $u$ and its 
relevant derivatives decrease to machine precision (roughly 
$10^{-16}$ in double precision) and $N$ will be chosen such that the $\widehat \u_k$ decreases to machine precision for 
large  $|k|$.
 In the results presented here, multiplying $N$ or $LN$ 
by a factor $2$ or $4$ neither improves nor deteriorates the 
accuracy. For more information on Fourier spectral 
methods, the reader is referred to~\cite{trefethen} and the 
literature cited therein. 

The system of $N$ nonlinear equations~\eqref{discrete} will 
be solved iteratively by a standard Newton iteration,
\begin{equation}
    \u^{(m+1)} = \u^{(m)} - \delta\u^{(m)}, \label{newton}
    \end{equation}
with
\begin{equation}
    \mbox{Jac}(\mathcal{F}(\u))|_{\u=\u^{(m)}}\delta\u^{(m)} = \mathcal{F}(\u^{(m)})
    \label{Jac},
\end{equation}
where $\u^{(m)}$ denotes the $m$th iterate, and where 
$\mbox{Jac}\mathcal{F}(\u)$ is the Jacobian of $\mathcal{F}(\u)$ with respect 
to $\u$. 
The initial iterate $\u^{(0)}$ will be chosen as the \ref{SGN} solitary wave given by~\eqref{zetaSGN} at collocation points. It may also be constructed by extrapolating from previously computed solitary waves and Lagrange polynomials.
For the iteration, we apply and compare two strategies: first we test a Newton-GMRES method, 
i.e., we solve~\eqref{Jac} with the Krylov subspace iterative method GMRES~\cite{GMRES} as for instance in~\cite{KS}. 
Alternatively we solve~\eqref{Jac} through standard LU factorization. Notice that due to the translation invariance of the problem,
the Kernel of the Jacobian of the continuous (infinite-dimensional) vector-field is non-empty when evaluated at a non-trivial solution $u$, since $\partial_x u$ is an element of its nullspace.  The corresponding spectral projection can be inferred from the symmetry property mentioned in Remark~\ref{R.Lagrange}. Correspondingly, $\mbox{Jac}(\mathcal{F}(\u^{(m)}))$ has an extremely small eigenvalue as $\u^{(m)}$ converges towards the desired solution. However, by symmetry considerations, we can ensure that at each iterate $\u^{(m)}$ and therefore $\mathcal{F}(\u^{(m)})$ are even, and as a consequence its spectral projection onto the corresponding eigenspace vanishes (up to machine precision).    Yet we find it advisable to add the aforementioned spectral projection to the Jacobian when solving~\eqref{Jac}, although in practice we mostly observe a slight phase shift on the numerical approximation if the spectral projection is not added.

\begin{figure}[htb!]
  \includegraphics[width=0.49\textwidth]{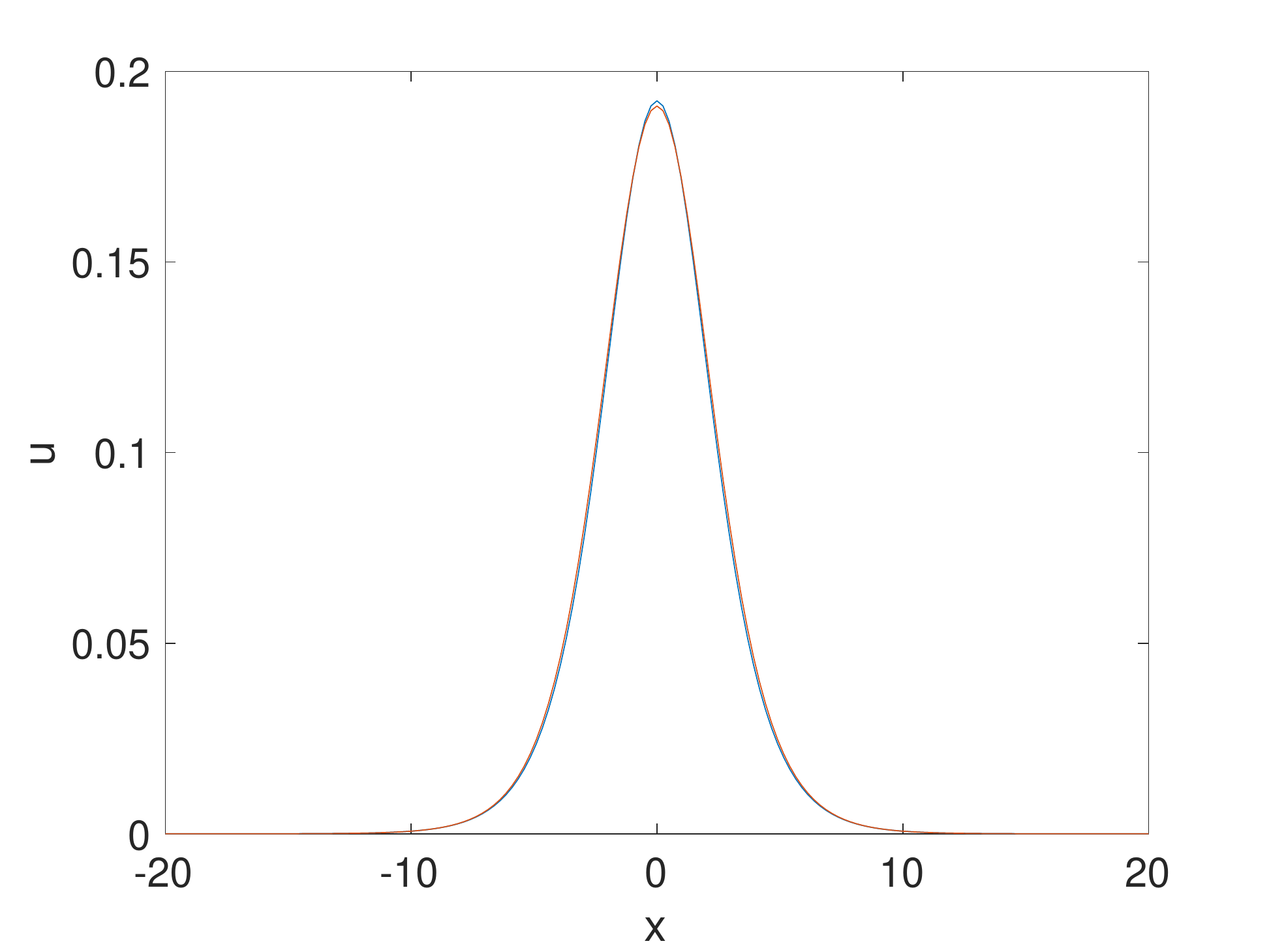}
  \includegraphics[width=0.49\textwidth]{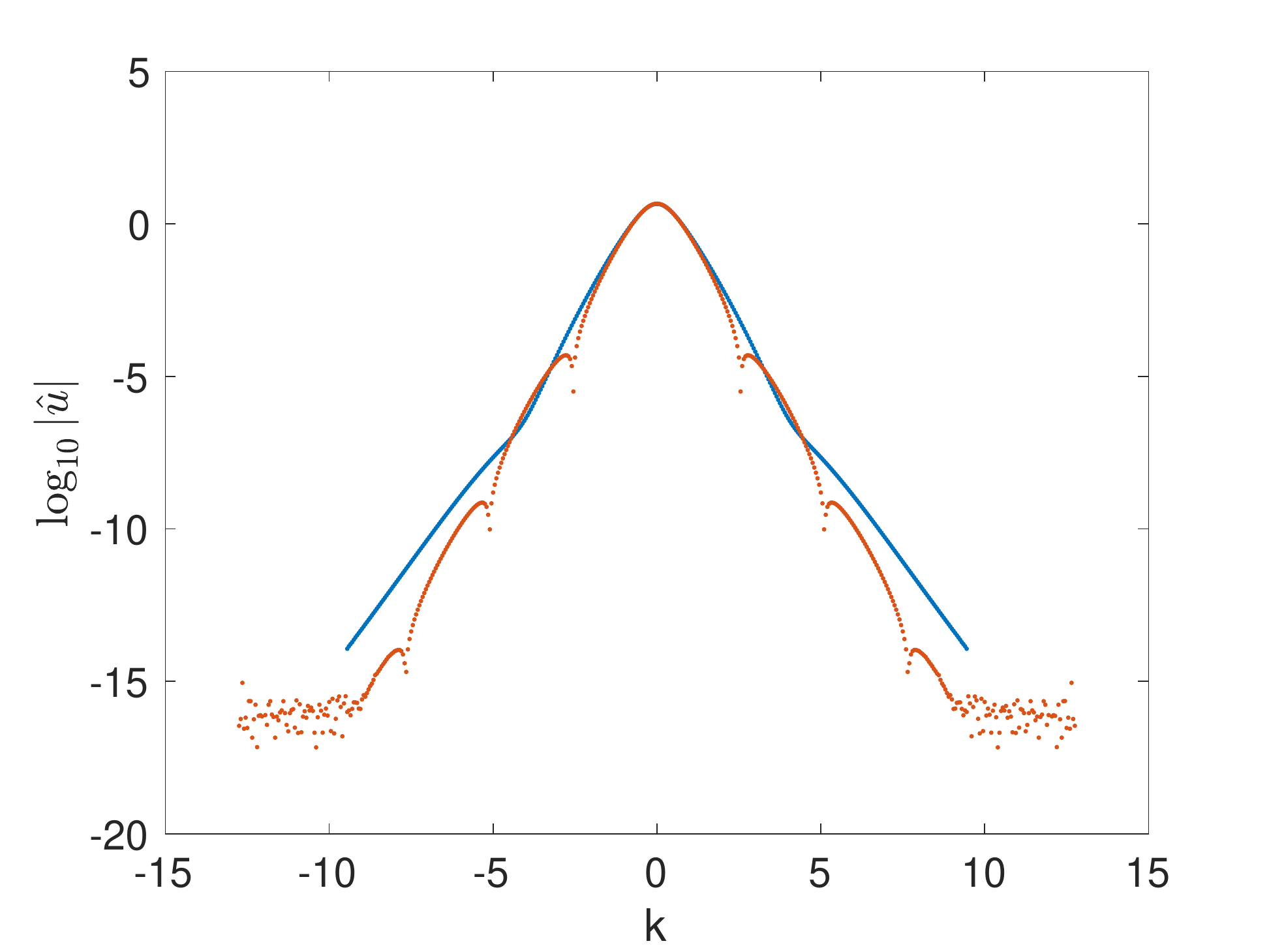}
 \caption{Left: solitary wave for the \ref{WGN} equations for $c=1.1$ in blue 
 and the \ref{SGN} equations for the same velocity in red; right: Fourier 
 coefficients for both solitary waves on the left.}
 \label{fdSGNsolc11}
\end{figure}

Let us now present our results, starting with the case of small and slow solitary waves. For $c=1.1$, 
we treat the equation~\eqref{utravel-2} for $x\in 20[-\pi,\pi]$ and use 
$N=2^{9}$ collocation points. We use a Krasny filter 
of the order of $10^{-14}$, which puts to zero Fourier coefficients with modulus smaller than $10^{-14}$. 
We also apply a preconditioner of the 
form $M = {\rm Diag}((1+\k^{2}/3))$ ---which is motivated by the linear dispersion of 
the \ref{SGN} equations--- i.e., instead of solving iteratively with GMRES  $A{\bf x}={\bf b}$ the 
equation~\eqref{Jac}  in  Fourier space, we solve  $M^{-1}A{\bf 
x}=M^{-1}{\bf b}$.
The Newton-GMRES code converges within 3 iterations with a residual 
$\big\Vert\mathcal{F}(\u^{(3)})\big\Vert_{\ell^\infty}$ of the order of $10^{-13}$. The residual 
of the initial iterate, $\big\Vert\mathcal{F}(\u^{(0)})\big\Vert_{\ell^\infty}$ with $\u^{(0)}$ the \ref{SGN} solitary wave, is of the order of 
$10^{-2}$, GMRES converges in 21 iterations with a relative residual 
of $10^{-11}$. The resulting Newton iterate, $\big\Vert\mathcal{F}(\u^{(1)})\big\Vert_{\ell^\infty}$ has a residual of $10^{-4}$, GMRES again 
converges within 20 iterations to a relative residual of 
$10^{-11}$. The residual of the subsequent  Newton iteration, $\big\Vert\mathcal{F}(\u^{(2)})\big\Vert_{\ell^\infty}$ is then of the order of 
$10^{-7}$. In the next step, GMRES stagnates with a 
relative residual of the order of $10^{-7}$, which is explained by the smallness
of the previous residual of the Newton iteration. The ensuing residual of the 
Newton iteration, $\big\Vert\mathcal{F}(\u^{(3)})\big\Vert_{\ell^\infty}$, is of the order of $10^{-13}$, and 
the iteration is stopped. The iteration thus shows 
the well-known quadratic convergence of a Newton iteration, loosely 
speaking the number of correct digits doubles in each 
iteration.
The resulting solitary wave can be seen in 
Fig.~\ref{fdSGNsolc11} on the left. The \ref{WGN} solitary wave for $c=1.1$ is 
very close to the \ref{SGN} solitary wave for the same velocity shown in red in 
the same figure. On the right of the same figure, the modulus of the 
Fourier coefficients shows that the solution is numerically resolved 
to the order of the Krasny filter. We also observe that the exponential 
decay rate of Fourier coefficients with large Fourier modes is slightly
smaller for the \ref{WGN} solitary wave than it is for the \ref{SGN} solitary wave. 

\begin{figure}[htb!]
  \includegraphics[width=0.49\textwidth]{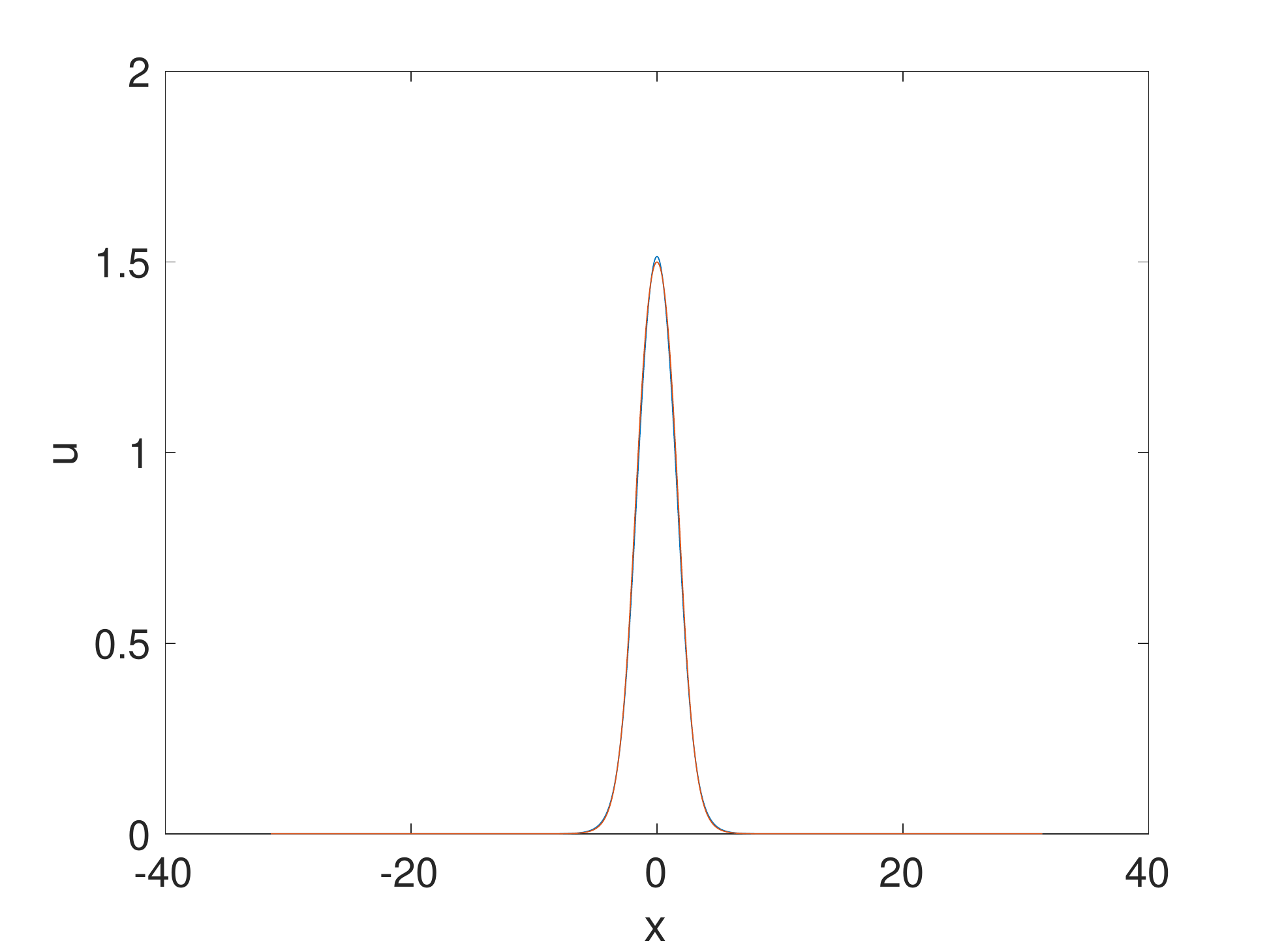}
  \includegraphics[width=0.49\textwidth]{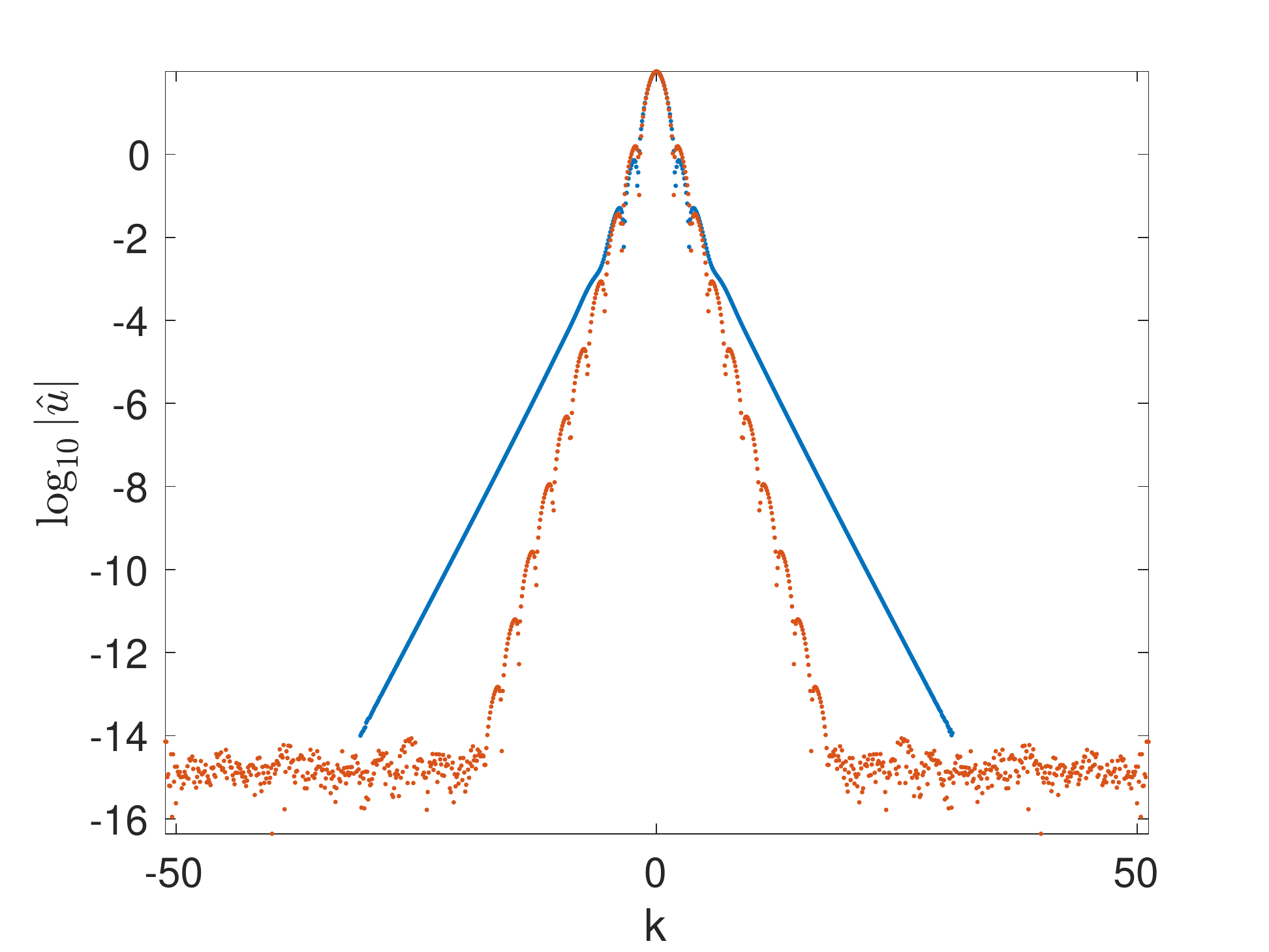}
 \caption{Left: solitary wave for the \ref{WGN} equations for $c=2$ in blue 
 and the \ref{SGN} equations for the same velocity in red; right: Fourier 
 coefficients for the  solitary waves on the left.}
 \label{fdSGNsolc2}
\end{figure}

Now we consider a larger value of the velocity, $c=2$ and use 
$N=2^{10}$ collocation points for $x\in 20[-\pi,\pi]$. The 
Newton-GMRES code converges in 4 iterations with a residual of the 
order of $10^{-12}$. The resulting solitary wave is shown in 
Fig.~\ref{fdSGNsolc2} on the left. Again the solution is very close 
to the \ref{SGN} solitary wave shown in the same figure in red. The solution is 
well resolved in Fourier space as can be seen on the right of 
Fig.~\ref{fdSGNsolc2}, and  we observe once more that the exponential 
decay rate of Fourier modes is smaller. The Newton-GMRES iteration behaves similarly 
to what is described before. Note that it does not converge without a
 preconditioner due to issues for the high Fourier modes.  
 For even larger 
 values of the velocity such as $c=3$, the iteration no longer converges 
 because of GMRES problems for the high Fourier modes. 

\begin{figure}[htb!]
  \includegraphics[width=0.49\textwidth]{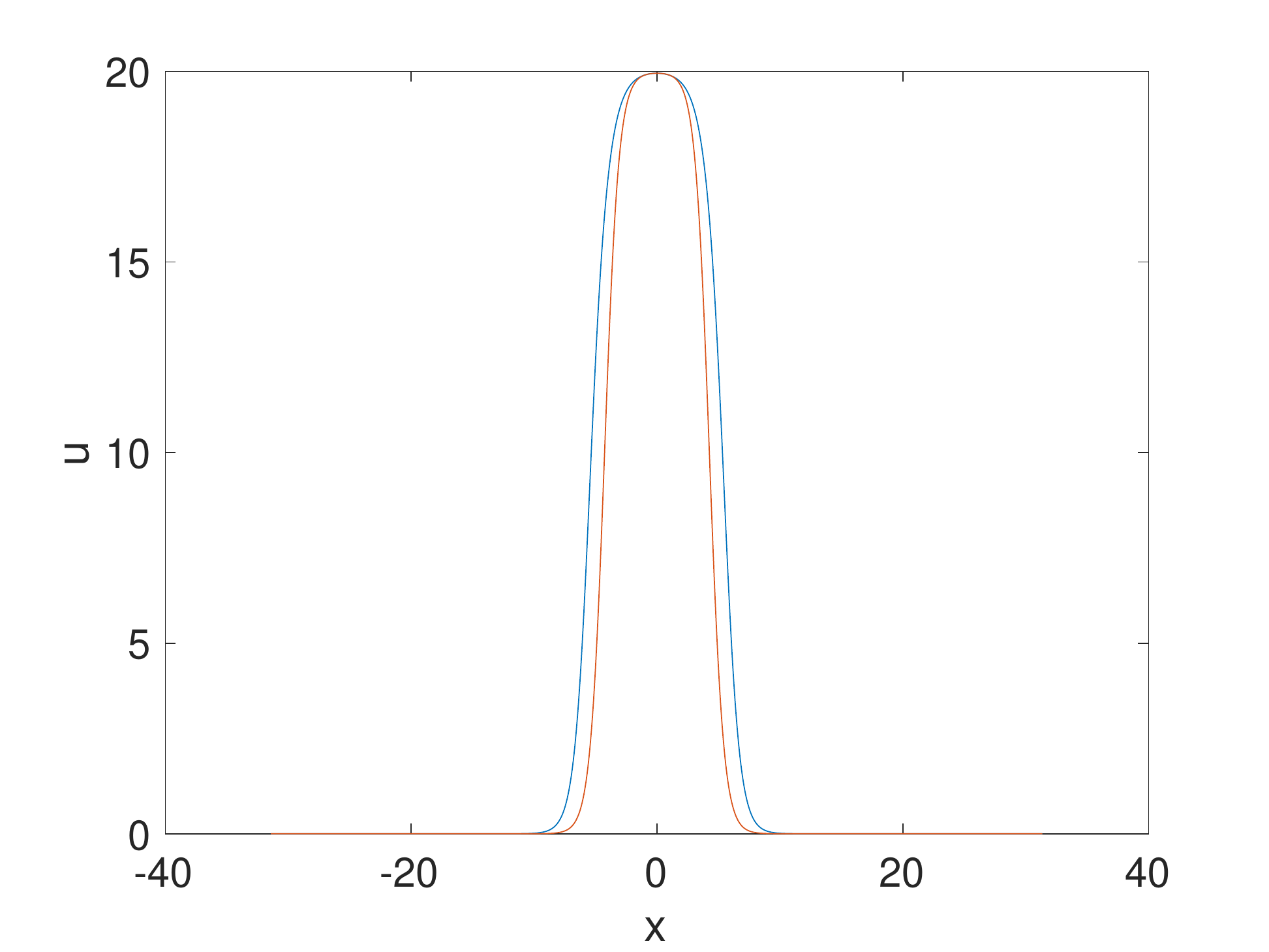}
  \includegraphics[width=0.49\textwidth]{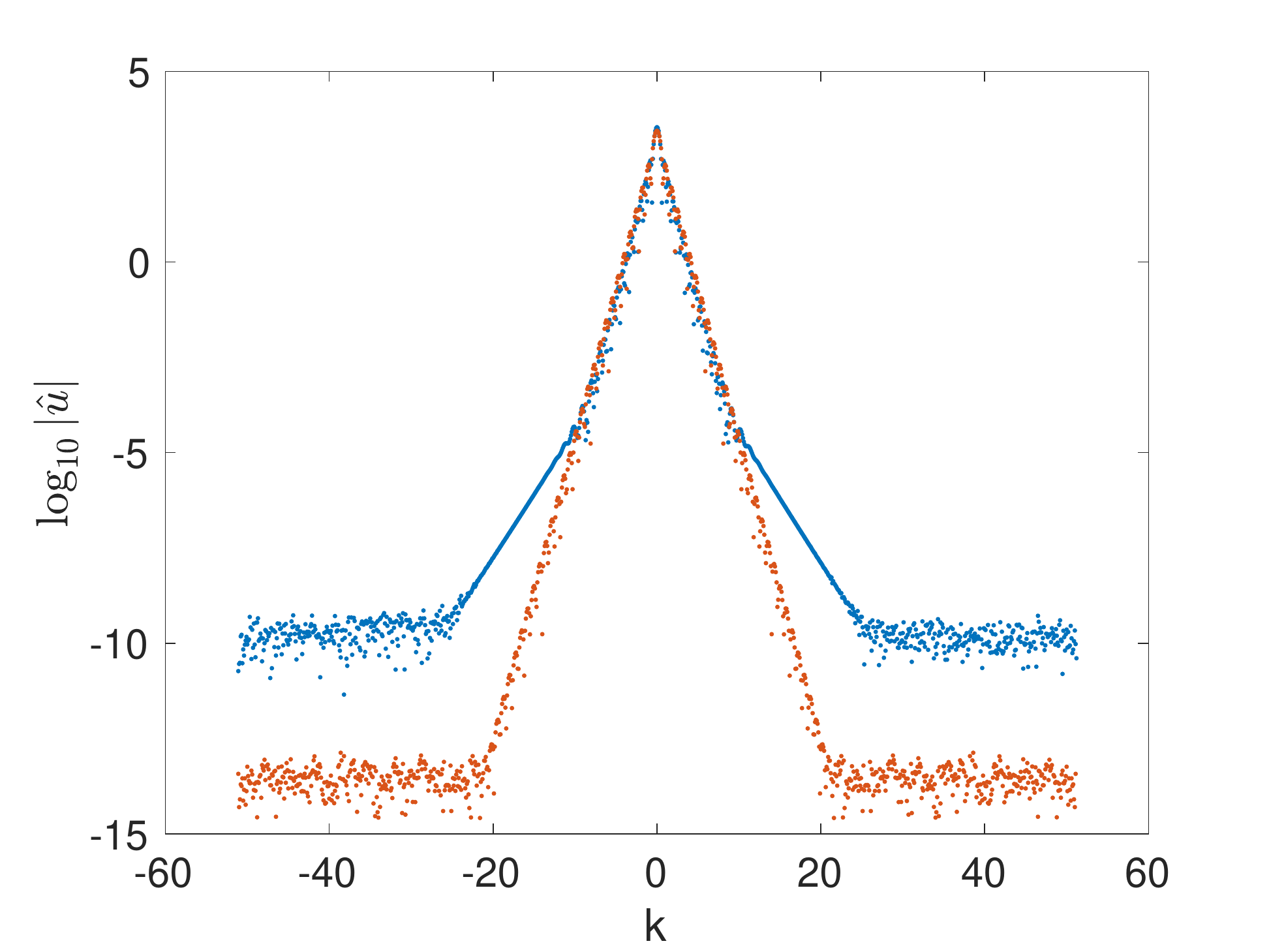}
 \caption{Left: solitary wave for the \ref{WGN} equations for $c=20$ in blue 
 and the \ref{SGN} equations for the same velocity in red; right: Fourier 
 coefficients  on the left.}
 \label{fdSGNsolc20}
\end{figure}

Therefore we switch for larger values of the velocity $c$ to a Newton 
iteration with a direct numerical factorization of the Jacobian. Still with 
$N=2^{10}$ collocation points on $x\in 10[-\pi,\pi]$, we
consider the case $c=3$. The Newton iteration converges with a 
direct inversion of the Jacobian in 3 iterations to a residual of the 
order of $10^{-9}$. After further iterations, the minimal residual reachable with this method appears 
to be of the order of $10^{-13}$. For $c=20$ and the same parameters 
as before, iteration converges normally to the order of $10^{-7}$.
The solution and its Fourier coefficients can 
be seen in Fig.~\ref{fdSGNsolc20}. The deviation from the \ref{SGN} solitary wave 
in red is now clearly visible. We again observe a lower rate of decay of 
Fourier coefficients of the \ref{WGN} solitary wave, however only for higher Fourier modes (the lower decay rate is apparent starting from much smaller Fourier modes when plotting the corresponding figures for the surface deformation, $\zeta_c$).
The Fourier coefficients saturate at 
the order of $10^{-10}$ indicating problems in the conditioning of 
the Jacobian which we investigate later on. Part of this can be attributed 
to the presumption motivated by formula~\eqref{zetaSGN} that $\zeta_c$ 
---and hence $h_c=1+\zeta_c=\tfrac{c}{c-u_c}$--- scales as $c^{2}$ for large $c$. 

\begin{figure}[htb!]
  \includegraphics[width=0.49\textwidth]{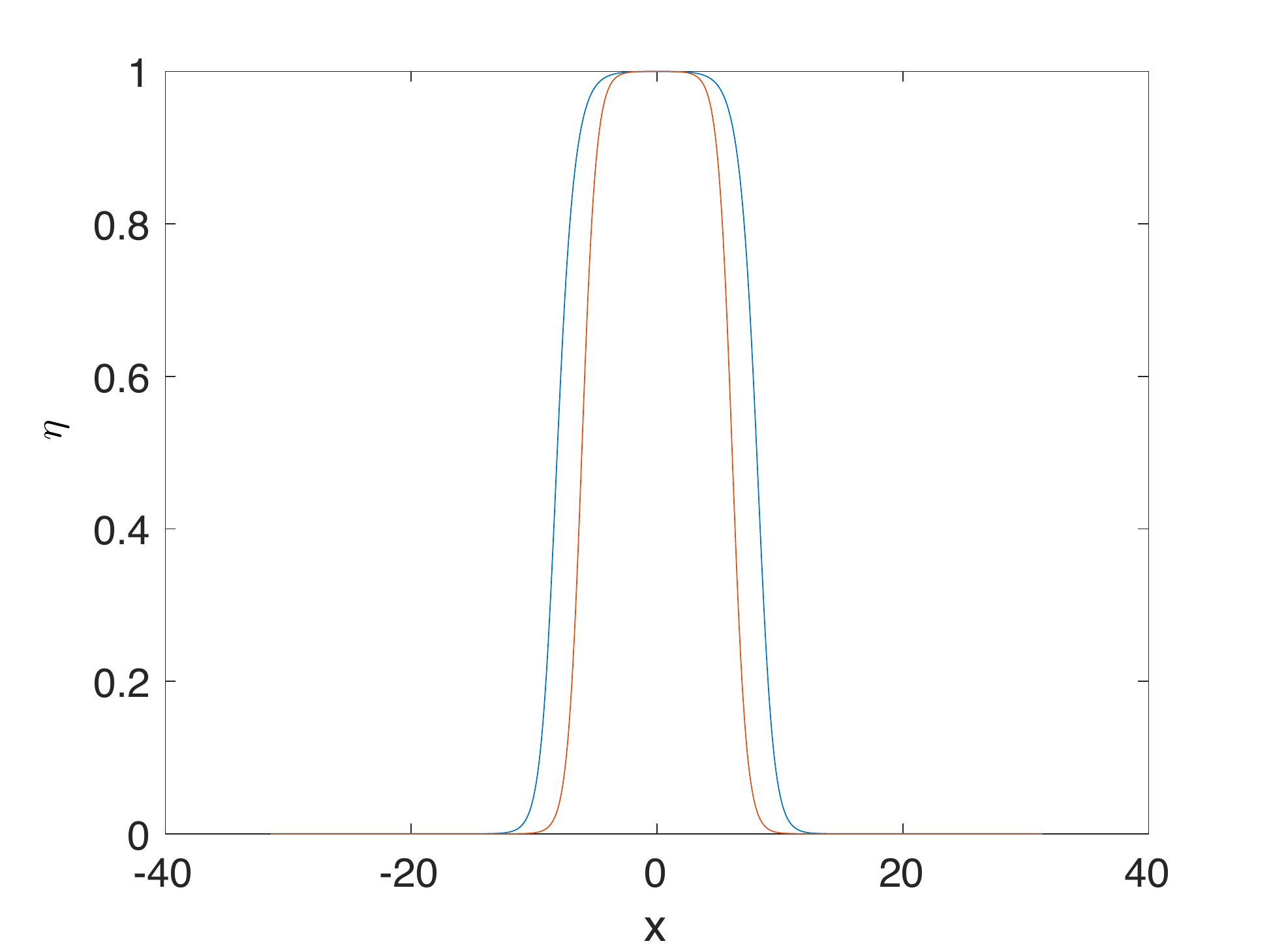}
 \includegraphics[width=0.49\textwidth]{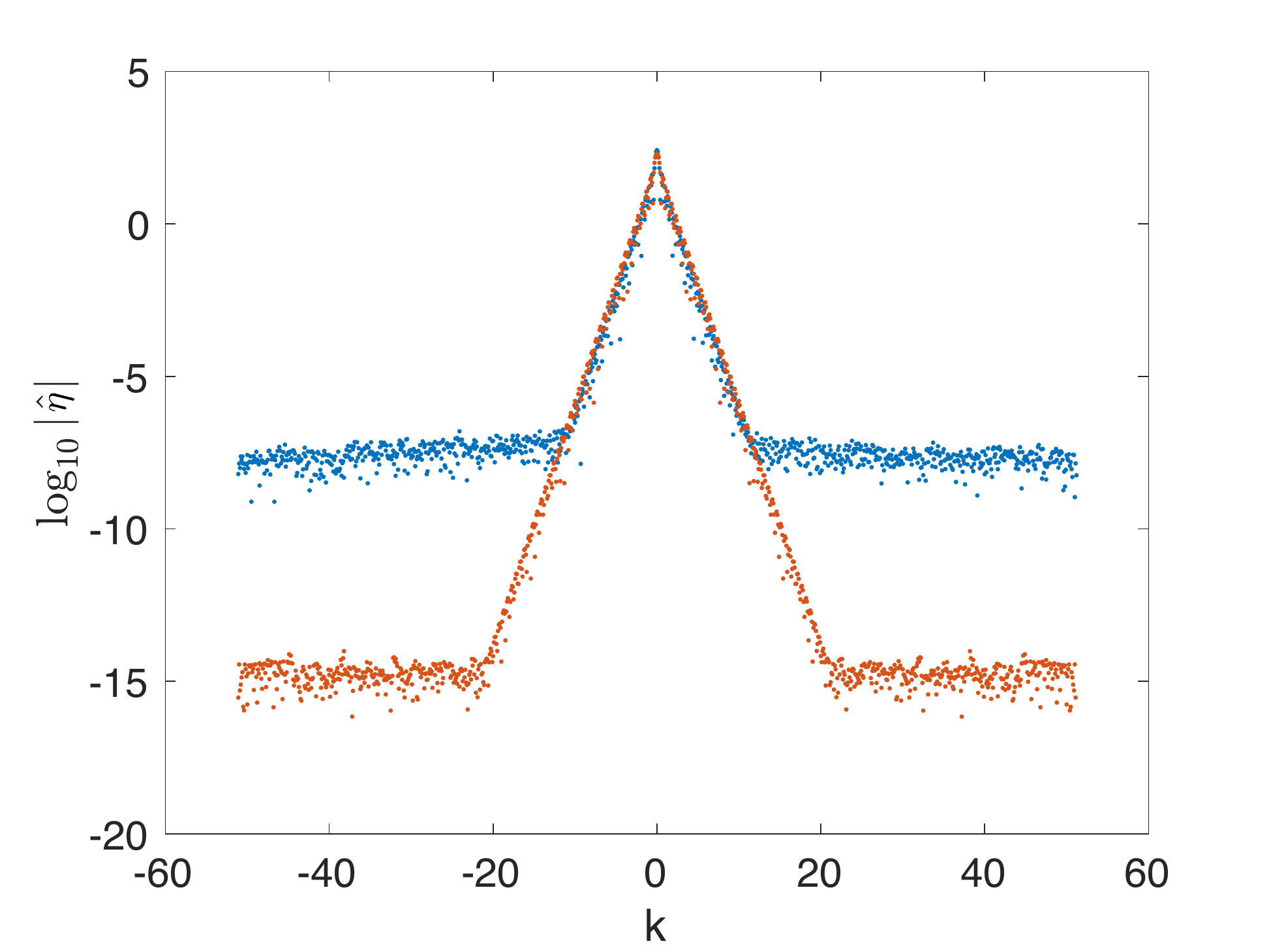}
 \caption{Left: solitary wave for the \ref{WGN} equations for $c=100$ in blue 
 and the \ref{SGN} equations for the same velocity in red; right: Fourier 
 coefficients  on the left.}
 \label{fdSGNsolc100}
\end{figure}

For even larger values of 
the velocity as $c=100$, these problems become worse and the Newton iteration does not converge. 
There is however no indication that the solitary wave might become singular or might not exist for such values. 
Since $u_c$ is proportional to $c$, we look at 
the rescaled quantity $\eta_c=u_c/c=\zeta_c/(1+\zeta_c)$ and consider the equation~\eqref{utravel-3}.
Using again the same parameters 
as before, the optimal residual for the iteration in this case is of the order of 
$10^{-6}$. The 
solution can be seen in Fig.~\ref{fdSGNsolc100} on the left. The 
Fourier coefficients on the right of the same figure saturate at the 
order of $10^{-8}$ which partly explains why no lower residual can 
be achieved. The lower exponential rate of decay of the \ref{WGN} 
solitary wave Fourier coefficients is no longer visible, presumably because it occurs for higher Fourier modes than the numerically well-resolved ones.
 It is apparent when plotting the corresponding figures for the surface deformation, $\zeta_c$.
 
Let us briefly comment on this decay rate. It is well-known that it relates to the maximal width of the strip around the real axis for which the analytic extension of the function is free of singularities (see {\em e.g.}~\cite[Theorem~7.1]{AlazardBurqZuily}). The authors have no clear understanding on how such properties can be related to the balancing of dispersion and nonlinearity of the corresponding equation. The exponential decay rate of the Fourier transform of the surface deformation of the \ref{SGN} solitary waves can be explicitly inferred from formula~\eqref{zetaSGN}. Such is not the case for the velocity variable, to the authors' knowledge, due to the fact that the function $z\mapsto z/(1+z)$ is not entire (see {\em e.g.}~\cite[Proposition~7.10]{AlazardBurqZuily}), although numerical computations indicate that the decay rate is the same. The decay rate for the surface deformation and velocity variables of the \ref{WGN} solitary waves also appear to be the same. They are close to the ones of the \ref{SGN} solitary waves for small values of the velocity $c\approx 1$, and quickly deviate for larger values (about $c=2$) to approximately half the ones of the \ref{SGN} solitary waves. This interesting question will certainly need to be studied more in the future.
 
\begin{figure}[htb!]
  \includegraphics[width=0.49\textwidth]{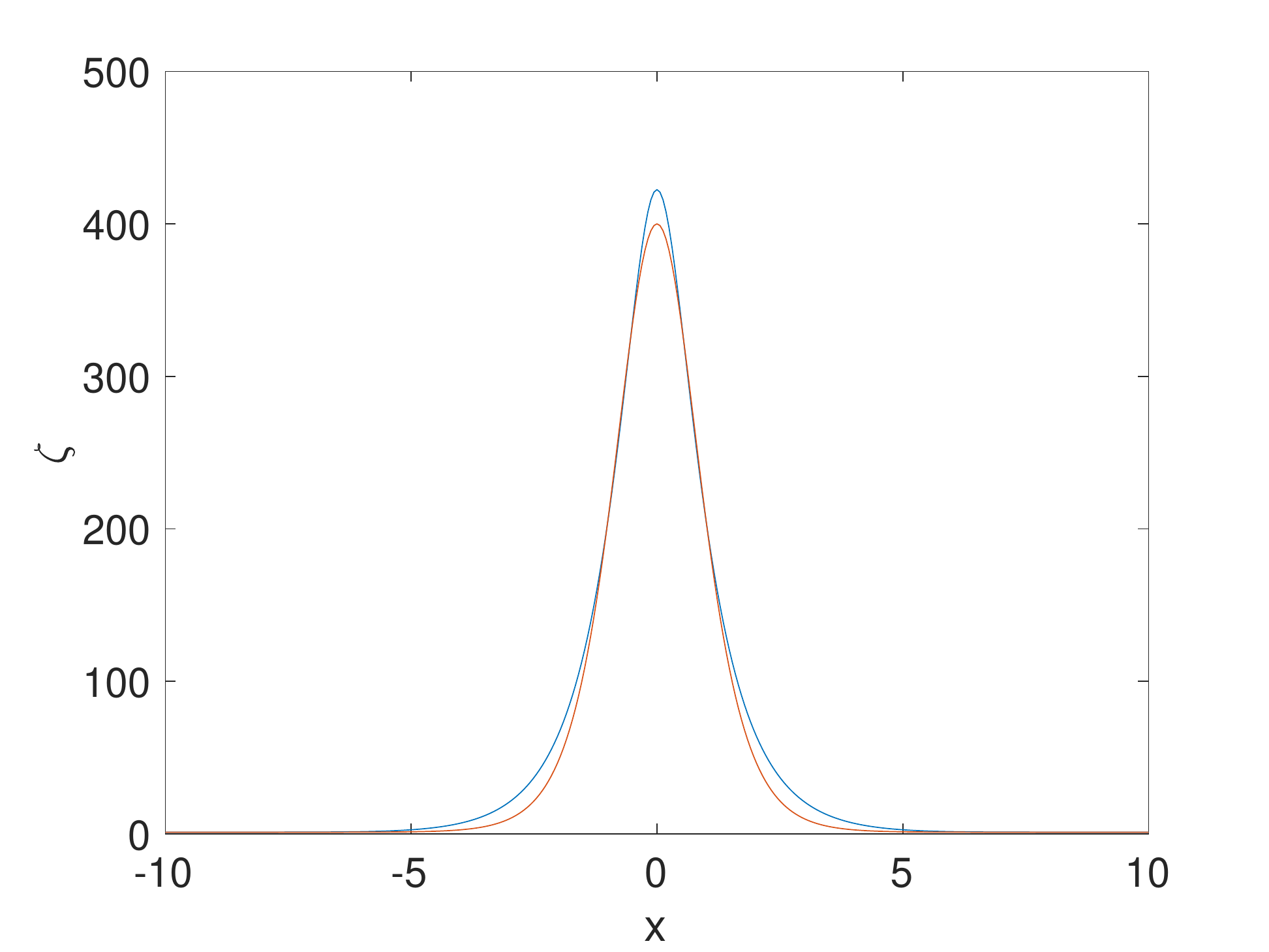}
  \includegraphics[width=0.49\textwidth]{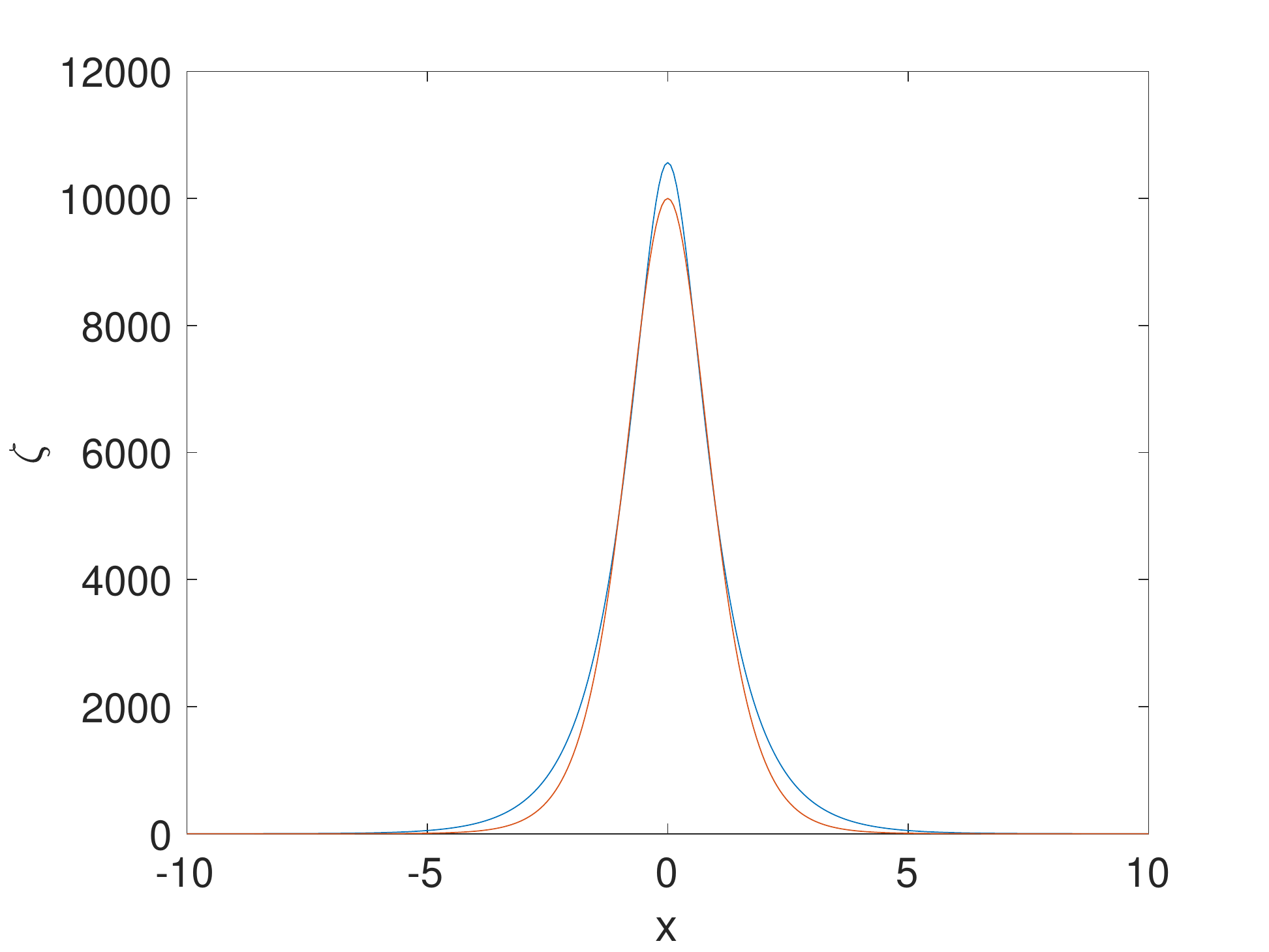}
 \caption{The function $\zeta$ for the solitary wave for the \ref{WGN} equations  in blue 
 and the \ref{SGN} equations for the same velocity in red: left $c=20$, 
 right $c=100$.}
 \label{fdSGNsolzeta}
\end{figure}

One reason for the problems in the iterations are due to the function $\zeta_c$ being 
 of the order $c^{2}$. We show these functions for $c=20$ and 
$c=100$ in Fig.~\ref{fdSGNsolzeta}. It can be seen that this function 
is always more peaked than the corresponding one for the \ref{SGN} equations with the same 
velocity, given by the explicit formula~\eqref{zetaSGN}. 

\begin{figure}[htb!]
  \includegraphics[width=0.49\textwidth]{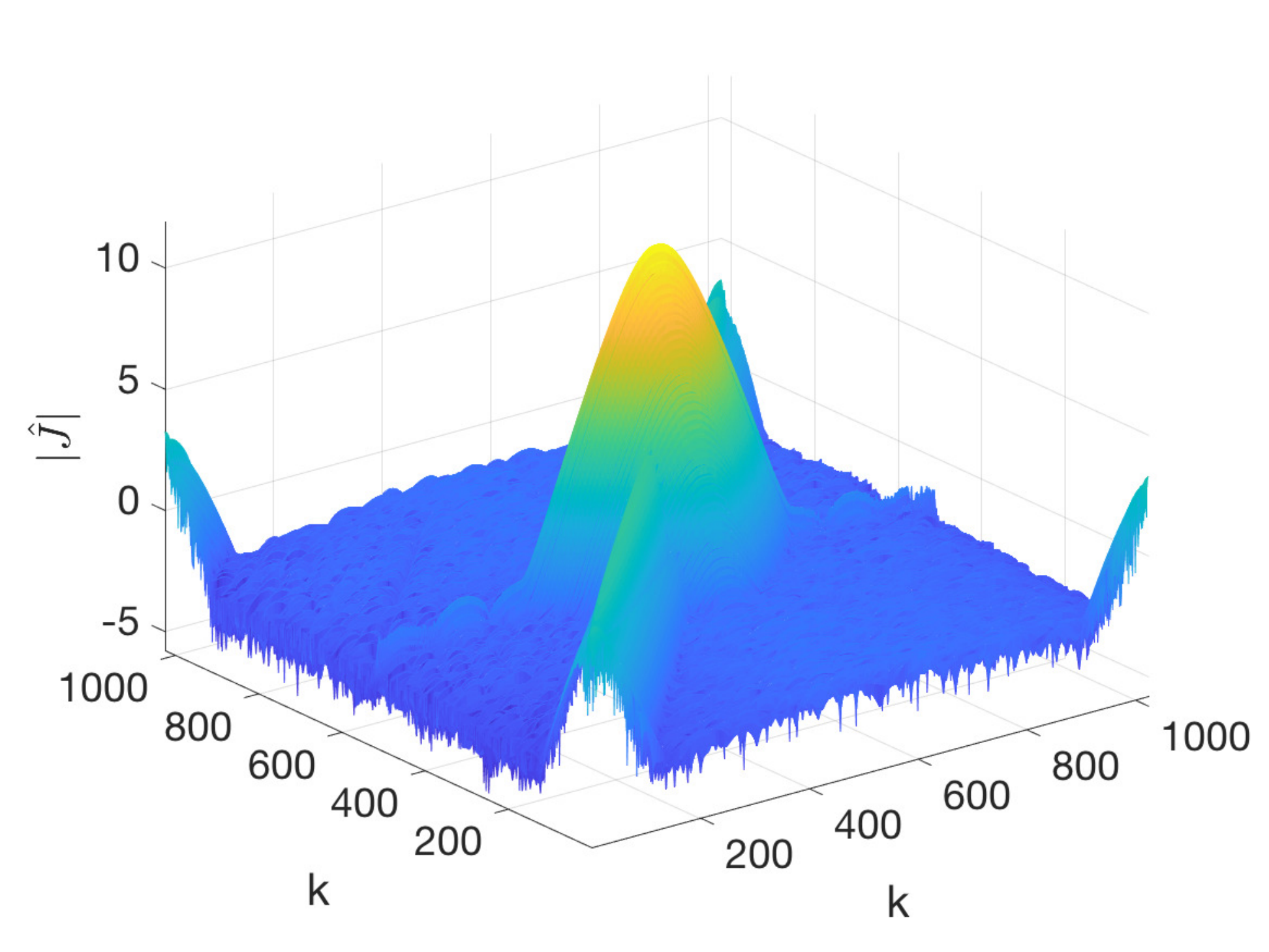}
  \includegraphics[width=0.49\textwidth]{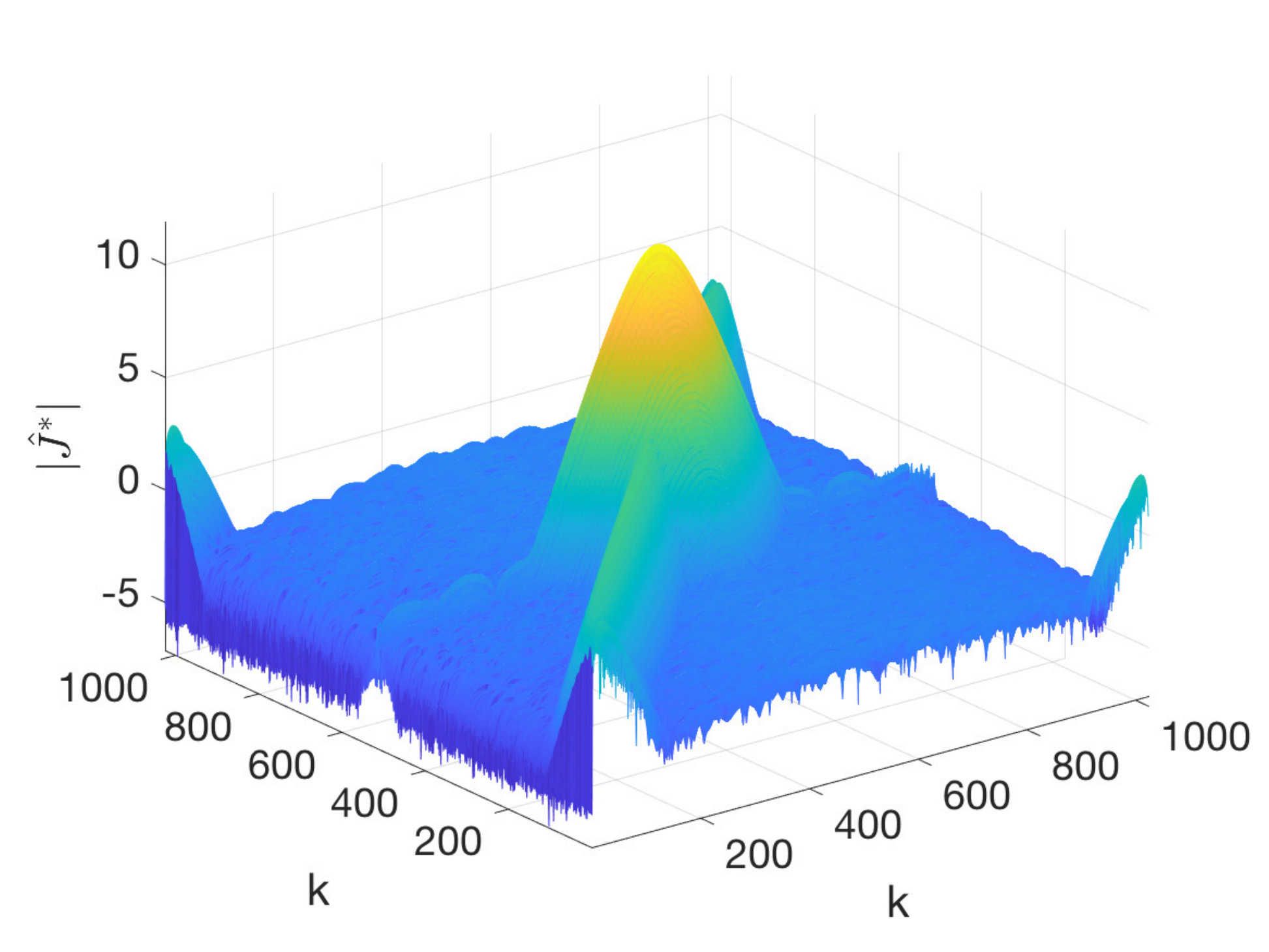}
 \caption{Left: discrete Fourier transform of the Jacobian for the initial 
 iterate (\ref{SGN} solitary wave) for $c=20$; right: non diagonal part $\hat J^{*}$ on the left.}
 \label{fdSGNsolc100jac}
\end{figure}

To understand the difficulties in the Newton iterations, we look at the Jacobian for the initial iterate (the \ref{SGN} solitary wave) for $c=100$ which is shown in Fig.~\ref{fdSGNsolc100jac}. We 
denote the discrete Fourier transform of this Jacobian (from the left 
and from the right) by $\hat{J}$. It can be seen that it is, except for a 
peak near the center, essentially constant along the diagonal. In 
addition there is a plateau of the order of $10^{-5}$ seemingly due 
to rounding errors since the maximum is of the order of $10^{10}$. 
Note that the dominant contribution to the Jacobian is due to 
$\hat J^{*}:= \frac{-(1-\eta_c)^2}{3}\partial_{x}\mathsf{F}\Big(\frac{1}{(1-\eta_c)^3}\partial_{x}\mathsf{F}\Big)$
which we refer to as the `non diagonal' part. 

Thus it appears that the problems in the computation of the solitary waves 
for very high velocities are due to machine precision errors being
increased by the multiplication in physical space by the function
$\frac{1}{(1-\eta_c)^3}$ taking large values (of order $c^6$), and
by the application of $\partial_{x}\mathsf{F}$ which multiplies 
Fourier modes by $i\k_k\, F(\k_k)\sim i\sqrt{3|\k_k|} $ for large $k$.
But there is no 
indication of a maximal velocity for the solitary wave solutions to~\eqref{WGN}. The conclusions of the numerical investigation in this section are summarized in the form of Conjecture~\ref{conj1}.

\section{Time evolution}\label{S.evolution}
    In this section we present and validate the numerical scheme we 
	employ for integrating in time the equations~\eqref{SGN} 
	and~\eqref{WGN}. In principle we adapt here a 
	classic code for the KdV equation (code 27 in~\cite{trefethen}) 
	to the equations studied here. The main new aspect is the 
	inversion of the elliptic operator via GMRES, and we concentrate 
	on this aspect whilst refering to~\cite{trefethen} for general 
	properties of the numerical approach. Let us first recall the 
	SGN and WGN equations with a slight reformulation: sufficiently regular solutions to~\eqref{WGN} satisfy (once again, we set $g=d=1$)
 \begin{equation}\label{SGN-WGN}
      \left\{ \begin{array}{l}
      \partial_t\zeta+\partial_x(hu) =0,\\ \\
   \partial_t v +\partial_x\big(\zeta+ uv-\tfrac12 u^2-\tfrac12 h^2(\partial_x {\sf F}u)^2\big)=0
      \end{array}\right.
      \end{equation}
      where we recall that $h=1+\zeta$ and ${\sf F}$ is the Fourier multiplier operator with symbol $F(k)=\sqrt{\frac{3}{|k|\tanh(|k|)}-\frac3{ |k|^2}}$; and $v$ and $u$ are related through the elliptic equation 
 \begin{equation}\label{elliptic}
     v =  u-\frac{1}{3h}\partial_x\big(h^3\partial_x{\sf F} u\big).
 \end{equation}
 Sufficiently regular solutions to~\eqref{SGN} satisfy the above, replacing ${\sf F}$ by the identity.
  By standard elliptic theory ~\cite[Lemma~5.45]{Lannes}, $u$ is uniquely determined by~\eqref{elliptic} from sufficiently regular $(v,\zeta)$ with $\inf_\RR (1+\zeta)>0$, and we can solve~\eqref{SGN-WGN} as evolution equations for $(\zeta,v)$. Incidentally, notice $v$ represents a rescaled tangential fluid velocity at the free interface, the fourth quantity in~\eqref{SGNcons} and~\eqref{WGNcons}.
    
        \subsection{The numerical scheme}
    
    Many numerical schemes have been proposed for solving the \ref{SGN} equations; see for instance~\cite{CienfuegosBarthelemyBonneton06,LeGavrilyukHank10,BonnetonChazelLannesEtAl11,DutykhClamondMilewskietal13,LiGuyenneLietal14,MitsotakisIlanDutykh14,MitsotakisDutykhCarter17,PittZoppouRoberts18,GavrilyukNkongaShyueEtAl20,AissioueneBristeauGodlewskiEtAl20,DorodnitsynKaptsovMeleshko,AntonopoulosDougalisMitsotakis}.
The presence of Fourier multipliers in the \ref{WGN} equations naturally leads to Fourier pseudospectral methods, already employed in~\cite{DutykhClamondMilewskietal13} and described thereafter.
     One of the difficulties when integrating in time the \ref{SGN} 
	 equations (and in a similar way the \ref{WGN} equations) is that 
	 we are led to solve the elliptic problem~\eqref{elliptic} at 
	 each time-step.  The aforementioned issue is addressed in 
	 particular 
	 in~\cite{DiasMilewski10,LannesMarche15,DuranMarche15} (see also~\cite{FavrieGavrilyuk17} and references in \cite{MM4WW} for relaxation approaches) where different approximate models are introduced. In this work, we stick with the original equations and simply note that, thanks to the efficiency of pseudospectral methods, it is not too costly ---at least in our one-dimensional framework--- to solve the elliptic problem at each time step while maintaining high resolution. To this aim, we observe in our experiments that the aforementioned Krylov subspace iterative method GMRES is highly efficient and converges within a few iterations to the desired accuracy,  although the choice of preconditioner may turn out to be crucial. Only for extreme situations not studied here and far from the range of applicability of the equations an inversion {\em via} standard LU factorization is found necessary.   

Let us now be more precise.
We use the same Fourier pseudospectral approach as outlined in the previous 
section, i.e., we approximate the solution $u$, $\zeta$ via discrete 
Fourier transforms. This means we can treat initial data which are 
smooth and periodic or in the Schwartz class of rapidly decreasing functions 
(the latter can be treated within the finite numerical precision 
as periodic on sufficiently large 
domains) with \emph{spectral accuracy}, i.e., with a numerical error 
exponentially decaying with the number $N$ of Fourier modes. 

With this spatial discretization, both \ref{SGN} and 
\ref{WGN} are finite dimensional systems of ODEs 
coupled with a system of equations
 of the form
\begin{equation}
\left\{\begin{array}{l}
\displaystyle	\frac{\dd\hat\zeta}{\dd t}=\mathcal{G}_1(\hat \zeta,\hat \u),\\  \\
\displaystyle    \frac{\dd\hat \v}{\dd t}=\mathcal{G}_2(\hat\zeta, \hat \u,\hat \v),\\ \\
\displaystyle   \mathcal{M}(\hat\zeta)\hat \u = \hat \v
    \end{array}\right.
    \label{uform},
\end{equation}
where $\hat \zeta(t),\hat \u(t), \hat \v(t)$ are $N$-dimensional vectors, 
and $\mathcal{M}(\hat\zeta)$ is an $N$-by-$N$ matrix.
The two ODEs in system~\eqref{uform} will be integrated  with the 
standard explicit fourth order Runge-Kutta method (RK4). Note that this is 
not trivial since the system will be \emph{stiff} because of the 
three derivatives with respect to $ x $ on the right hand side of~\eqref{SGN}. The \emph{stiffness}
 implies that explicit schemes as the ones applied 
here can become inefficient because of restrictive  stability conditions on the 
time step. We briefly recall the basic concepts of Dahlquist's 
stability theory, for more details see Chapter 10 of~\cite{trefethen} 
and references therein. 
The basic idea is to consider a linear model problem 
$y'=\lambda y$, $\lambda\in\mathbb{C}$, where $\lambda$ 
is some characteristic parameter.
% (e.g., an 
%eigenvalue) of the linearization of the ODE. 
%If $\Re \lambda< 0$, the exact solution tends to zero for $t\to\infty$, and a 
%numerical solution is called stable if it also tends to zero in this limit. 
If $\Re \lambda\leq 0$, the exact solution of the model problem is bounded for positive times, and a 
numerical solution is called stable if it also bounded for positive times. 
For explicit time evolution schemes as the one used here, this 
 condition defines a domain of stability in the complex 
$z:=\lambda h$ plane (with $h$ the time step), see output~25 in Chapter 10 of~\cite{trefethen} for 
RK4. 
Considering now an ODE system ${\bf y}' 
= f(t,{\bf y})$ (in our case, ${\bf y}$ would be formed by the vectors of the Fourier coefficients $\hat \zeta$ and $\hat \v$), 
 one has to choose a time step $h$  
so that for all eigenvalues $\lambda$ of linearizations, $\lambda h$ is in the stability domain. 
Since for the SGN equations (and similarly for the WGN equations), the 
dominant contribution to the linearization is due to
order-one operators\footnote{while the \ref{SGN} and \ref{WGN} equations 
involve operators of order 3, the formulation~\eqref{SGN-WGN} (thanks to the regularizing effect of inverting the elliptic problem
\eqref{elliptic} and similarly to the Camassa-Holm equation~\eqref{CH}) shows
that
the \ref{SGN} and \ref{WGN} equations can be considered as (quasilinear)
systems involving operators of order at most one.}
and since we use a Fourier discretisation where 
$\mathcal{F}\partial_{x}=ik$ with $k=-N/2+1,\ldots,N/2$, the dominant 
contribution to $\lambda$ is of the order of $N/2$. Thus we have to 
make sure that $h=\mathcal{O}(1/N)$ which is always ensured in our 
examples. In our numerical experiments we do not 
encounter any instability issues  which would be marked by 
exponentially growing modes in time, see again Chapter~10 in 
\cite{trefethen}.  Note that if the stiffness lies in the linear part of the equations as 
for the Korteweg-de Vries equation,  
many efficient time integration schemes are known, see for instance 
\cite{etna} and references therein and the mentioned code 27 in 
\cite{trefethen}, but here the stiffness is (also) in the 
nonlinear part.

The system of linear equations in~\eqref{uform} is a convolution in the space of Fourier coefficients:
the matrix $\mathcal M(\hat \zeta)$ is constructed using (inverse) Fast Fourier Transform and 
multiplication in physical space. As already mentioned, the inversion will be done with the Krylov
 approach GMRES~\cite{GMRES}. For high 
accuracy, we use GMRES with up to 100 iterations and a stopping 
criterion of the  iteration of a relative residual of the order of 
machine precision.  As we will show below at examples, GMRES is less of a problem than in the 
inversion of the Jacobian in  the Newton iteration, though we use, unless otherwise stated, the 
same preconditioner as there. 

\subsection{Validation}

The accuracy of the code will be validated as discussed in 
\cite{etna}: the spatial resolution is controlled via the decay of 
the Fourier coefficients which is known to be exponential for 
analytical functions. Thus the order of magnitude of the highest Fourier 
coefficients gives an indication of the numerical error. The 
resolution in time will be controlled via conserved quantities of the equations. 
As mentioned in the introduction, both equations studied here have conserved quantities which will 
depend during the numerical computation on time because of 
unavoidable numerical errors. The numerical conservation of these 
quantities will give an estimation of the numerical error (it 
generally underestimates this error by 1-2 orders of magnitude, see 
the discussion in~\cite{etna}). We shall use here the third 
quantity in, respectively,~\eqref{SGNcons}  and~\eqref{WGNcons}.
 
\begin{figure}[htb!]
  \includegraphics[width=0.49\textwidth]{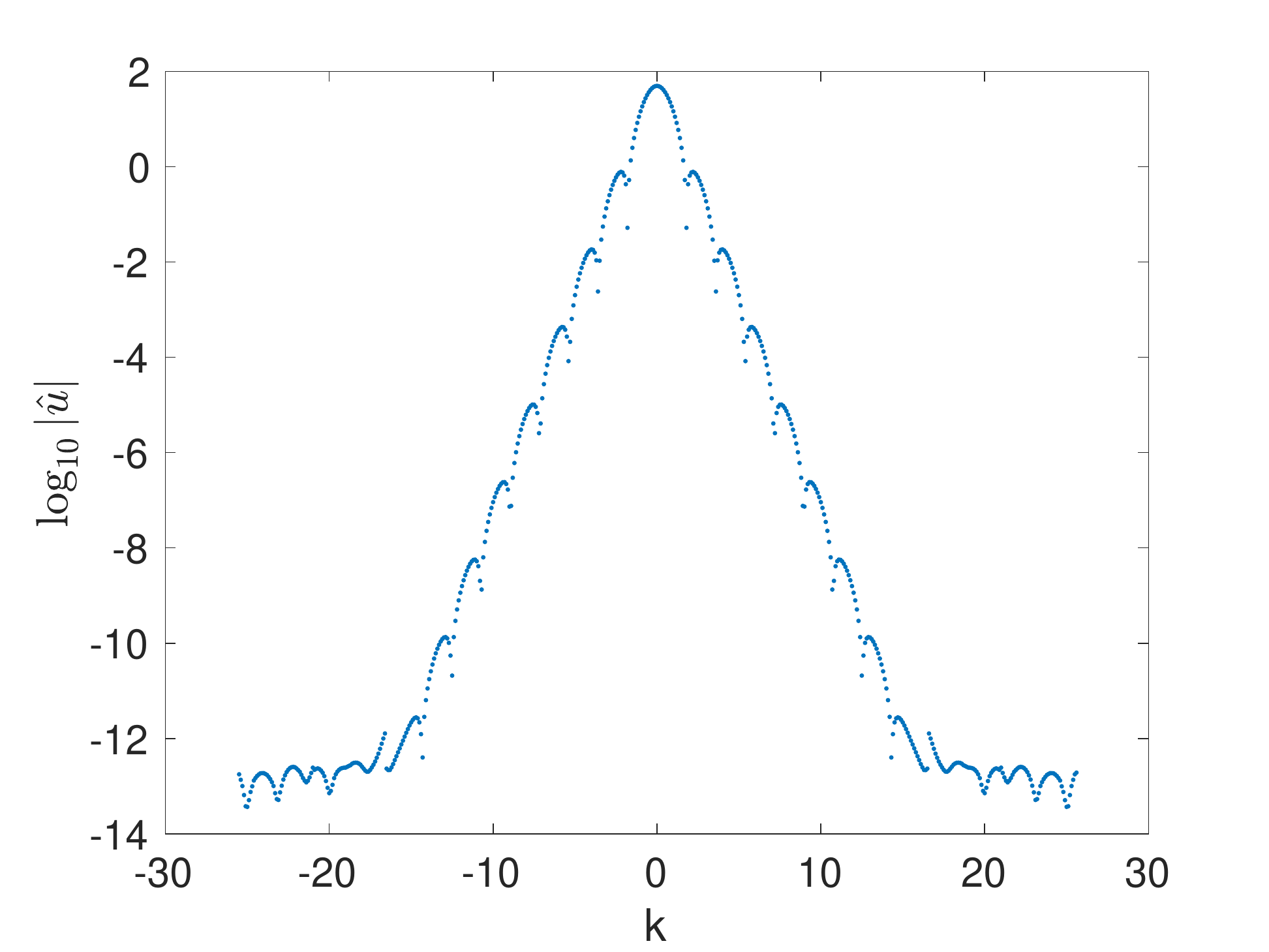}
  \includegraphics[width=0.49\textwidth]{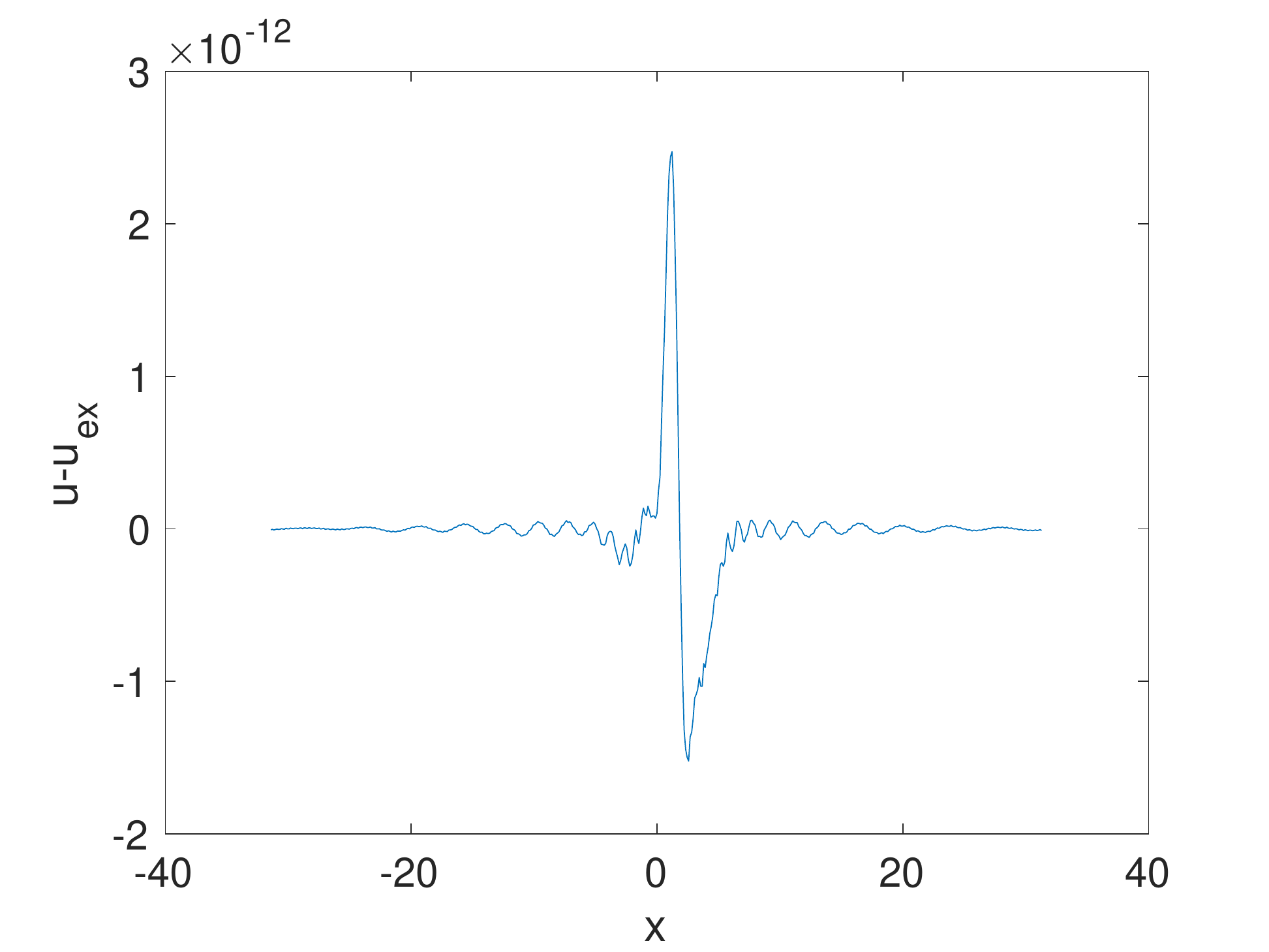}
 \caption{Left: modulus of the Fourier coefficients of the numerical 
 solution to the \ref{SGN} equations for solitary wave initial data with $c=2$ at 
 t=1; right: the difference between this solution to the \ref{SGN} equations 
 and the exact solution.}
 \label{SGNc2}
\end{figure}

To test the code we consider the solitary waves with velocity $c=2$ as an example. 
First we consider the \ref{SGN} equations with the initial 
data~\eqref{zetaSGN} for $c=2$. We use $N=2^{9}$ Fourier modes for $ 
x\in10[-\pi,\pi] $ and 
$N_{t}=2000$ time steps for $t\in[0,1]$. The relative conservation of 
the  quantities in~\eqref{SGNcons} is of the order of $10^{-14}$.
The Fourier coefficients can be seen on the left of Fig.~\ref{SGNc2}. 
They decrease to the order of $10^{-13}$. The difference between the 
numerical and the exact solution for $t=1$ can be seen on the right 
of Fig.~\ref{SGNc2}. It is of the order of $10^{-12}$,  and the 
numerical error is thus just 
an order of magnitude larger than what is indicated by the Fourier 
coefficients and the conserved quantities. 

\begin{figure}[htb!]
  \includegraphics[width=0.49\textwidth]{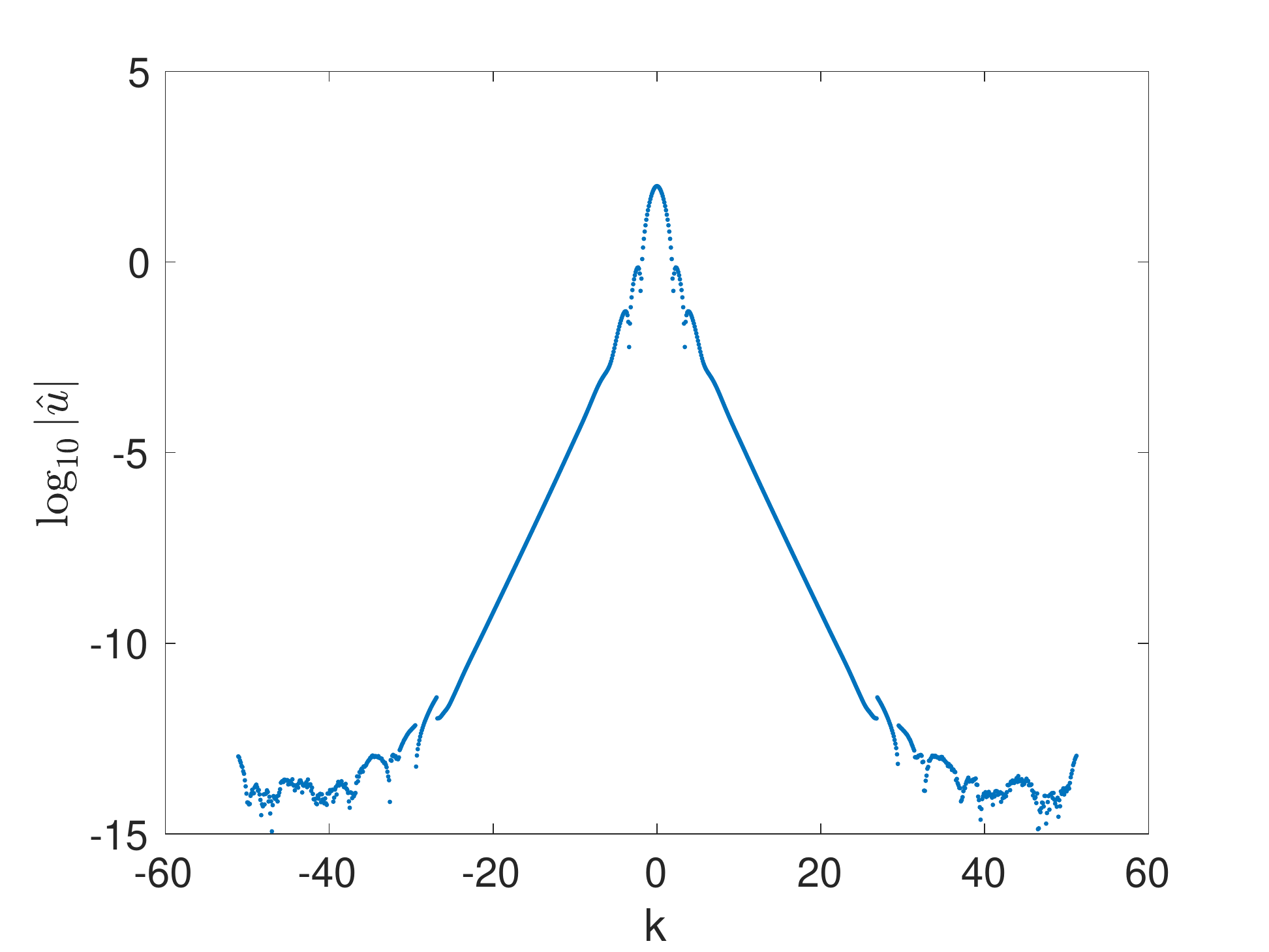}
  \includegraphics[width=0.49\textwidth]{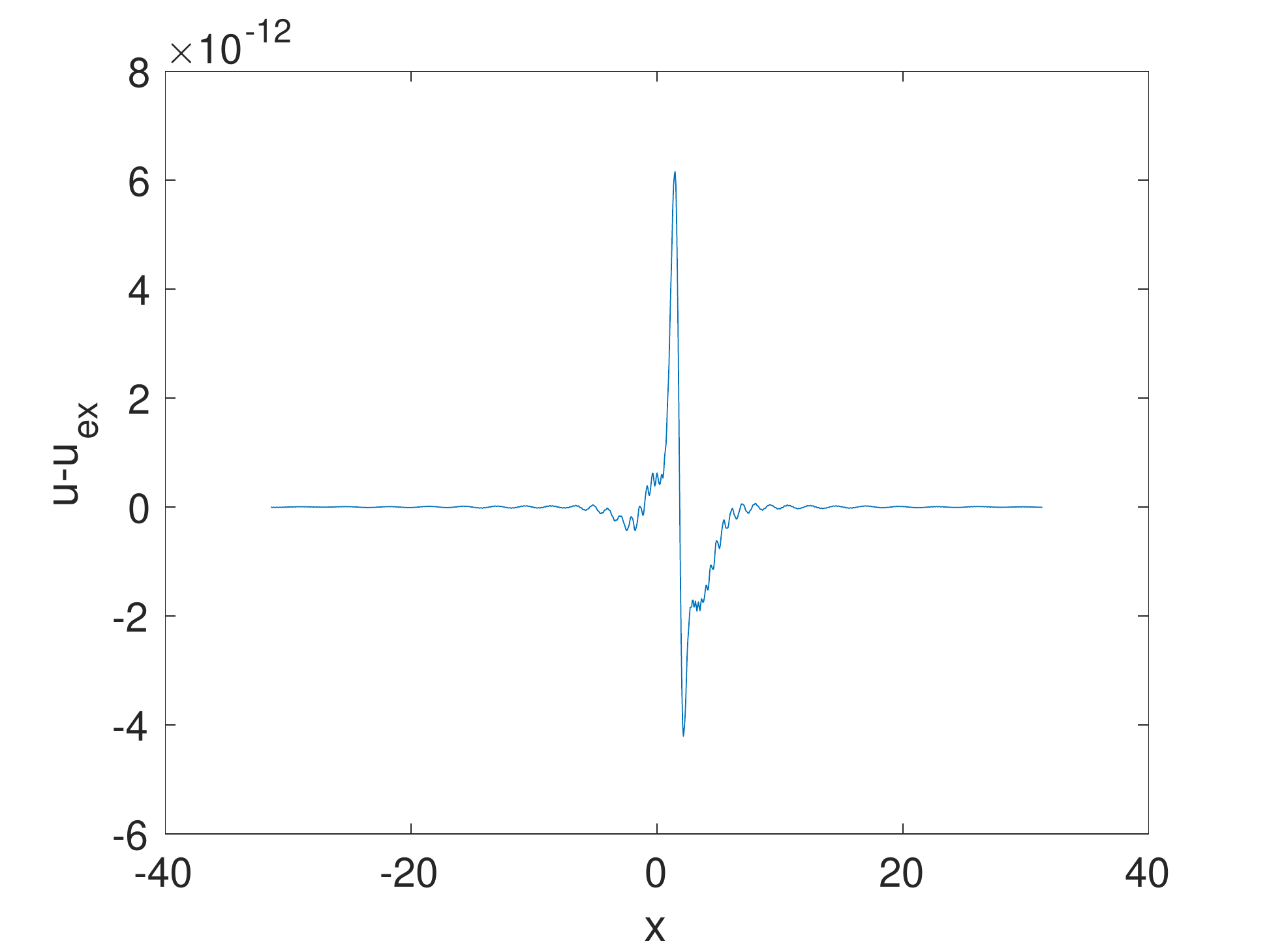}
 \caption{Left: modulus of the Fourier coefficients of the numerical 
 solution to the \ref{WGN} equations for solitary wave initial data with $c=2$ at 
 t=1; right: the difference between this solution to the \ref{WGN} equations 
 and the exact solution.}
 \label{fdSGNc2}
\end{figure}

To test the \ref{WGN} equations in a similar way, we first numerically 
construct the solitary wave for $c=2$ with $N=2^{10}$ Fourier modes. Then we 
use this numerical solution as initial data for the \ref{WGN} equations. 
This also assesses the accuracy with which the solitary wave is 
numerically constructed.  
Again we apply $N_{t}=2000$ time steps for $t\in[0,1]$. The conserved 
quantities are relatively conserved to the order of $10^{-13}$. 
The Fourier coefficients of the solution at the final time can be seen on 
the left of Fig.~\ref{fdSGNc2}, the difference with the numerically 
constructed solitary wave on the right. Obviously we reach the same 
accuracy as in the \ref{SGN} case both in terms of resolution in space as 
indicated by the decay of the Fourier coefficients for large
Fourier modes and the difference between numerical and exact 
solution. 

Note 
that in both examples here, a further increase in resolution both in 
space or time does not lead to higher accuracy.  

\section{Stability of the solitary waves}\label{S.stability}

In this section we study the behavior of solutions to both the \ref{SGN} and the \ref{WGN} equations
set initially as  various perturbations of the solitary waves. Let us describe existing results 
in the literature. In~\cite{Li01,Li02}, Li investigates 
the linear stability of the explicit solitary waves of the \ref{SGN} equations.
The first observation is that while solitary waves are critical points of a functional 
immediately stemming from the Hamiltonian structure of the equation, 
 the second variational derivative of the functional has an infinite-dimensional
(essential) negative spectrum, hence standard tools for the nonlinear
stability analysis do not apply. This comment also applies to the \ref{WGN} equations.
Then Li proves that solitary waves of the \ref{SGN} equations with sufficiently small
velocity $0<c-1\ll 1$ are (orbitally) {\em linearly stable} (see details therein)
for infinitely small and exponentially decaying perturbations. The proof is not
easily extended to the \ref{WGN} equations since the differential nature of the operators
is used. Let us also mention that Carter and Cienfuegos numerically studied in~\cite{CarterCienfuegos11}
the linear stability of {\em cnoidal waves} and found that sufficiently large or steep cnoidal waves 
exhibit linear instability, with relatively small growth rate. By nature, since the unstable modes 
are periodic with the period being a multiple of the period of the cnoidal wave, the results do not apply to solitary waves;
see the discussion in~\cite{CarterCienfuegos11}. Finally, let us also mention the work~\cite{ElGrimshawSmyth06}
where the {\em modulational stability} of small-amplitude {\em bores} of the \ref{SGN} equations is found.

In this section we consider the case of perturbations of solitary waves  with moderate up to large amplitudes (we show figures for velocities $c=2$, $c=4$ and $c=10$) and
---consistently with the related numerical experiments provided in~\cite{MitsotakisDutykhCarter17}--- 
 find that these solitary waves appear asymptotically orbitally stable. 

We consider first perturbations of the solitary waves with velocity $c=2$. We work with $N=2^{10}$ Fourier modes for $ x\in10[-\pi,\pi] $  and 
$N_{t}=2000$ time steps for the time interval $t\in[0,10]$. 
The relative conservation of  the third  quantity in~\eqref{SGNcons} for the \ref{SGN} equations remains valid up to the order of $10^{-10}$, and the corresponding one in~\eqref{WGNcons} for the \ref{WGN} equations up to the order of $10^{-9}$. 

\begin{figure}[htb!]
  \includegraphics[width=0.7\textwidth]{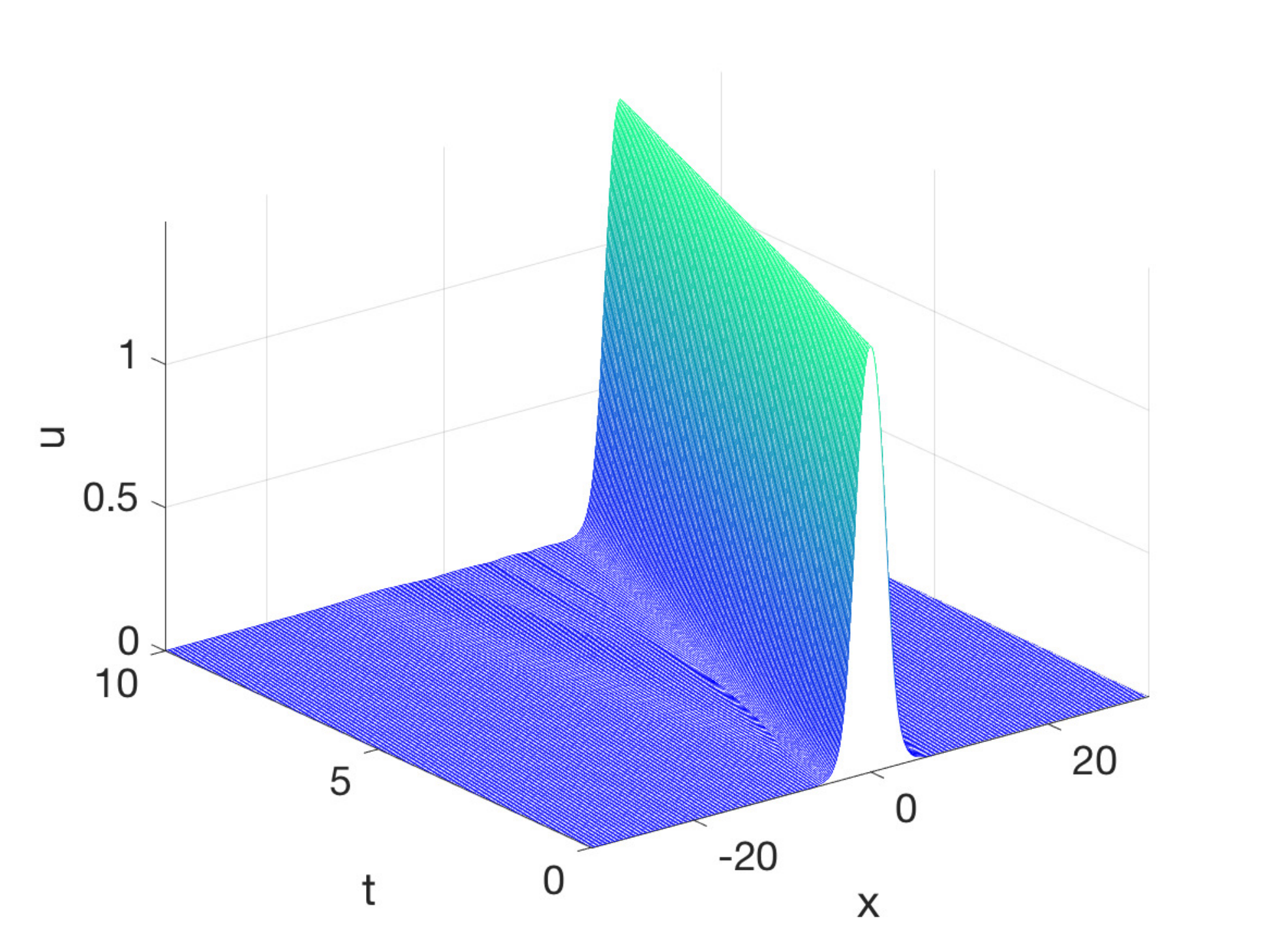}
 \caption{Solution to the \ref{SGN} equations for the initial data 
 $u(x,t=0)=0.99 u_c(x)$ for $c=2$ and $\zeta(x,t=0)=\zeta_2(x)$.}
 \label{SGNc2_099}
\end{figure}

We first study initial data of the form $(\zeta(x,t=0),u(x,t=0):=(\zeta_c(x),\lambda u_c(x))$ with 
$\lambda\in \mathbb{R}$  where $(\zeta_c,u_c)$ is the solitary wave with velocity $c$.  Similar perturbations of the 
initial data for $ \zeta $ instead of $ u $ lead to similar results, not represented here.
The solution to the \ref{SGN} equations for these initial data (with $c=2$)
and $\lambda=0.99$ can be seen in in Fig.~\ref{SGNc2_099}. There is 
some radiation propagating to the left, but the final state appears 
to be a solitary wave of slightly different mass than the perturbed solitary wave. 

\begin{figure}[htb!]
  \includegraphics[width=0.49\textwidth]{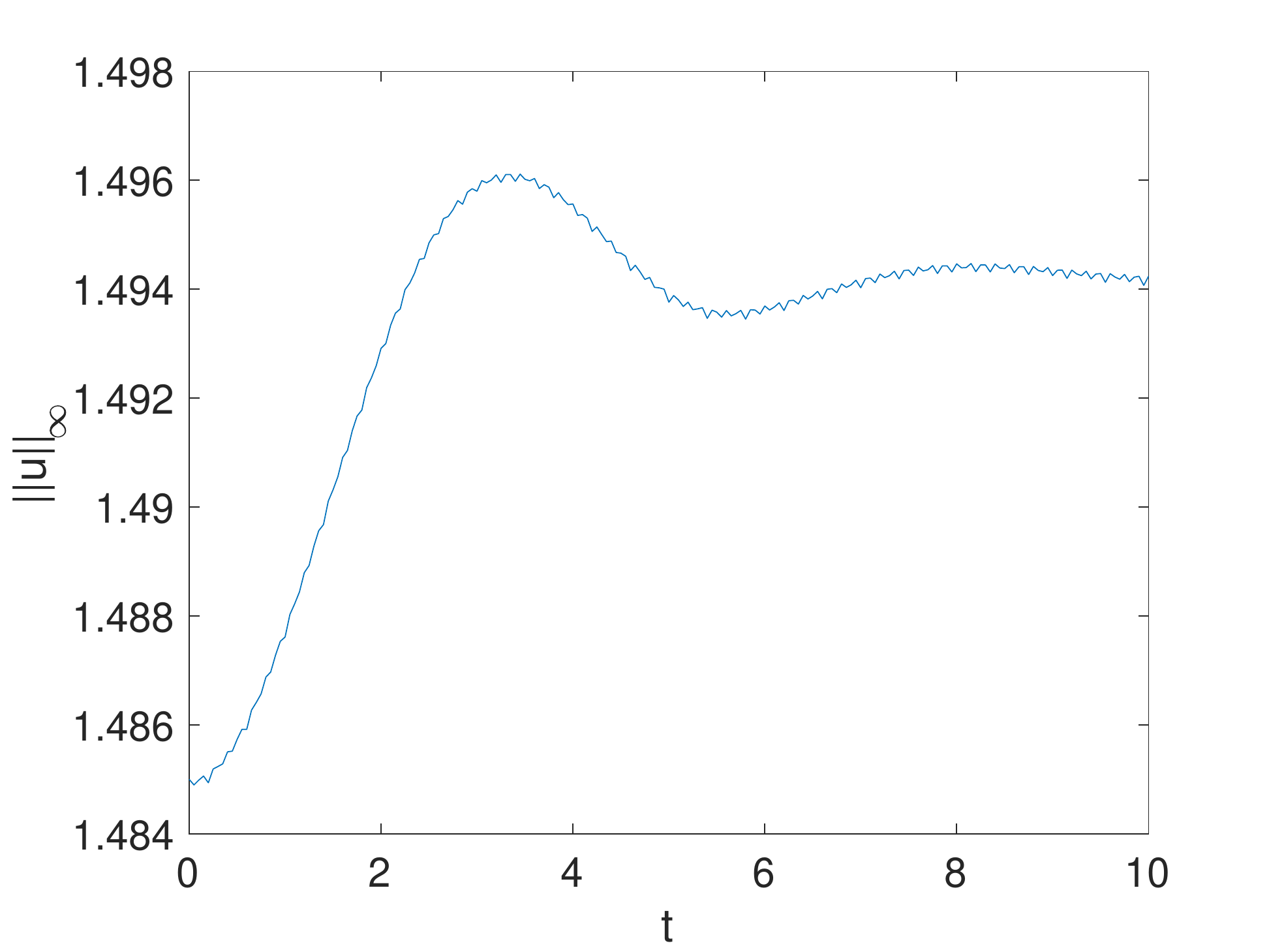}
  \includegraphics[width=0.49\textwidth]{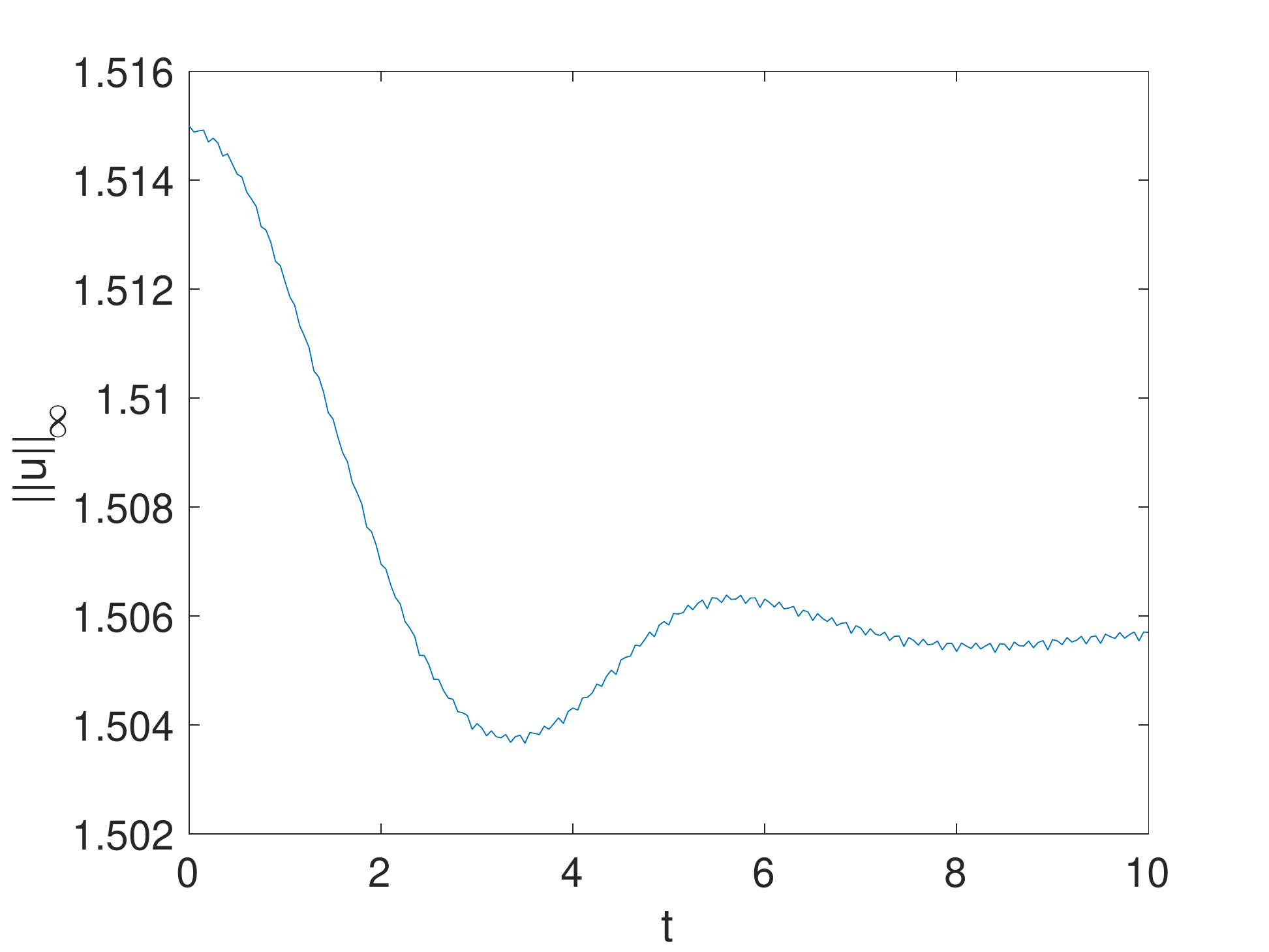}
 \caption{$L^{\infty}$ norms of the solutions to the \ref{SGN} equations for 
 the initial data  $\zeta(x,t=0)=\zeta_2(x)$ and $u(x,t=0)=\lambda u_2(x)$, on the left for 
 $\lambda=0.99$, on the right for $\lambda=1.01$.}
 \label{SGNc2_099inf}
\end{figure}

The fact that a solitary wave is approached is even more obvious looking at the 
$L^{\infty}$ norm of the solution 
(computed  on collocation points and thus only an approximation of 
the $L^{\infty}$ norm since the maximum might not be taken on a grid 
point. This explains the apparent fluctuations in the subsequent figures.)
 plotted on the left of Fig.~\ref{SGNc2_099inf}.
 The $L^{\infty}$ norm is increasing at the 
beginning and after a few oscillations appears to reach a final value. Since we 
approximate a situation on the real line by a periodic setting, the 
radiation cannot escape to infinity here which means that a final 
state cannot be reached. Extending the computation
 to a larger time interval, we observe visible oscillations due to the interaction
 with radiation starting about $t=17$.  The case $\lambda=1.01$  is very 
similar in appearance to Fig.~\ref{SGNc2_099}, therefore we do not show 
the corresponding figure, just the $L^{\infty}$ norm on the right of 
Fig.~\ref{SGNc2_099inf}. It can be seen that the norm decreases here 
before reaching its final value. 

\begin{figure}[htb!]
  \includegraphics[width=0.49\textwidth]{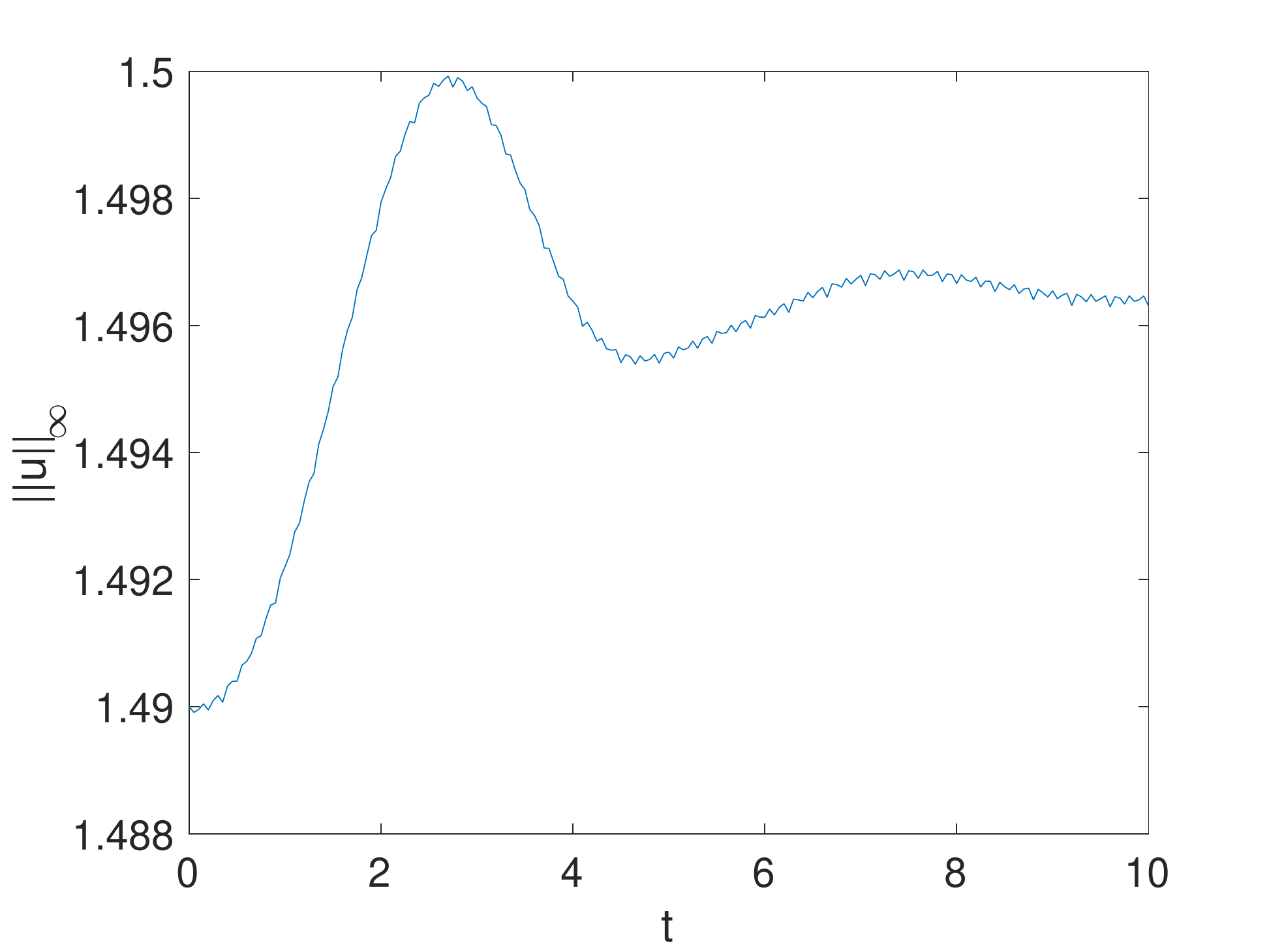}
  \includegraphics[width=0.49\textwidth]{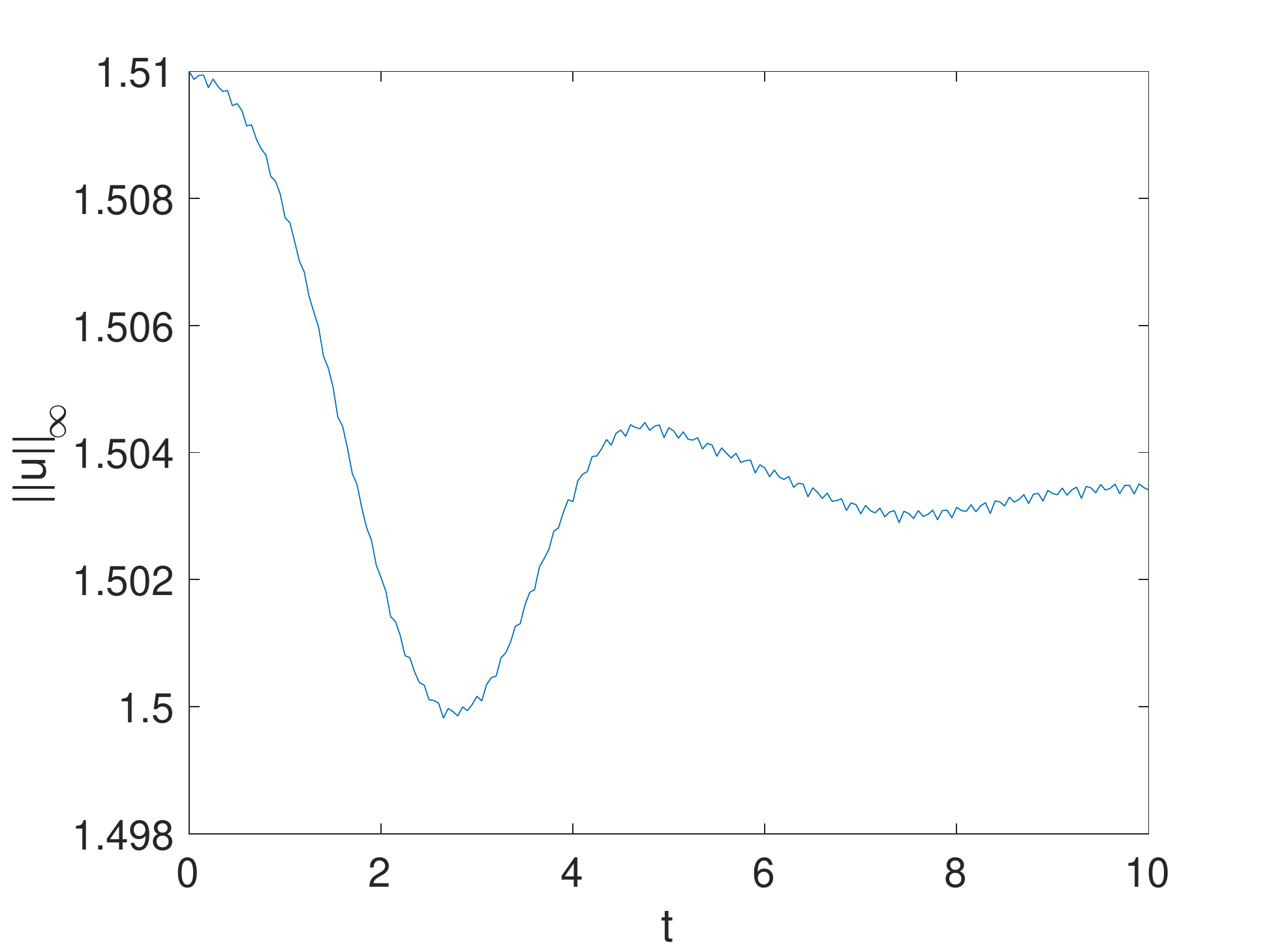}
 \caption{$L^{\infty}$ norms of the solutions to the \ref{SGN} equations for 
  initial data $\zeta(x,t=0)=\zeta_2(x)$, $u(x,t=0)=u_2(x)- 0.01\exp(-x^{2})$ on the left and 
 $u(x,t=0)=u_2(x)+ 0.01\exp(-x^{2})$ on the right.}
 \label{SGNc2gaussinf}
\end{figure}

In Fig.~\ref{SGNc2gaussinf} we show the $L^{\infty}$ norm of the 
\ref{SGN} solution with initial data given by
$(\zeta(x,t=0),u(x,t=0)):=(\zeta_2(x),u_2(x)\pm 0.01\exp(-x^{2}))$, 
on the left for the minus sign in the initial data, 
on the right for the plus sign. The situation is very similar to the 
situation of Fig.~\ref{SGNc2_099inf}, the final state appears to be a 
solitary wave of slightly different amplitude than the initial one. 

\begin{figure}[htb!]
  \includegraphics[width=0.49\textwidth]{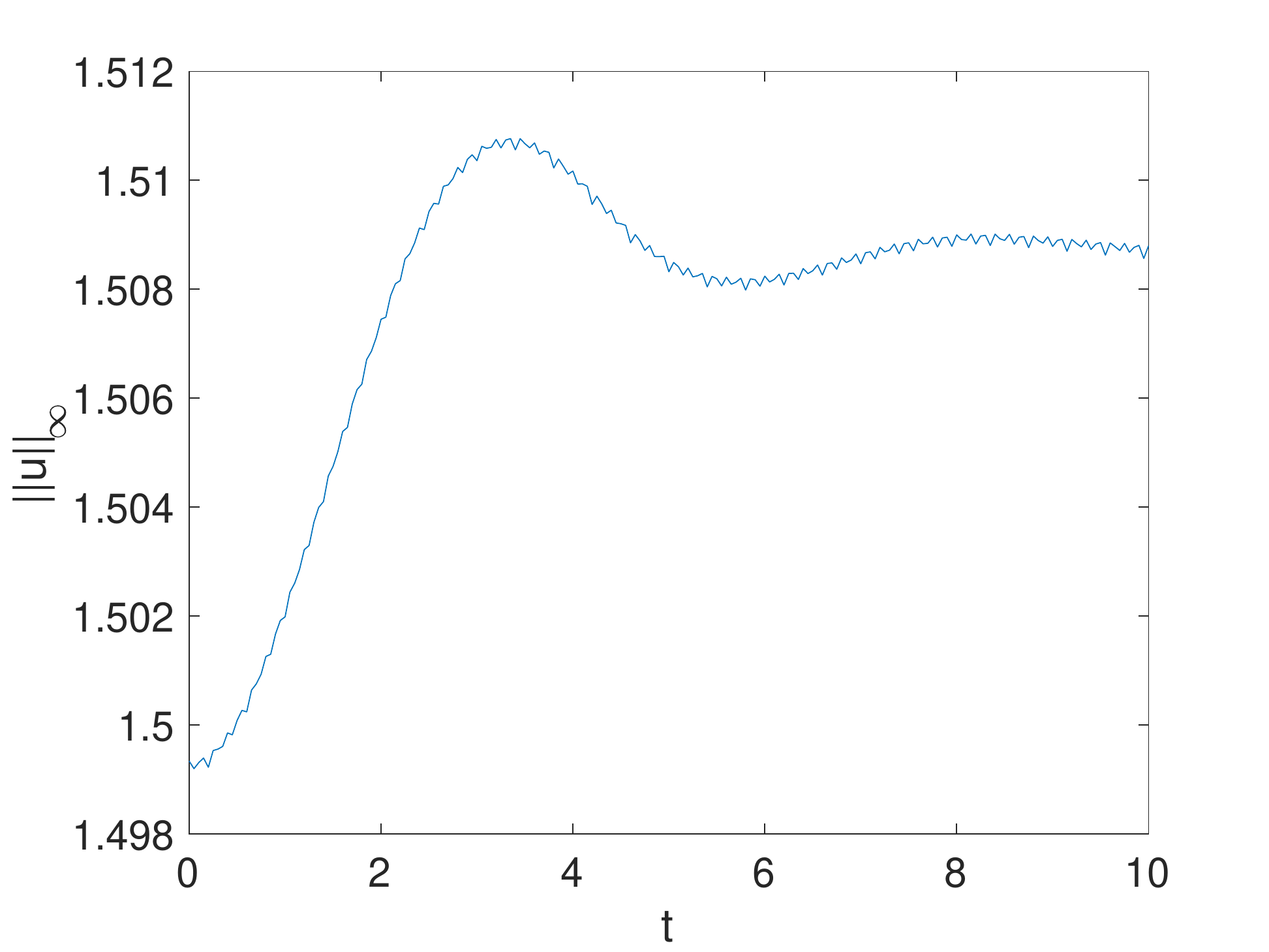}
  \includegraphics[width=0.49\textwidth]{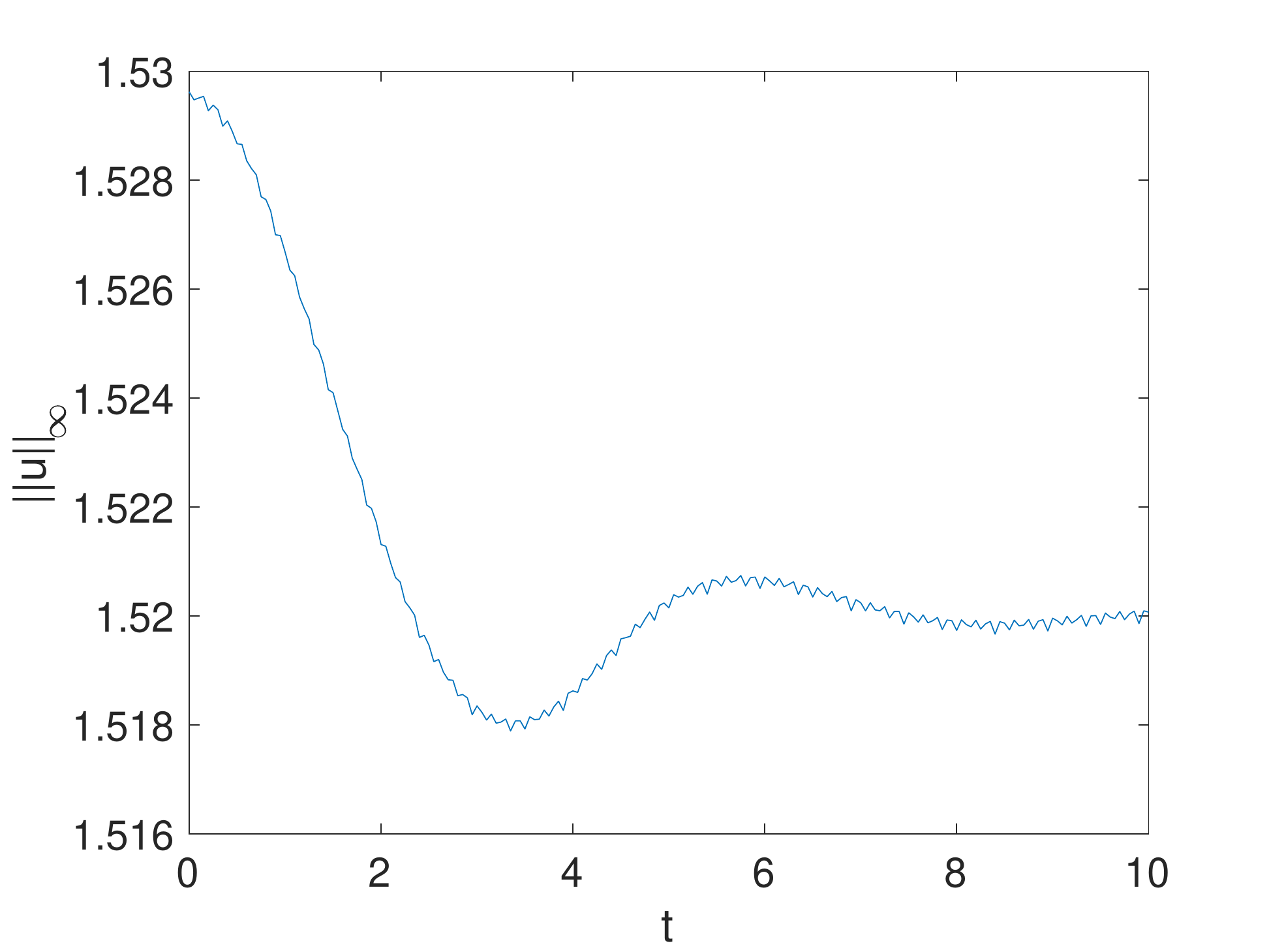}\\
  \includegraphics[width=0.49\textwidth]{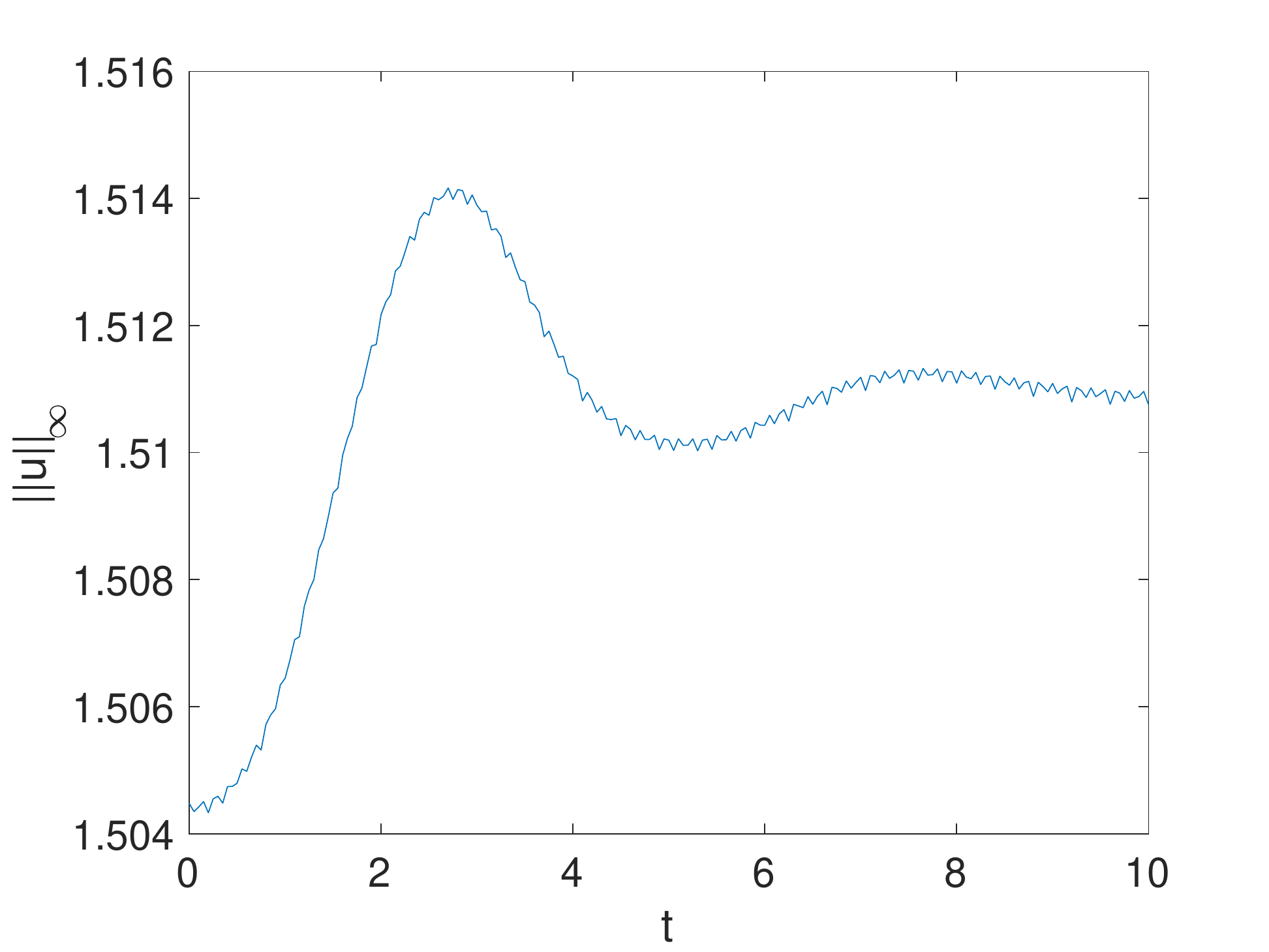}
  \includegraphics[width=0.49\textwidth]{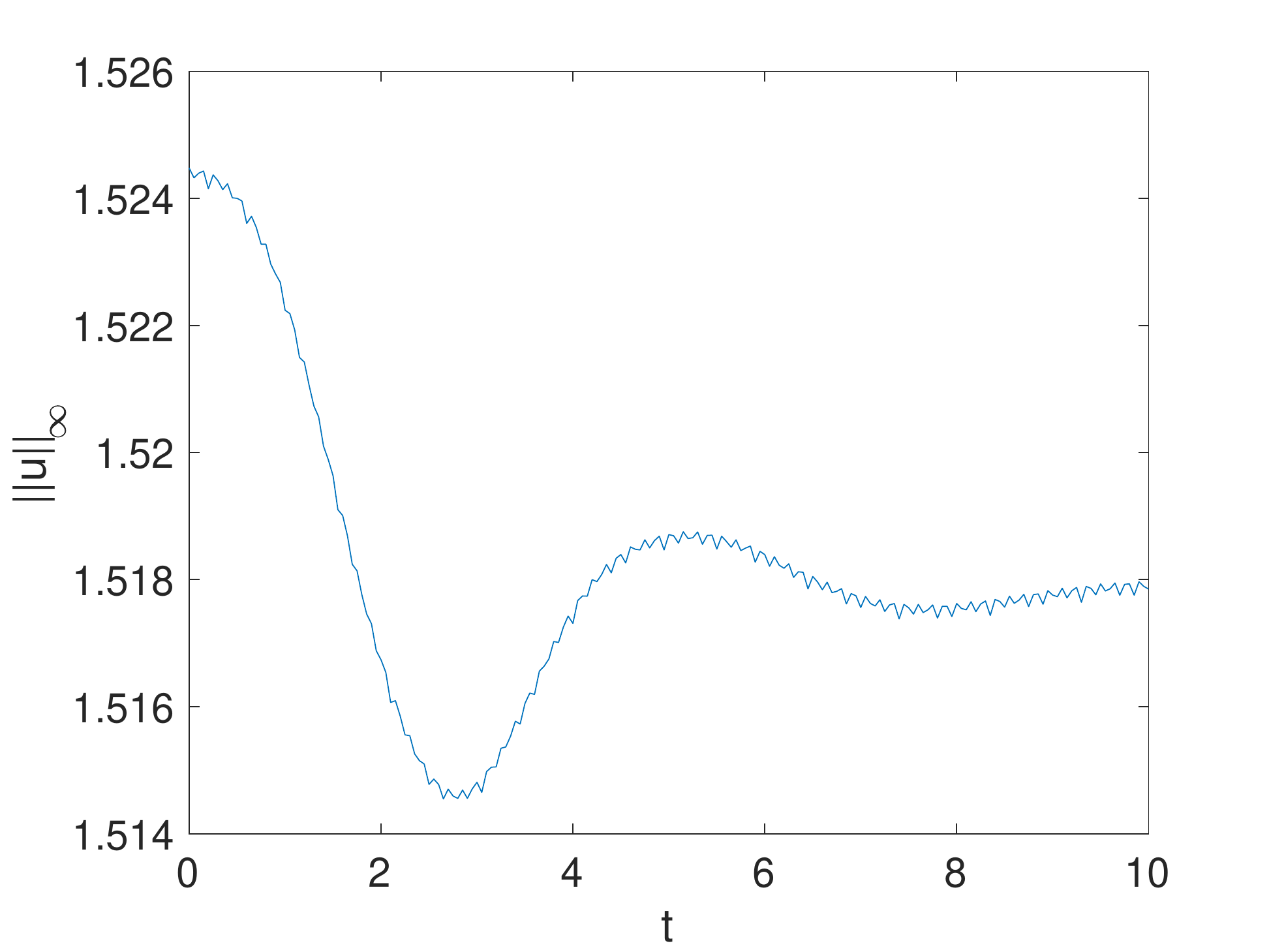}
 \caption{$L^{\infty}$ norms of the solutions to the \ref{WGN} equations for 
 the initial data  $\zeta(x,t=0)=\zeta_2(x)$ and $u(x,t=0)=\lambda u_2(x)$ in the upper row, on the left for  $\lambda=0.99$, on the right for $\lambda=1.01$; and for  $u(x,t=0)=u_2(x)\pm 
0.01\exp(-x^{2})$ in the lower row, for the minus sign on the left and 
the plus sign on the right.}
 \label{fdSGNc2gaussinf}
\end{figure}

If we consider the same perturbations for the \ref{WGN} 
solitary wave (still with $c=2$), the resulting figures are very similar as can be seen in 
Fig.~\ref{fdSGNc2gaussinf}. The final state in each example 
appears to be a solitary wave of slightly different amplitude than the perturbed 
solitary wave.

\begin{figure}[htb!]
  \includegraphics[width=0.49\textwidth]{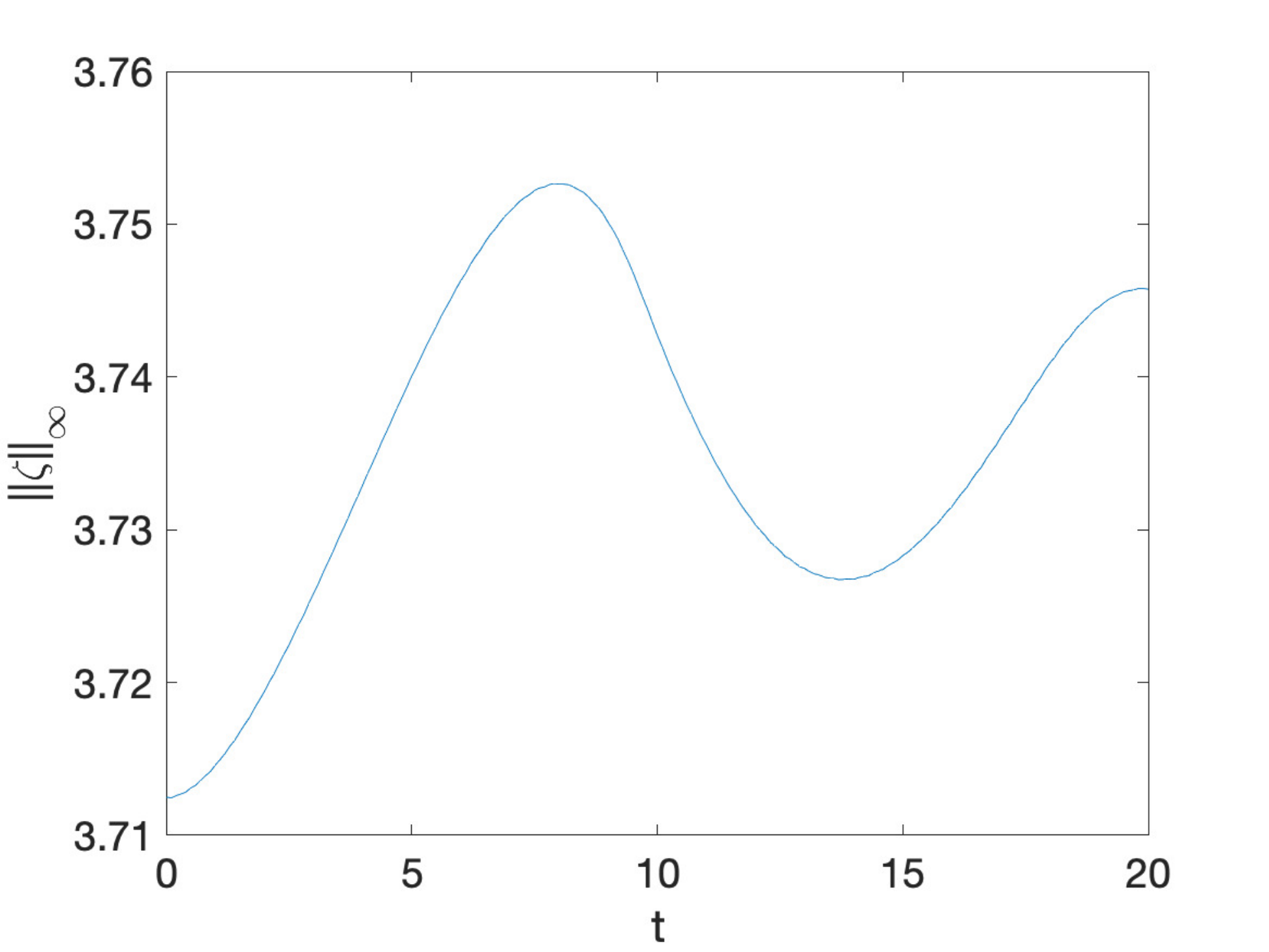}
  \includegraphics[width=0.49\textwidth]{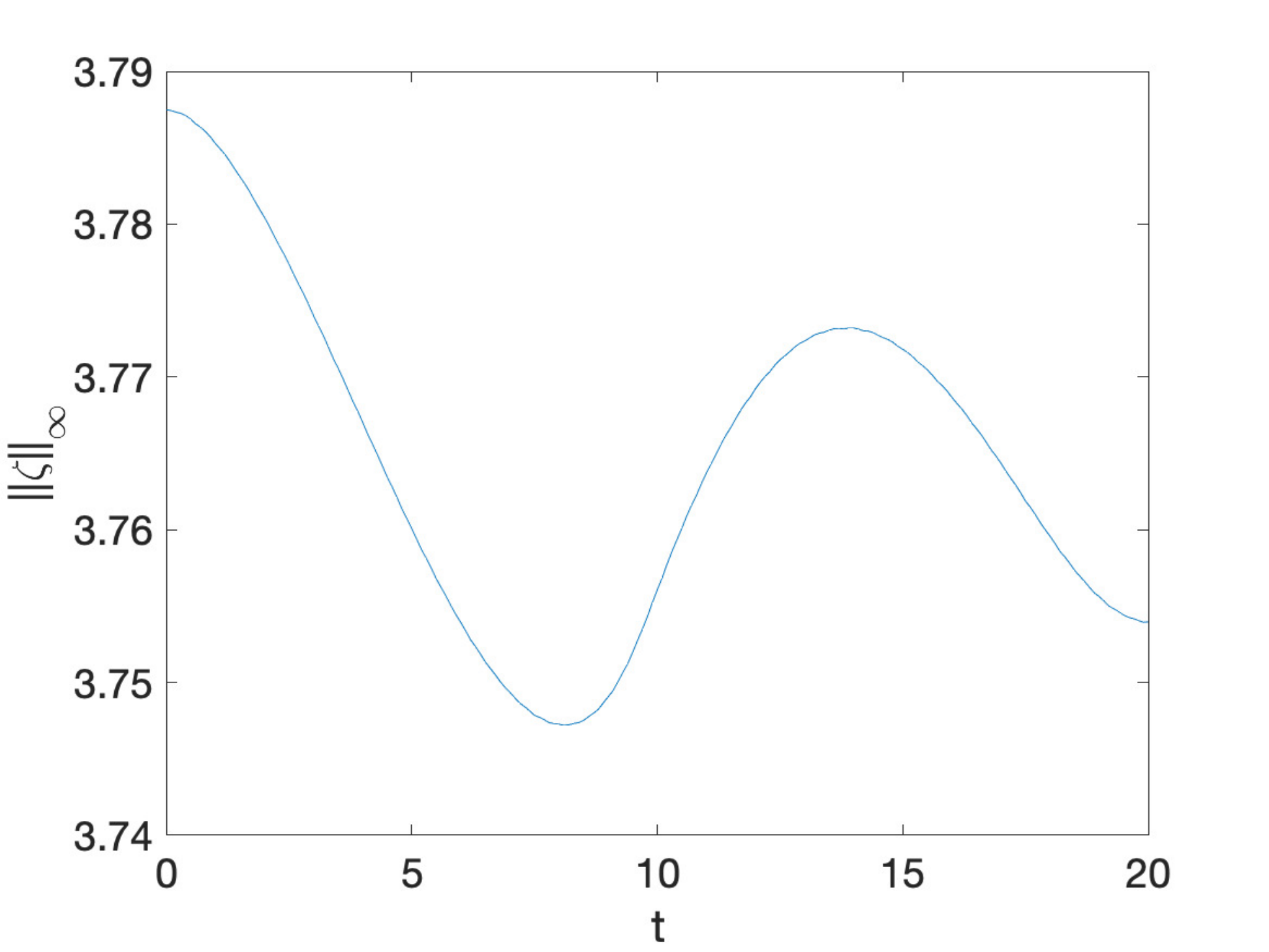}
 \caption{$L^{\infty}$ norms of the solutions to the \ref{SGN} equations for 
 the initial data $u(x,t=0)=\lambda u_4(x)$, on the left for 
 $\lambda=0.99$, on the right for $\lambda=1.01$.}
 \label{SGNc4inf}
\end{figure}

Obtaining results for $c=4$ is already much more computationally demanding.
We need to augment the space domain of computation, and hence the number of
modes in order to secure a sufficient accuracy, and also the final time of computation 
before the main wave reaches a final state. Specifically we use $ 
N=2^{12} $ Fourier modes for $ x\in 30[-\pi,\pi] $ and $ 
N_{t}=10^{4} $ time steps for $ t\in[0, 20] $. We show the $ 
L^{\infty} $ norms for the perturbations of the solitary wave of the form 
, $ \zeta(x,t=0)=\zeta_{4}(x) $ and $ u(x,t=0)=\lambda u_{4}(x)$ for $ \lambda=0.99 $ and $ 
\lambda=1.01 $ in Fig.~\ref{SGNc4inf}. It can be seen that the 
oscillations are more pronounced in this case than for the case $ c=2 
$ in Fig.~\ref{SGNc2_099inf}, but that they decrease in amplitude 
which seems to indicate that the solitary wave is again stable in this 
case. 

\begin{figure}[htb!]
  \includegraphics[width=0.49\textwidth]{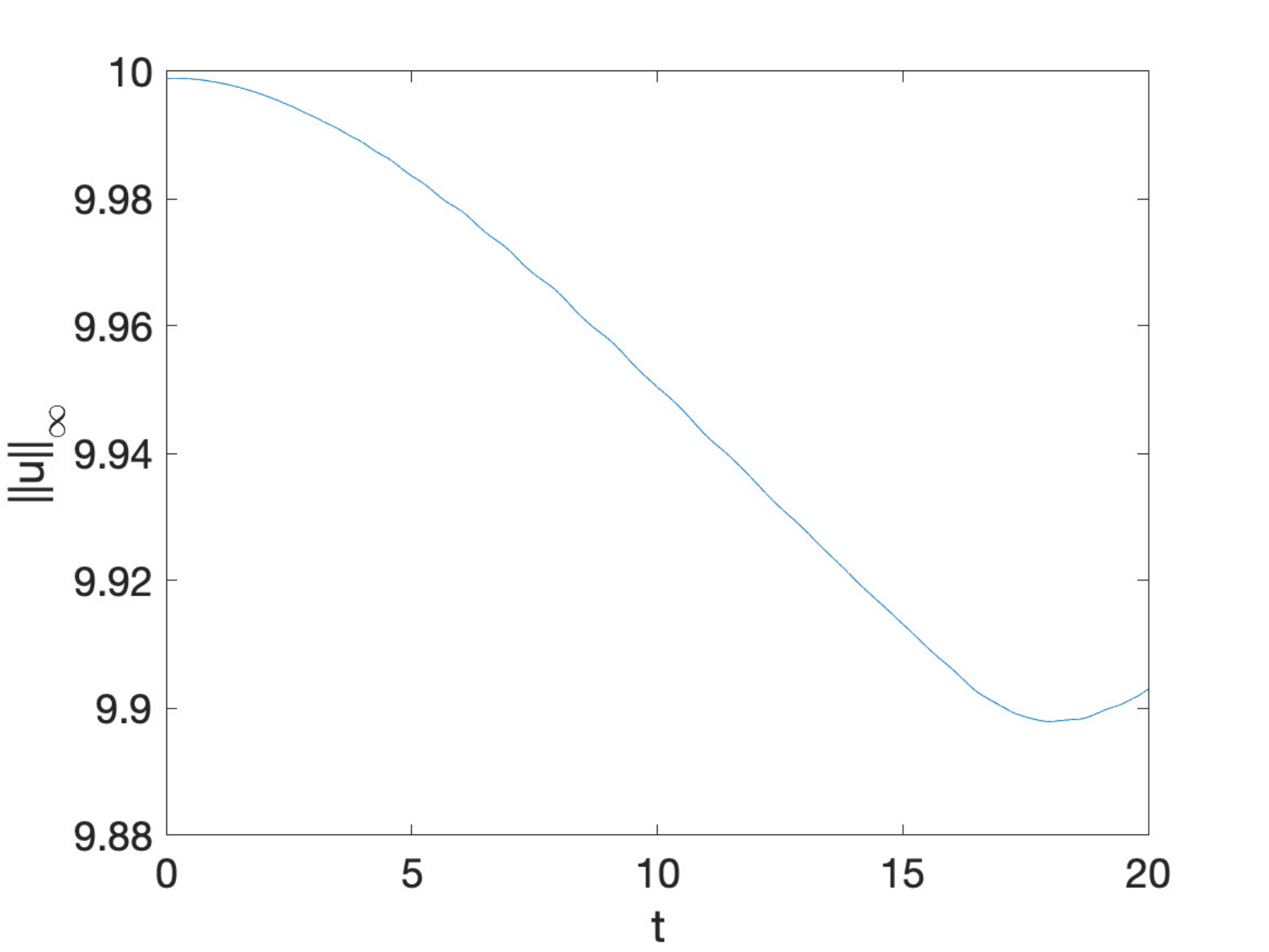}
  \includegraphics[width=0.49\textwidth]{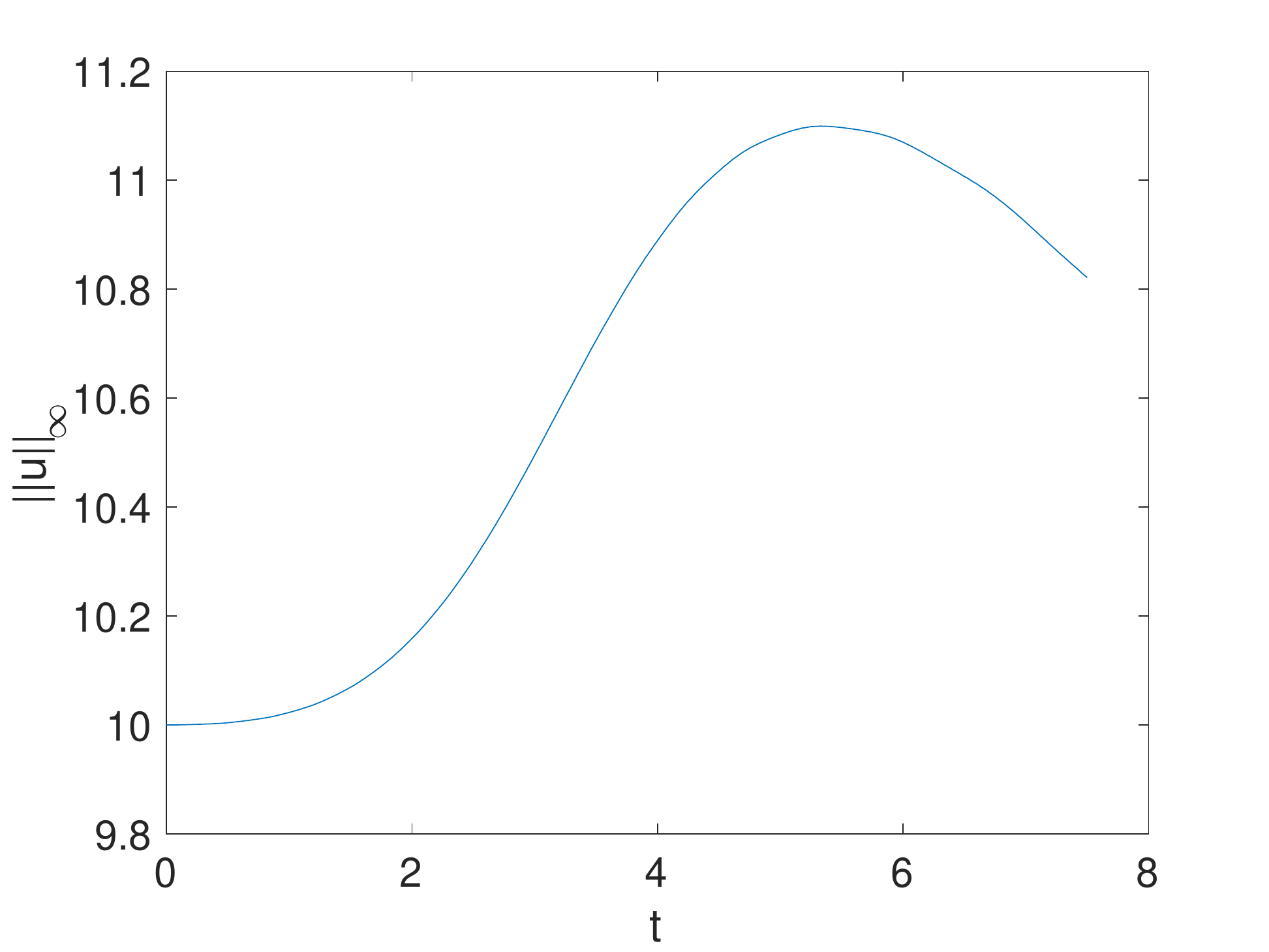}
 \caption{$L^{\infty}$ norms of the solutions to the \ref{SGN} equations for 
 the initial data $ \zeta(x,0)=\zeta_{10}(x) $, and $u(x,t=0)=1.01 u_{10}(x)$ on the left and  
 $u(x,t=0)=u_{10}(x)+0.1\exp(-x^{2})$ on the right.}
 \label{SGNc10inf}
\end{figure}

The study of even larger values of $ c $ becomes even more demanding, but can 
still be handled with the approach based on GMRES, 
which confirms it as a powerful tool (notice we use  the Fourier multiplier  $ 1+c^4k^2/3 $ as 
preconditioner). To study the case $ c=10 $, we 
use $ N=2^{11} $ Fourier modes for $ x\in 30[-\pi,\pi] $ and $ 
N_{t}=10^{4} $ time steps for the indicated time intervals. On the 
left of Fig.~\ref{SGNc10inf} we show the $ L^{\infty} $ norm of 
the solution to the \ref{SGN} equations for the initial data $ u(x,t=0)=1.01 
u_{10}(x) $, $ \zeta(x,t=0)=\zeta_{10}(x) $. The solitary wave appears to be 
again stable.
The same is true for a perturbation of the form 
$ u(x,0)=u_{10}(x)+0.01\exp(-x^{2}) $, although the observed oscillation of the amplitude 
is about ten times larger than the initial perturbation for $ u $, and with another factor of ten for the oscillation in $ \zeta $. If we consider a larger perturbation of the form 
$ u(x,0)=u_{10}(x)+0.1\exp(-x^{2}) $, the $ L^{\infty} $ norm on 
the right of Fig.~\ref{SGNc10inf}
seems to grow beyond what can be seen as a perturbative regime. We expect no blow-up since the norms decrease after some time, but the example 
makes clear that  \ref{SGN} solitary waves with large velocities $ c $ are  in applications
more easily affected by perturbations than, for instance, Korteweg-de Vries solitons.

\begin{figure}[htb!]
  \includegraphics[width=0.49\textwidth]{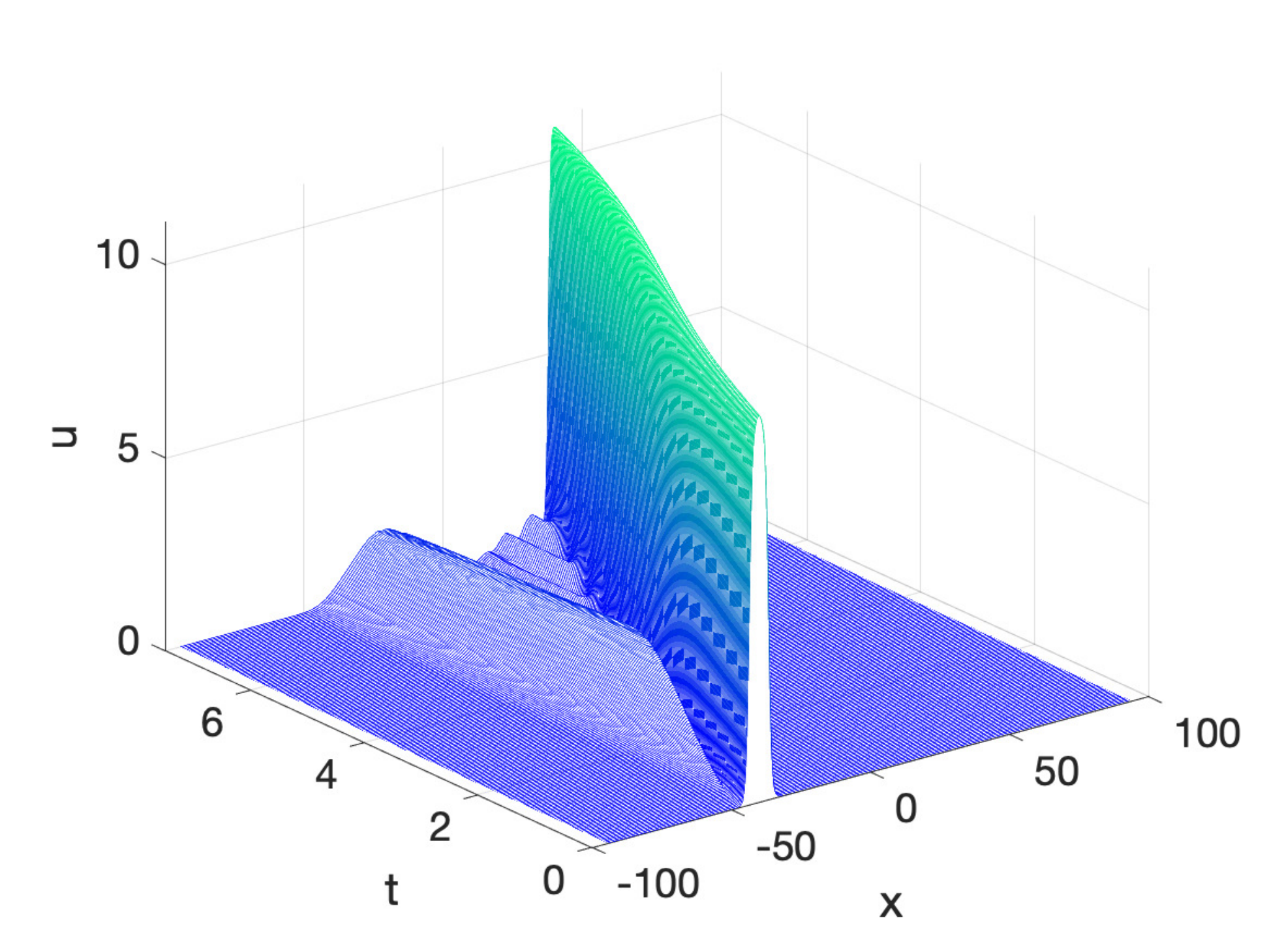}
  \includegraphics[width=0.49\textwidth]{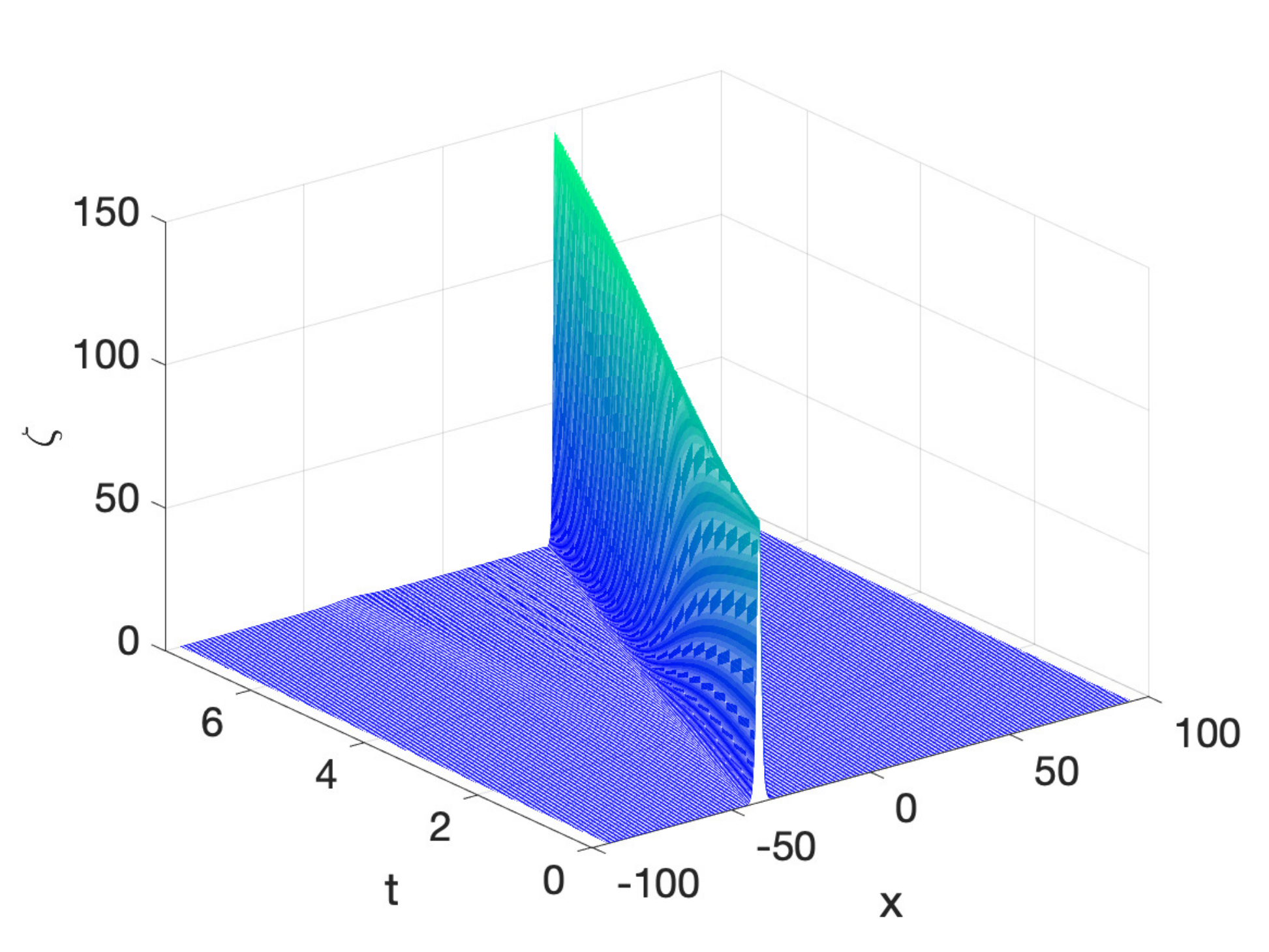}
 \caption{Solution to the \ref{SGN} equations for 
 the initial data  
 $u(x,t=0)=u_{10}(x+40))+0.1\exp(-((x+40)^{2})$, on the left $ u $, on the 
 right $ \zeta $.}
 \label{SGNc10water}
\end{figure}

This becomes even clearer if we look at the solutions to the \ref{SGN} equations for 
these initial data in Fig.~\ref{SGNc10water}. A strong growth 
especially in $ \zeta $ cannot be interpreted as a perturbation of the 
initial structure, yet the structures moving to the left do not 
appear to be solitary waves of smaller amplitude (the bumps in the figures cannot be 
fitted to a solitary wave). 

In all our experiments, considering the solitary waves of the \ref{WGN} equations
 instead of the \ref{SGN} equations
yields no qualitative difference. In particular, we observed no coherent state ---
other than the main wave--- emerging from the perturbed solitary wave, which
motivates our assertion that solitary waves appear dynamically stable.

We also performed a spectral stability analysis (not represented here)
essentially amounting to the large-period limit  of that in~\cite{CarterCienfuegos11} for cnoidal waves. 
We found that, while the linearized \ref{SGN} and \ref{WGN} equations about cnoidal waves 
always exhibit unstable modes, the rate of instability decreases as the period grows, and that no unstable mode for solitary waves is captured in this way.

\section{Emergence of modulated oscillations}\label{S.DSW}

Both the \ref{SGN} and \ref{WGN} equations reduce to the so-called
Saint-Venant or shallow-water system in the limit of infinitely long wavelength.
This is easily seen when considering the evolution equations  of a 
one-parameter family of initial data varying on a scale of order 
$1/\delta$,  for times of order $1/\delta$, where $\delta\ll1$.
Rescaling the coordinates $x\mapsto \delta x$, $t\mapsto \delta t$  
(and once again setting $g=d=1$ by similar rescaling) yields for the \ref{SGN} equations
\begin{equation}\label{SGN-delta}
      \left\{ \begin{array}{l}
      \partial_t\zeta+\partial_x(hu) =0,\\ \\
   \partial_t\big(u-\frac{\delta^2}{3h}\partial_x(h^3\partial_x u)\big)+\partial_x\zeta+u\partial_x u=\\
   \hspace{3cm}\delta^2\partial_x \big(\frac{u}{3h}\partial_x(h^3\partial_x u)+\frac12 h^2(\partial_x u)^2\big),
      \end{array}\right.
   \end{equation}
   where in an abuse of notation, we have kept the same notation for the 
   functions depending on $\delta$ as for the case $\delta=1$; and for the \ref{WGN} equations
   \begin{equation}\label{WGN-delta}
         \left\{ \begin{array}{l}
         \partial_t\zeta+\partial_x(hu) =0,\\ \\
      \partial_t\big(u-\frac{\delta^2}{3h}\partial_x\F^\delta(h^3\partial_x\F^\delta u)\big)+\partial_x\zeta+u\partial_x u\\
       \hspace{3cm}=\delta^2\partial_x \big(\frac{u}{3h}\partial_x\F^\delta(h^3\partial_x\F^\delta u)+\frac12 h^2(\partial_x\F^\delta u)^2\big),
         \end{array}\right.
      \end{equation}
   where $\F^\delta$ is the Fourier multiplier defined by
   \[\forall \varphi\in L^2(\RR) , \quad    \widehat{\F^\delta \varphi}(\xi)=F(\delta|\xi|)\widehat{\varphi}(\xi) \qquad \text{ where }  F(k)=\sqrt{\frac{3}{|k|\tanh(|k|)}-\frac3{ |k|^2}}.\]
Formally setting $\delta=0$, both~\eqref{SGN-delta} and~\eqref{WGN-delta} reduce to
the aforementioned shallow-water equations
\begin{equation}\label{SV}
\left\{\begin{array}{l}
      \partial_t\zeta+\partial_x(h u) =0,\\ \\
\partial_t u+\partial_x\zeta+u\partial_x u=0.
\end{array}\right.
\end{equation}
This formal derivation can be made rigorous for a class of sufficiently regular initial data satisfying the
non-cavitation assumption $\inf_\RR (1+\zeta)>0$, and  over a time interval of order 1
(i.e. $1/\delta$ in non-rescaled coordinates); see~\cite[Section~6.1.2]{Lannes}. For longer times, it is well-known that~\eqref{SV} 
may generate shock singularities, and the large wavelength assumption becomes invalid
before the singularity occurs. In this case it is well-documented that solutions to dispersive modified systems
may develop zones of rapid modulated oscillations in place of the shocks, 
which may eventually evolve into fully developed \emph{dispersive shock waves} or {\em solitary wave resolution},
depending on the asymptotic properties of the data as $|x|\to\infty$.
The literature on the subject is vast, and we refer to~\cite{GK12} for an overview and references in the case of
the Korteweg-de Vries (KdV) equation as a perturbation of the inviscid Burgers (iB) equation, where a complete
 asymptotic description is available, and to~\cite{ElHoefer16} for an introduction to the modulation theory for more general
  equations and an extensive list of references. Let us also specifically mention~\cite{ElGrimshawSmyth06,TkachenkoGavrilyukShyue20} 
  for a description of the Whitham modulation theory in the case of the \ref{SGN} equations, and 
\cite{LeGavrilyukHank10,DiasMilewski10,MitsotakisIlanDutykh14,MitsotakisDutykhCarter17,PittZoppouRoberts18,GavrilyukNkongaShyueEtAl20} 
for some numerical experiments. In contrast to these works, we do not consider here steplike 
initial data, but study the appearance of modulated oscillations from smooth rapidly decaying initial 
data, from which solitary wave resolution is the expected large-time asymptotic behavior. 

We will observe the emergence of modulated oscillations from initial data set as unidirectional waves for~\eqref{SV}. Based on Riemann invariants, the solution to~\eqref{SV} with initial data 
satisfying $u(x,t=0) = 2\sqrt{1+\zeta(x,t=0)} - 2$ can be described as a {\em simple wave} with
$r(x,t):=u(x,t) + 2\sqrt{1+\zeta(x,t)} -2$ satisfying the iB equation
\begin{equation}\label{iB}
 \partial_t r +(1+\tfrac34 r)\partial_x r =0.
 \end{equation}
Any spatially localized smooth solution to~\eqref{iB} develops a shock in finite time.
We consider in the following initial data of the form 
$\zeta(x,t=0)=\exp(-(x-x_{0})^{2})$, where $x_{0}$ is a constant whose role is to avoid the 
propagation of oscillations beyond the boundary of the computational 
domain.

\begin{figure}[htb!]
  \includegraphics[width=0.49\textwidth]{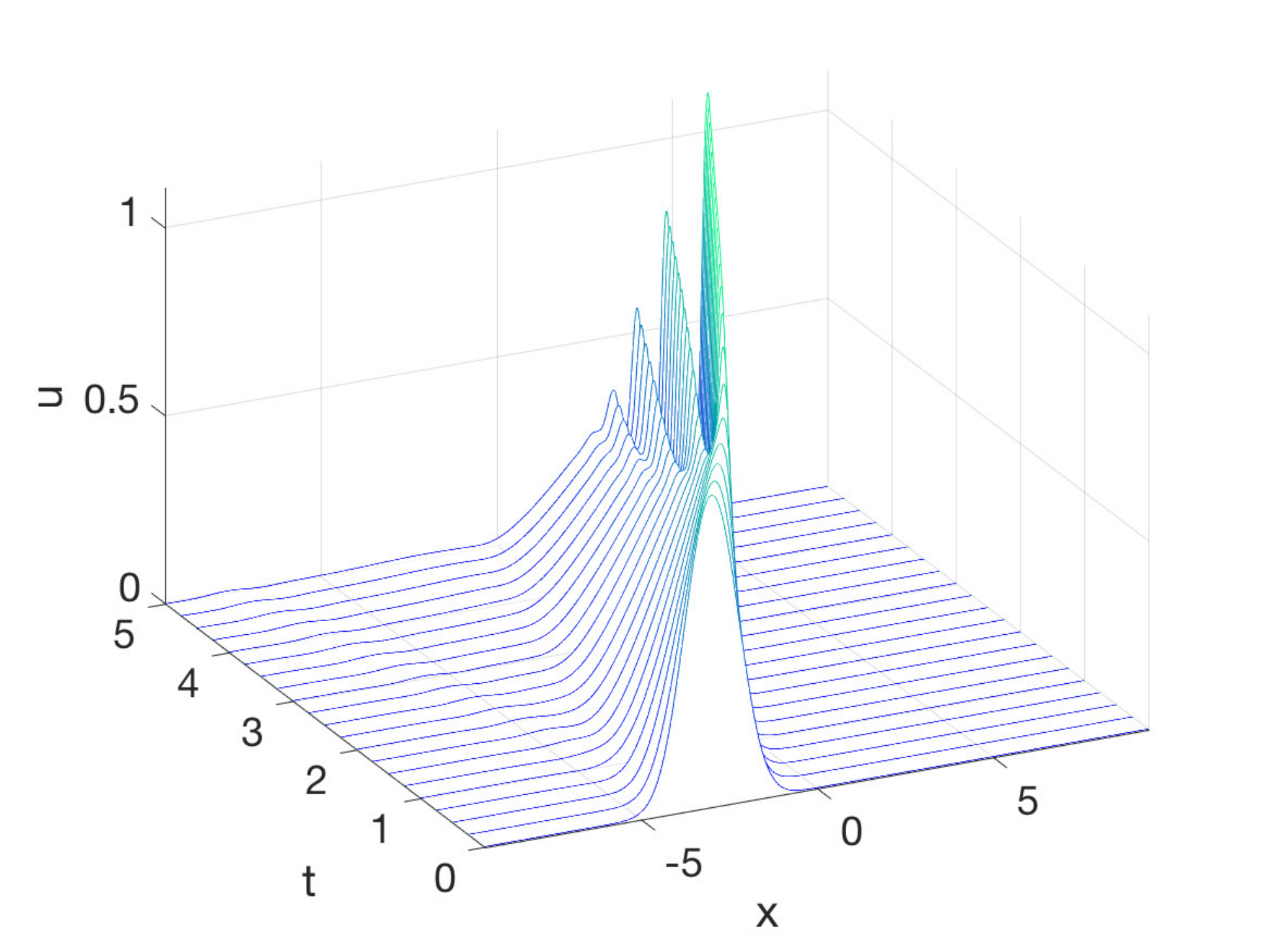}
  \includegraphics[width=0.49\textwidth]{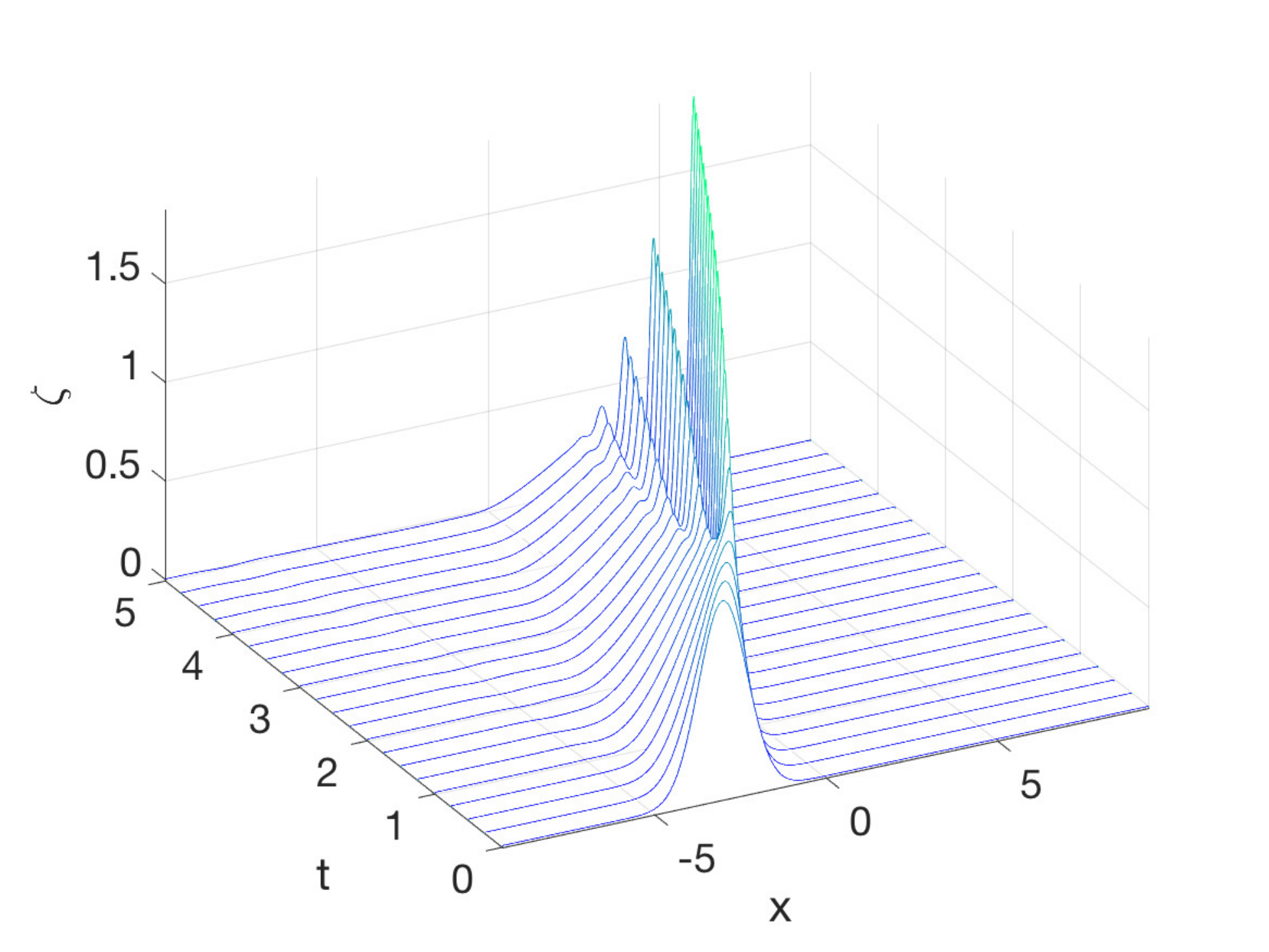}
 \caption{Solution to the SGN equations~\eqref{SGN-delta} with 
 $\delta=0.1$ for $\zeta(x,t=0)=\exp(-(x+3)^{2})$, 
 $u(x,t=0) = 2\sqrt{1+\zeta(x,t=0)} - 2$ in dependence of time; on the left 
 $u$, on the right $\zeta$.}
 \label{gauss1e1}
\end{figure}

We first consider the SGN equations~\eqref{SGN-delta} with $\delta=0.1$. 
We use  $N=2^{10}$ Fourier modes  on the computational domain $3[-\pi,\pi]$. 
We set $x_{0}=-3$, and compute $N_{t}=10^{4}$ time steps for $0\leq t\leq 5$. The 
solution can be seen as a function of time in 
Fig.~\ref{gauss1e1}, and the formation of modulated oscillations in place of the non-dispersive shock is clearly visible. 

\begin{figure}[htb!]
  \includegraphics[width=0.49\textwidth]{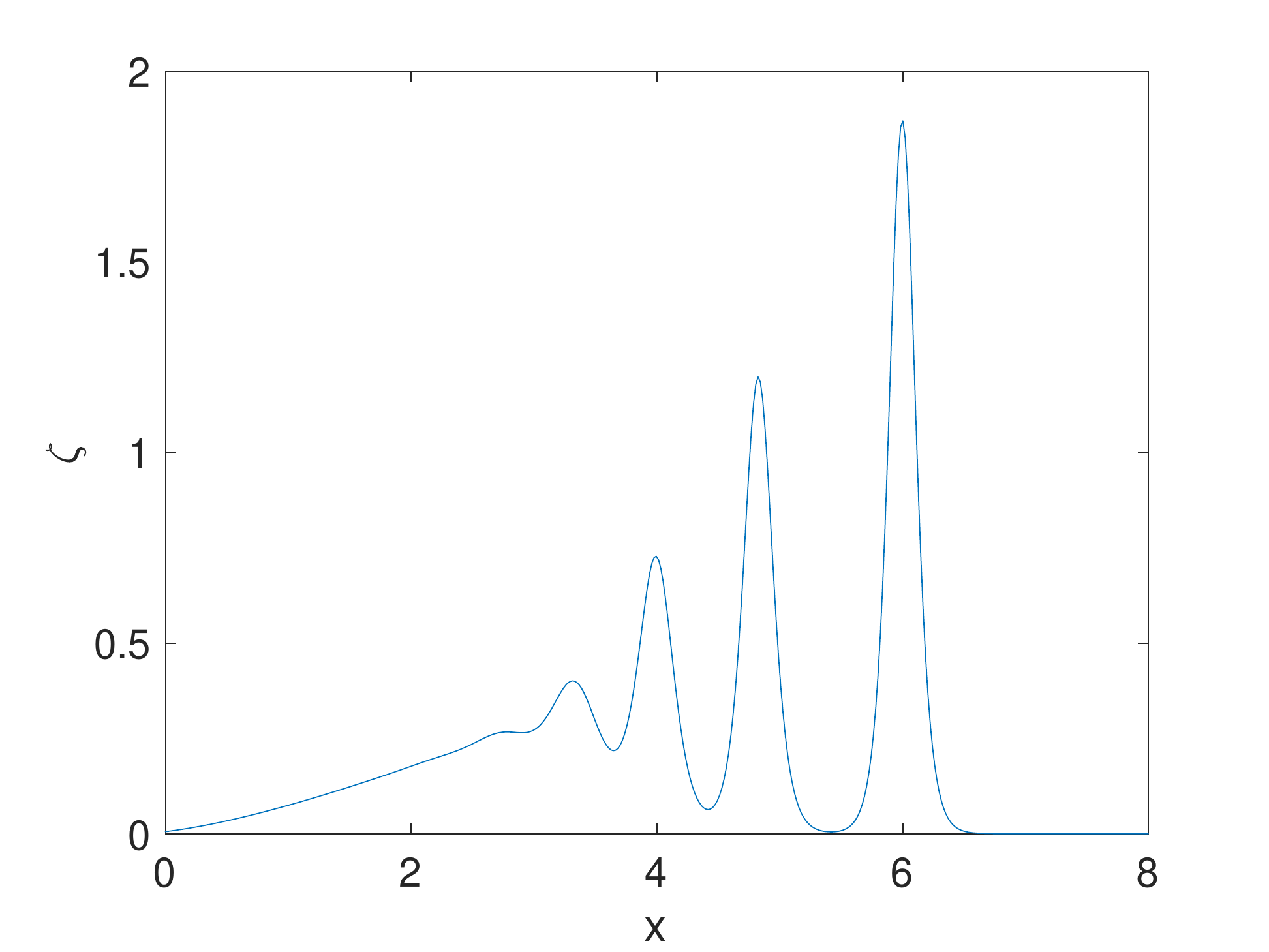}
  \includegraphics[width=0.49\textwidth]{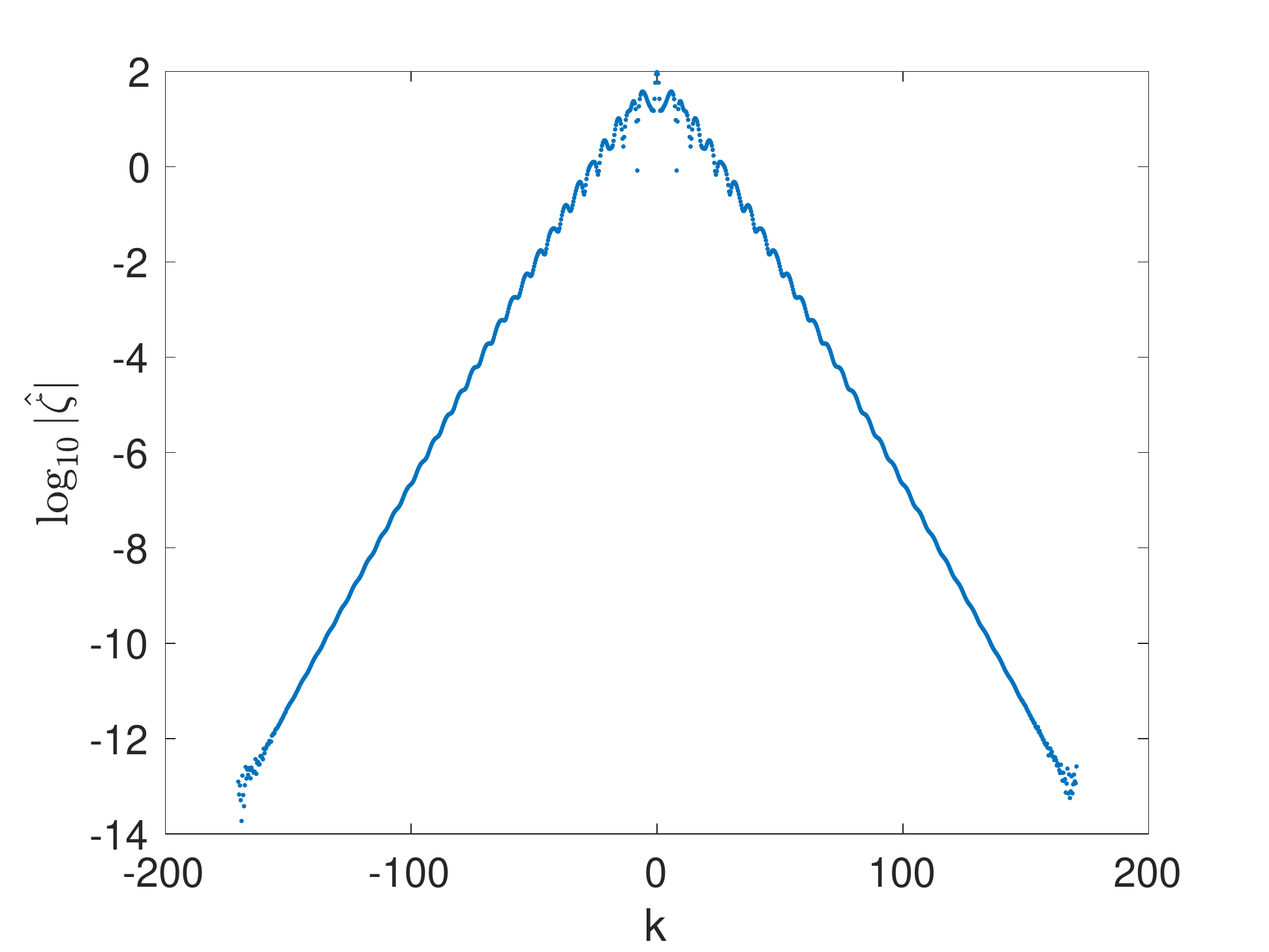}
 \caption{Function $\zeta$ of Fig.~\ref{gauss1e1}  at $t=5$ on the 
 left, and the Fourier coefficients  on the right.}
 \label{gauss1e12}
\end{figure}

A cross-section of the plot on the right of Fig.~\ref{gauss1e1},
specifically the solution $\zeta$ at time $t=5$, is shown on the left of Fig.~\ref{gauss1e12}.
One sees that the rightmost oscillation has already evolved into a localized solitary wave,
indicating the onset of solitary wave resolution. The Fourier coefficients on the right of the same figure show that the solution 
is numerically well resolved.

\begin{figure}[htb!]
  \includegraphics[width=0.49\textwidth]{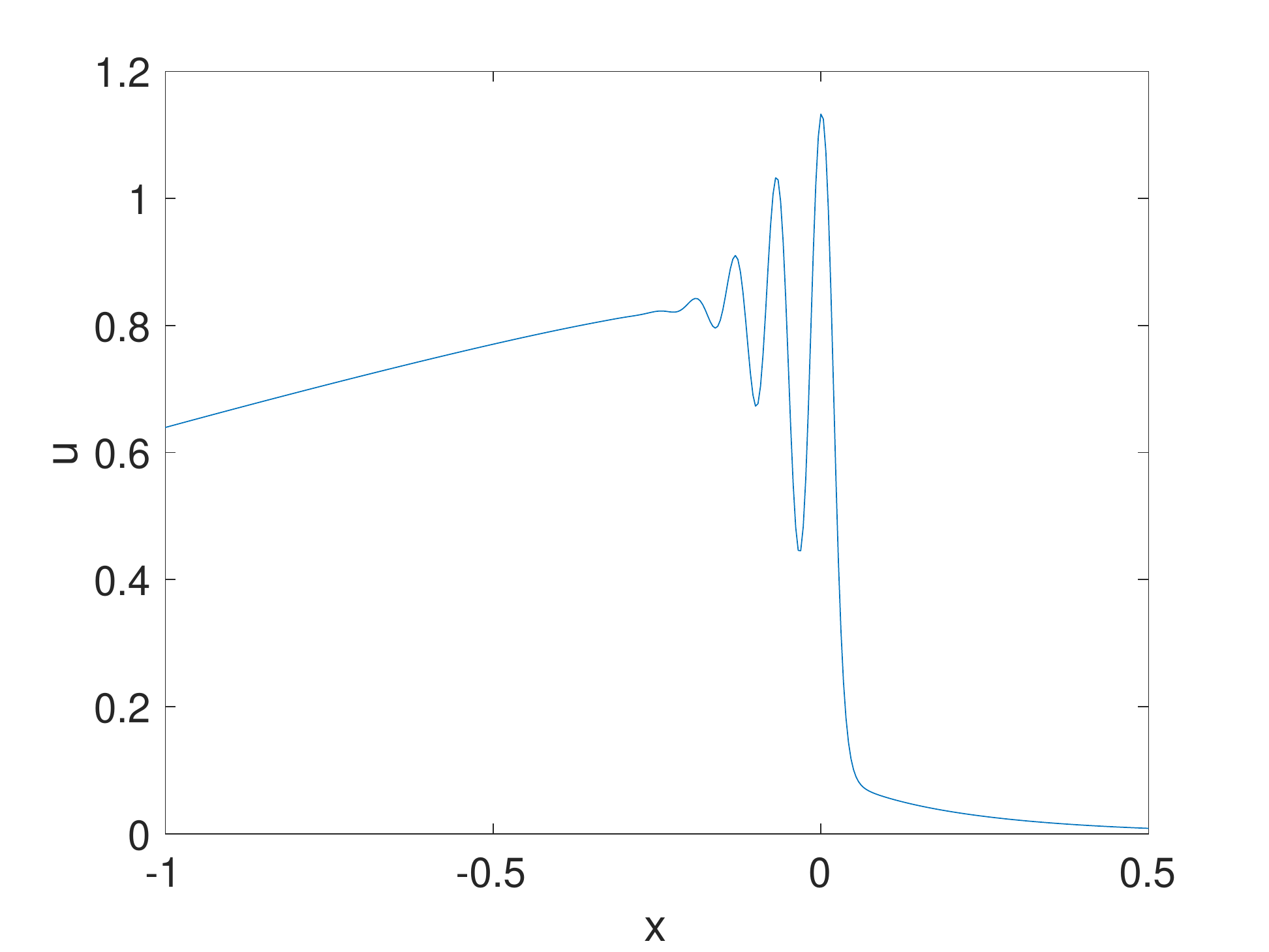}
  \includegraphics[width=0.49\textwidth]{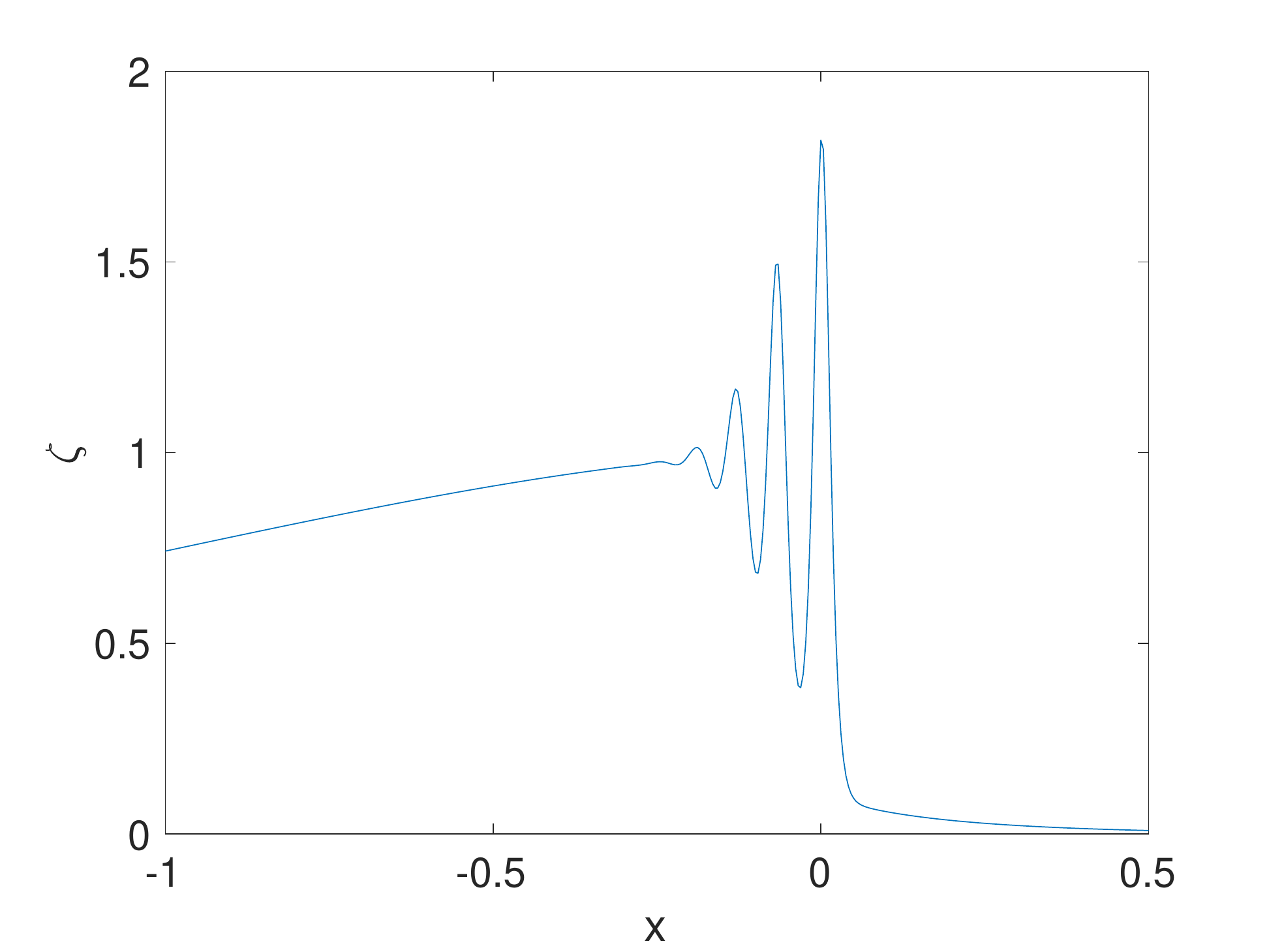}
 \caption{Solution to the SGN equations~\eqref{SGN-delta} with 
 $\delta=0.01$ for $\zeta(x,t=0)=\exp(-(x+3)^{2})$,
  $u(x,t=0) = 2\sqrt{1+\zeta(x,t=0)} - 2$ at $t=1.3$; on the left 
 $u$, on the right $\zeta$.}
 \label{gauss1e2}
\end{figure}

For smaller values of $\delta$, the modulated oscillations become more localized with 
stronger gradients. To treat the same initial data as in 
Fig.~\ref{gauss1e1} for $\delta=10^{-2}$, we use $N=2^{12}$ Fourier 
modes for $x\in2.5[-\pi,\pi]$ and $N_{t}=10^{4}$ time steps for $0\leq t\leq 
1.3$. The solution to~\eqref{SGN-delta} for $t=1.3$ can be seen in 
Fig.~\ref{gauss1e2}.

\begin{figure}[htb!]
  \includegraphics[width=0.49\textwidth]{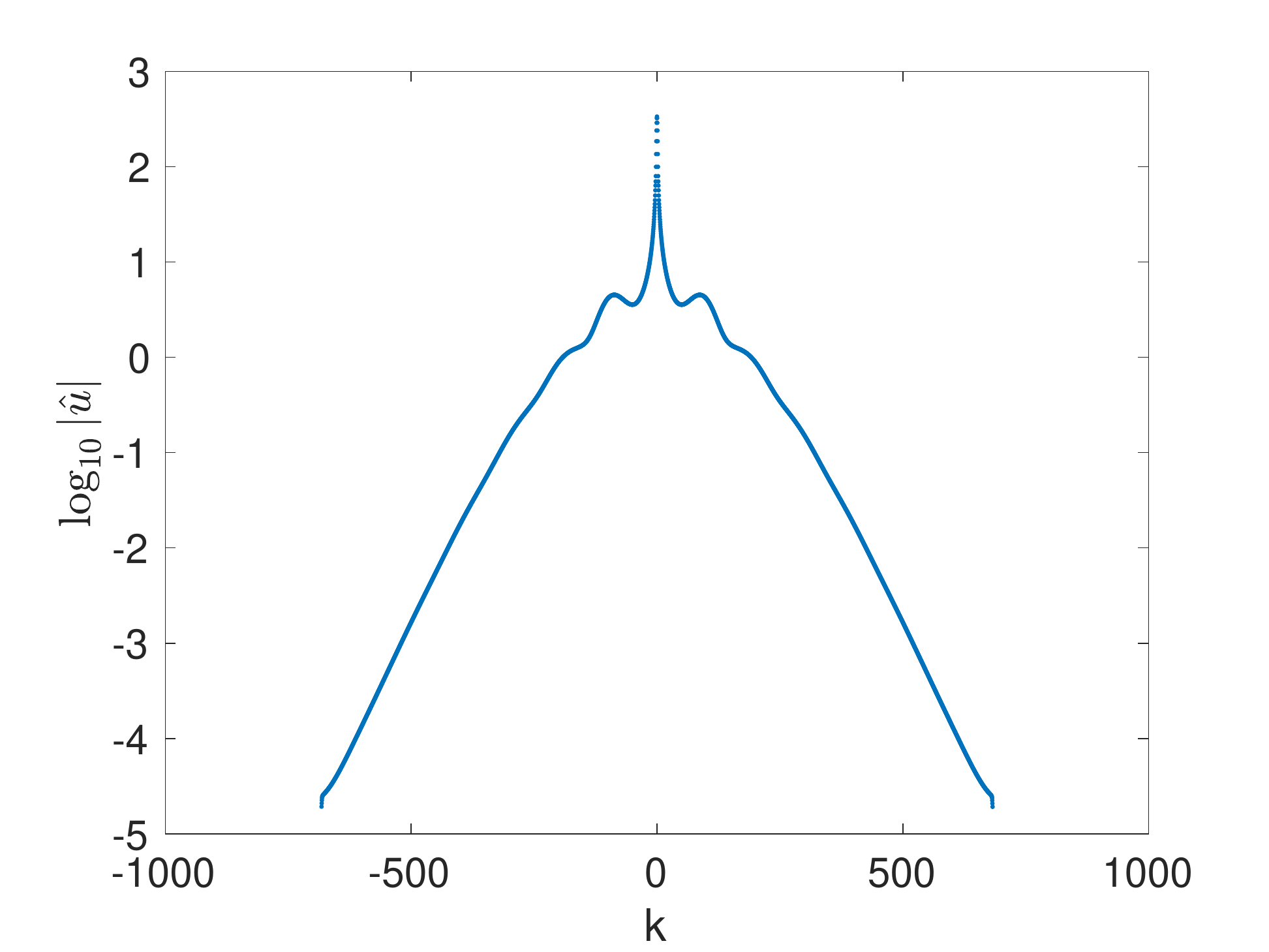}
  \includegraphics[width=0.49\textwidth]{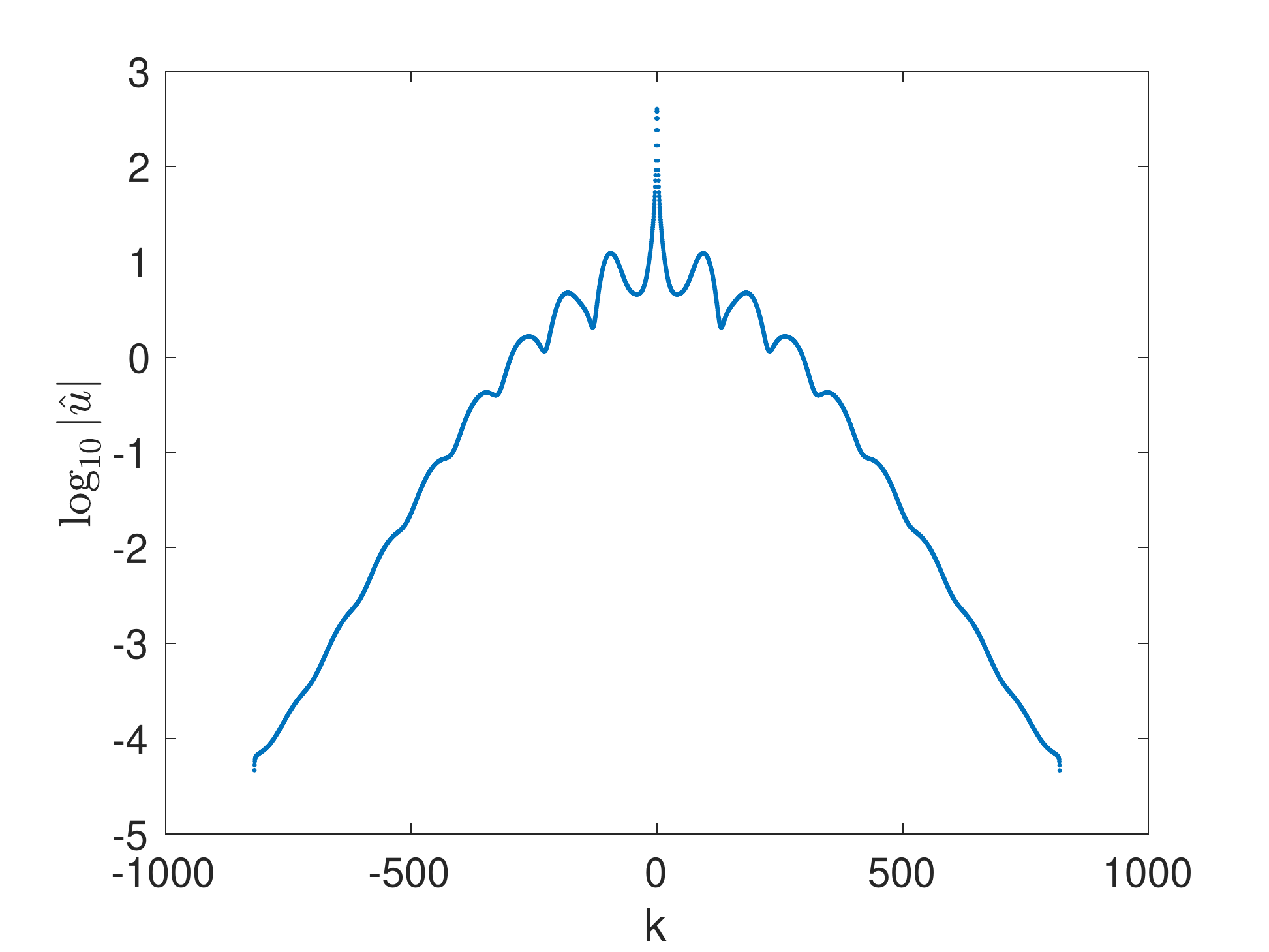}
 \caption{Fourier coefficients  of the solutions shown in Fig.~\ref{gauss1e2}.}
 \label{gauss1e2fourier}
\end{figure}

The fact that the solutions from Fig.~\ref{gauss1e2} are more 
demanding on computational resources is also clear from 
Fig.~\ref{gauss1e2fourier} where the Fourier coefficients for the 
solutions in Fig.~\ref{gauss1e2} are shown. Despite higher resolution 
the Fourier coefficients decrease to the order of $10^{-4}$ (the 
relative conservation of the third quantity in~\eqref{SGNcons} is of the order of 
$10^{-10}$) which implies a numerical error well below plotting 
accuracy.

\begin{figure}[htb!]
  \includegraphics[width=0.49\textwidth]{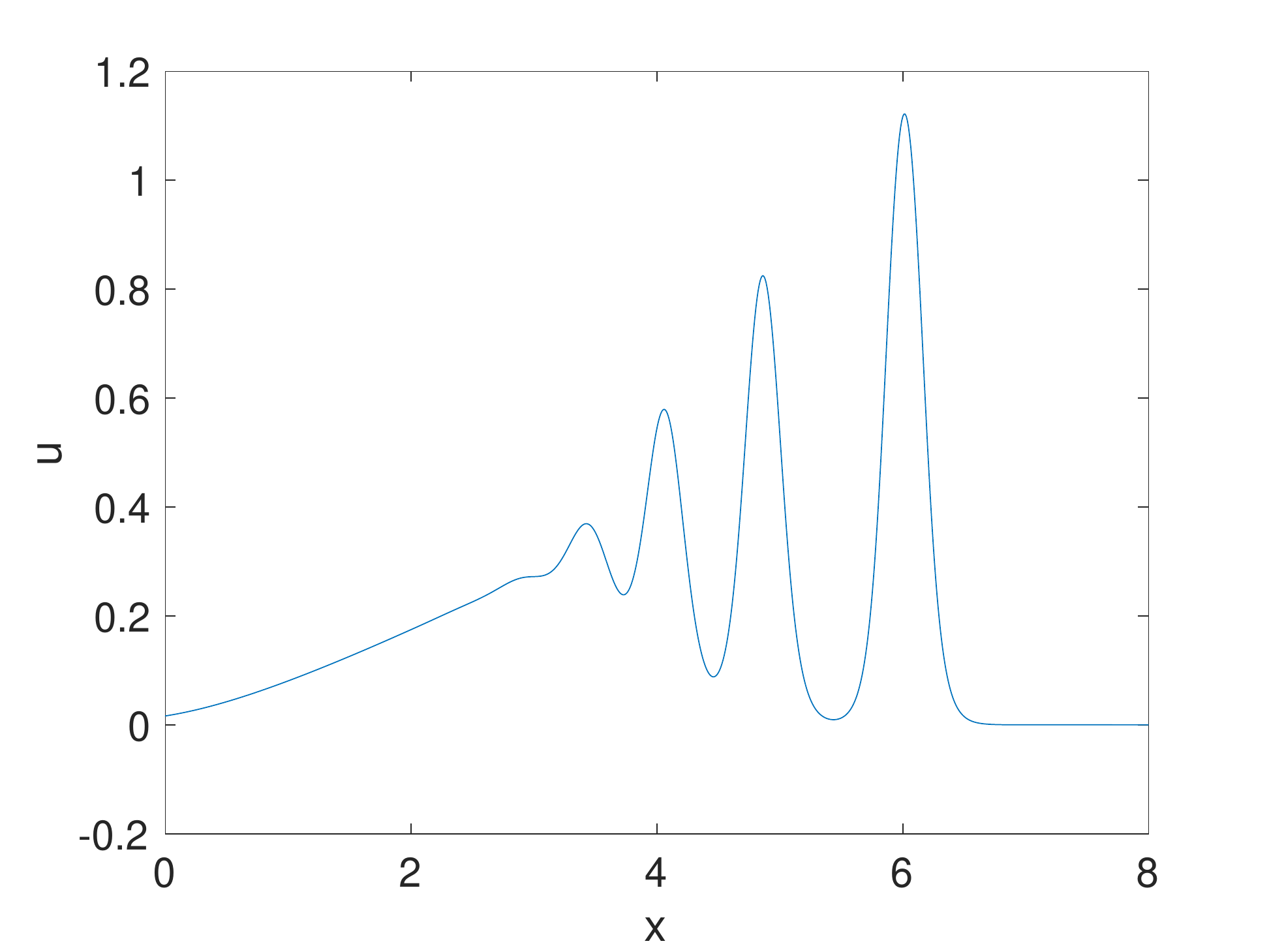}
  \includegraphics[width=0.49\textwidth]{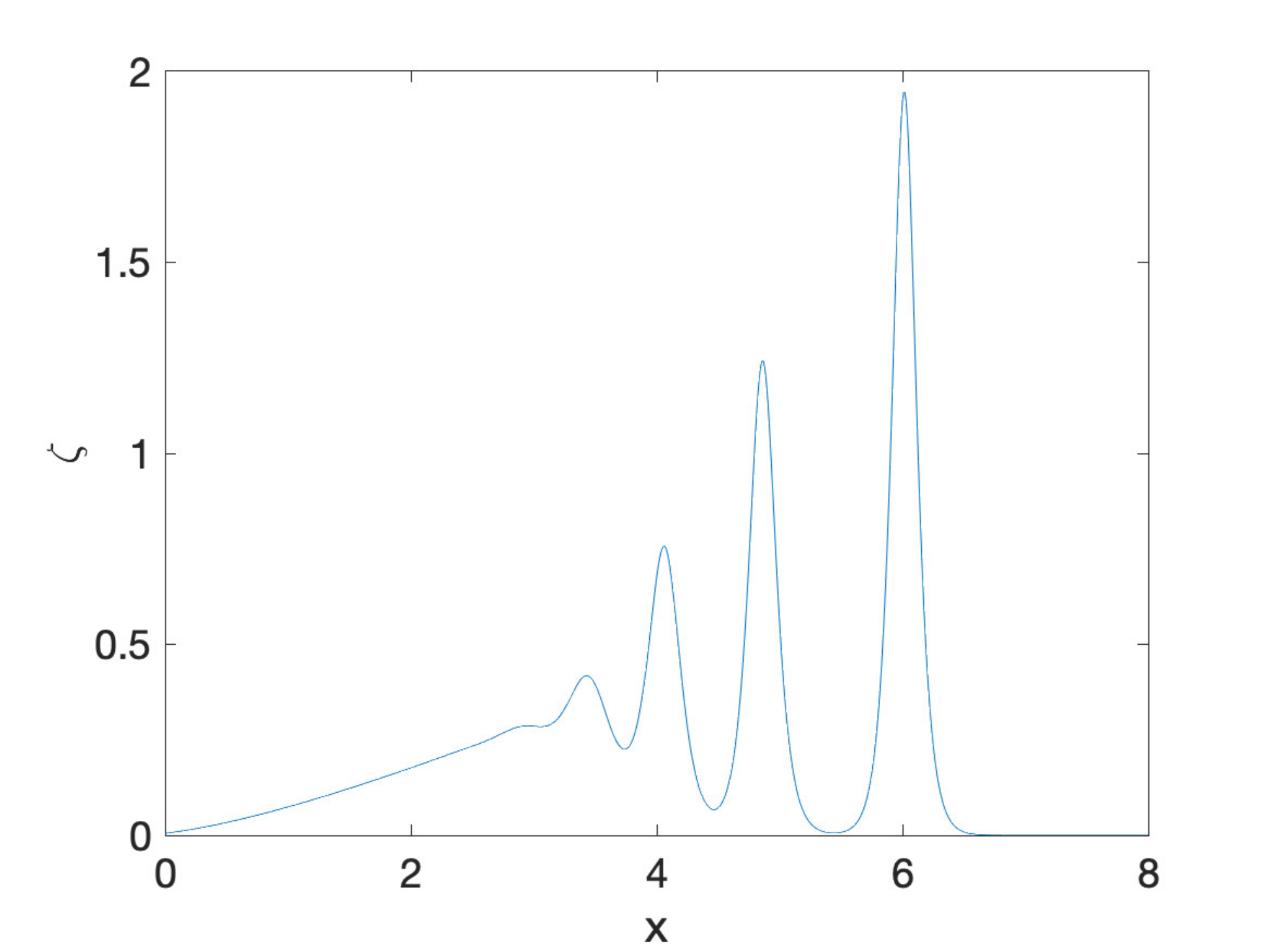}
 \caption{Solution to the WGN equations~\eqref{WGN-delta} with 
 $\delta=0.1$ for $\zeta(x,t=0)=\exp(-(x+3)^{2})$, 
 $u(x,t=0) = 2\sqrt{1+\zeta(x,t=0)} - 2$ in dependence of time; on the left 
 $u$, on the right $\zeta$.}
 \label{fdWGNgauss1e1}
\end{figure}

Using the same settings for the WGN equations~\eqref{WGN-delta} yields qualitatively
similar results. We note however that the exponential decay rate for the WGN equations is
smaller than the corresponding one for the SGN equation, consistently with the same
observation for solitary waves in the preceding section. 
Therefore to study the case of $\delta=0.1$ for the WGN equations, we use the same numerical 
parameters as for the SGN equations except for a higher number of Fourier modes ($ N=2^{11} $). 
The solution at $ t=5 $ can be seen in 
Fig.~\ref{fdWGNgauss1e1}. 

\begin{figure}[htb!]
  \includegraphics[width=0.49\textwidth]{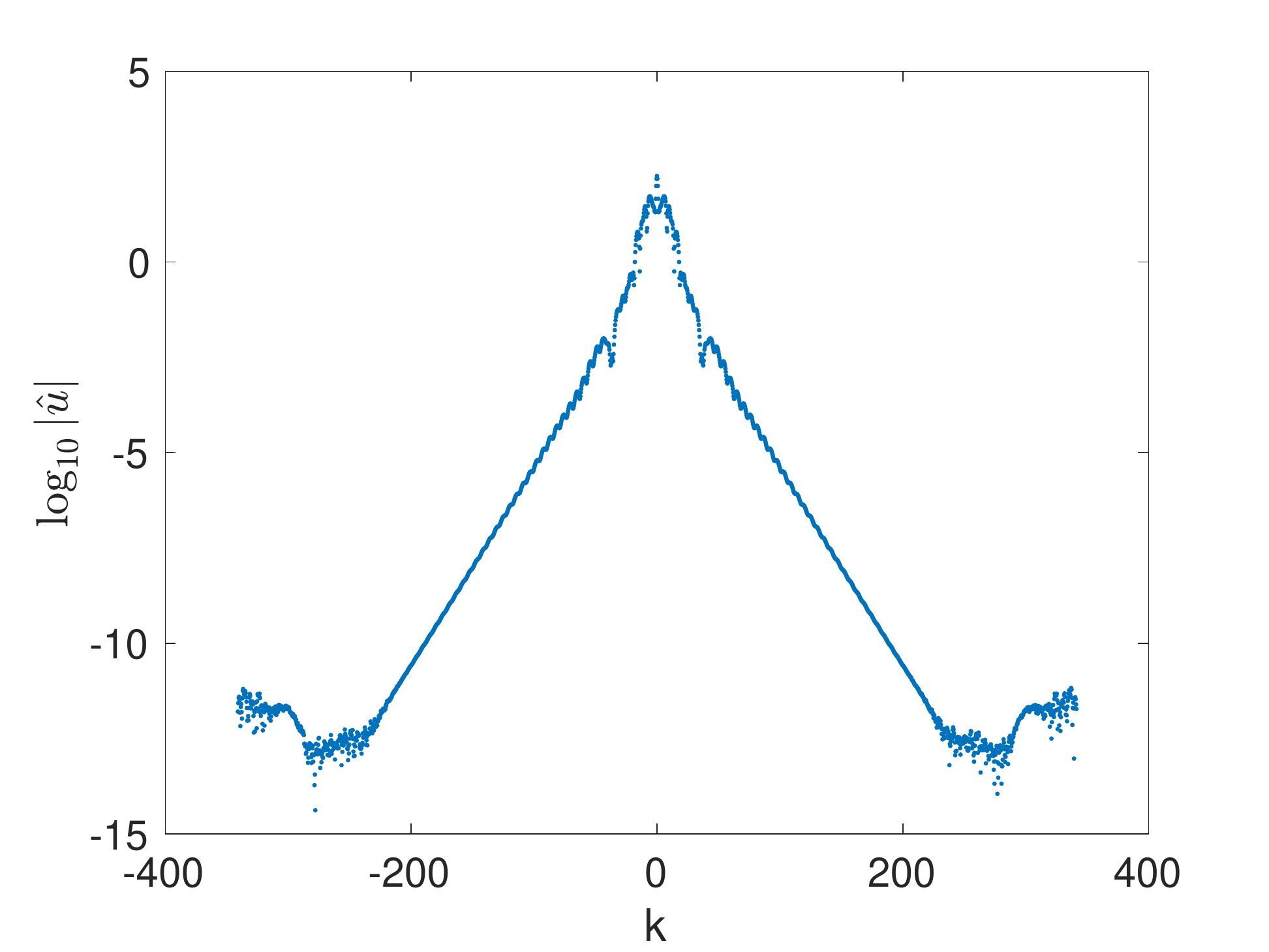}
  \includegraphics[width=0.49\textwidth]{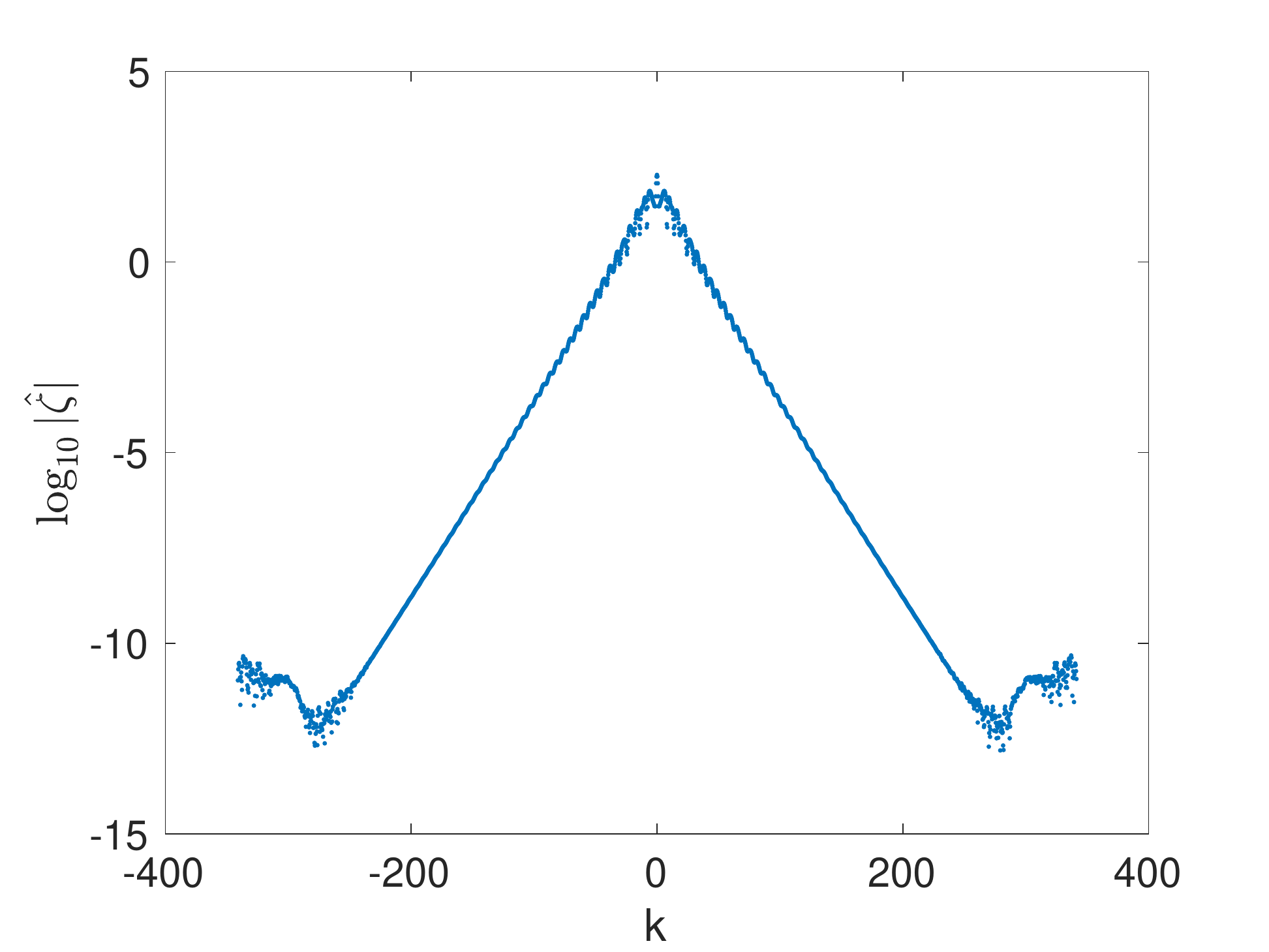}
 \caption{Fourier coefficients  of the solutions shown in Fig.~\ref{fdWGNgauss1e1}.}
 \label{fdWGNgauss1e1fourier}
\end{figure}

The behavior of Fourier coefficients as shown in 
Fig.~\ref{fdWGNgauss1e1fourier} is similar to that shown in Fig. 
\ref{gauss1e12},
although  twice the resolution in Fourier space is needed in order to achieve the same
decrease. 

\section{Comparison with the Camassa-Holm equation}\label{S.CH}

The emergence of modulated oscillations in the examples of the preceding section does 
not exclude the possibility of finite-time shock formation or other type of singularities
for solutions with different initial data.
An example of an equation with such a behavior is 
the Camassa-Holm equation~\cite{CamassaHolm93} 
\begin{equation}
  (1-\tfrac{5\delta^2}{12}\partial_{x}^2)  \partial_t w +\partial_x w+\tfrac32w\partial_x w -\tfrac{\delta^2}4 \partial_x^3 w-\tfrac{5\delta^2}{24}\big(2(\partial_x w)(\partial_x^2w) +w\partial_x^3w\big)=0.
    \label{CH}
\end{equation}
Here we chose coefficients from~\cite{ConstantinLannes09} (see (17) therein), where it is proved that~\eqref{CH} is a higher order 
(when compared with~\eqref{iB} or the KdV equation) unidirectional model for the SGN equations~\eqref{SGN-delta} (and hence the WGN equations~\eqref{WGN-delta} as well),
provided we set
\begin{equation}\label{CHtoSGN}
\left\{
\begin{array}{l}
u = w+\frac{\delta^2}{12} \partial_x^2 w+\frac{\delta^2}{6} w\partial_x^2 w,\\
\zeta = u+\frac14u^2+\frac{\delta^2}{6} \partial_t\partial_x u-\frac{\delta^2}{6} u\partial_x^2u-\frac{5\delta^2}{48} (\partial_xu)^2.
\end{array}\right.
\end{equation}
This equation generates both dispersive shock waves (see for instance~\cite{GK08,AGK}) and 
finite-time singularities in the form of {\em surging wavebreaking}~\cite{ConstantinEscher98}.
Let us finally mention that~\eqref{CH} and~\eqref{SGN-delta} are of the same kind,
namely quasilinear nonlocal dispersive equations involving only differential operators.

The set of initial data for which solutions to~\eqref{CH} lead to finite-time wavebreaking contain
any smooth and odd function $w_0=w(\cdot,t=0)$ such that $w_0'(0)<0$ and $w_0(x)<0$ for $x>0$~\cite{McKean04};
we choose here
\[ w_0(x) = -x\exp(-x^{2}).\]
Inferring initial data for $\zeta,u$ through~\eqref{CHtoSGN} would require to know $\partial_t\partial_x u(\cdot,t=0)$, however we can approximately solve the equations with a harmless $O(\delta^4)$ approximation by setting
\begin{equation}\label{CHtoSGN2}
\left\{
\begin{array}{l}
u(\cdot,t=0) =:u_0= w_0+\frac{\delta^2}{12} \partial_x^2 w_0+\frac{\delta^2}{6} w_0\partial_x^2 w_0\\
\zeta(\cdot,t=0) = u_0+\frac14u_0^2-\frac{\delta^2}{6} \partial_x^2 ( u_0+\tfrac34u_0^2)-\frac{\delta^2}{6} u_0\partial_x^2u_0-\frac{5\delta^2}{48} (\partial_xu_0)^2
\end{array}\right.
\end{equation}
where we used ~\eqref{CH} to infer that
\[\partial_t u = \partial_t w +O(\delta^2)=- (\partial_x w+\tfrac32w\partial_x w) +O(\delta^2)=- (\partial_x u+\tfrac32u\partial_x u) +O(\delta^2).\]

\begin{figure}[htb!]
  \includegraphics[width=0.49\textwidth]{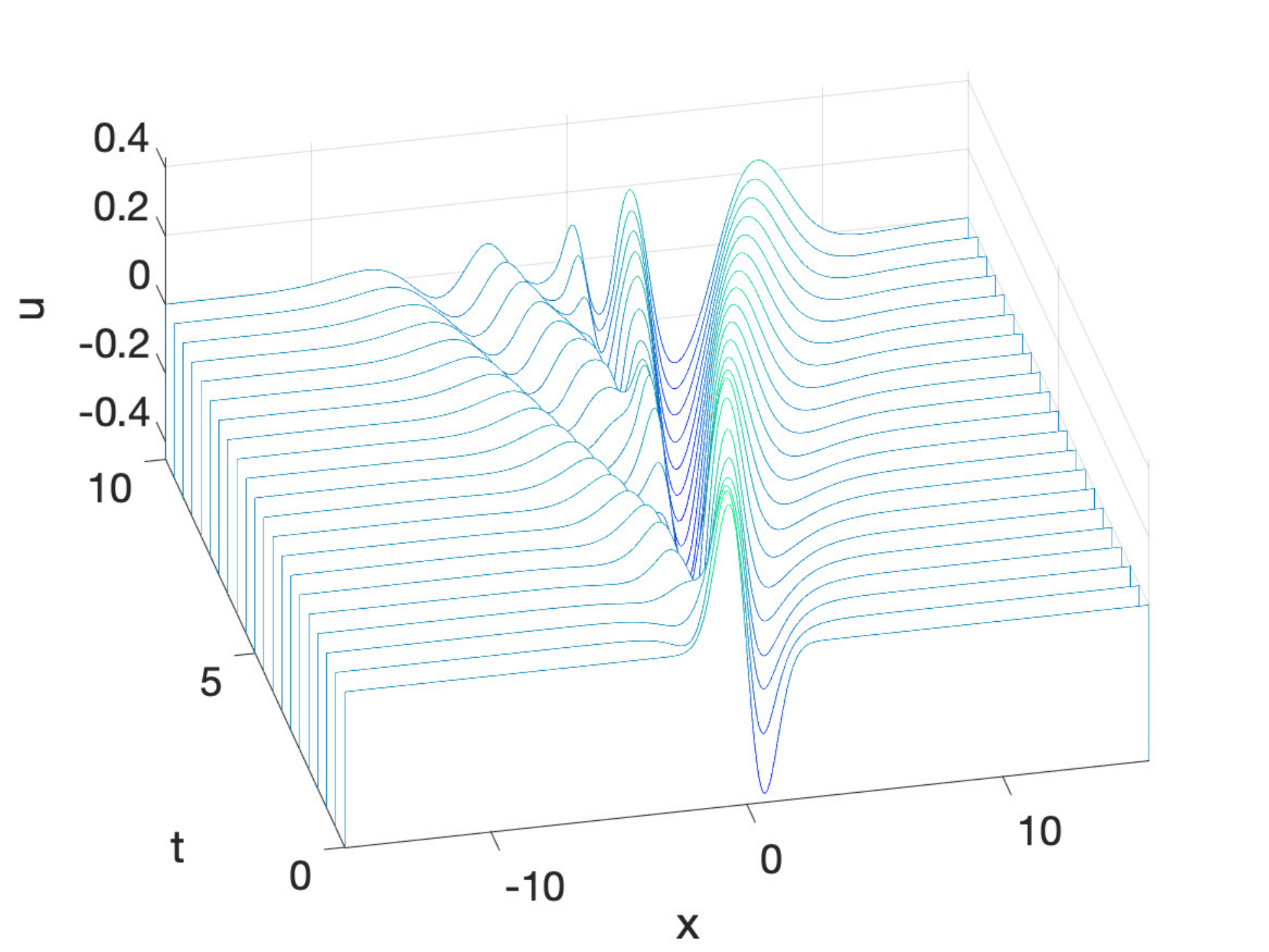}
  \includegraphics[width=0.49\textwidth]{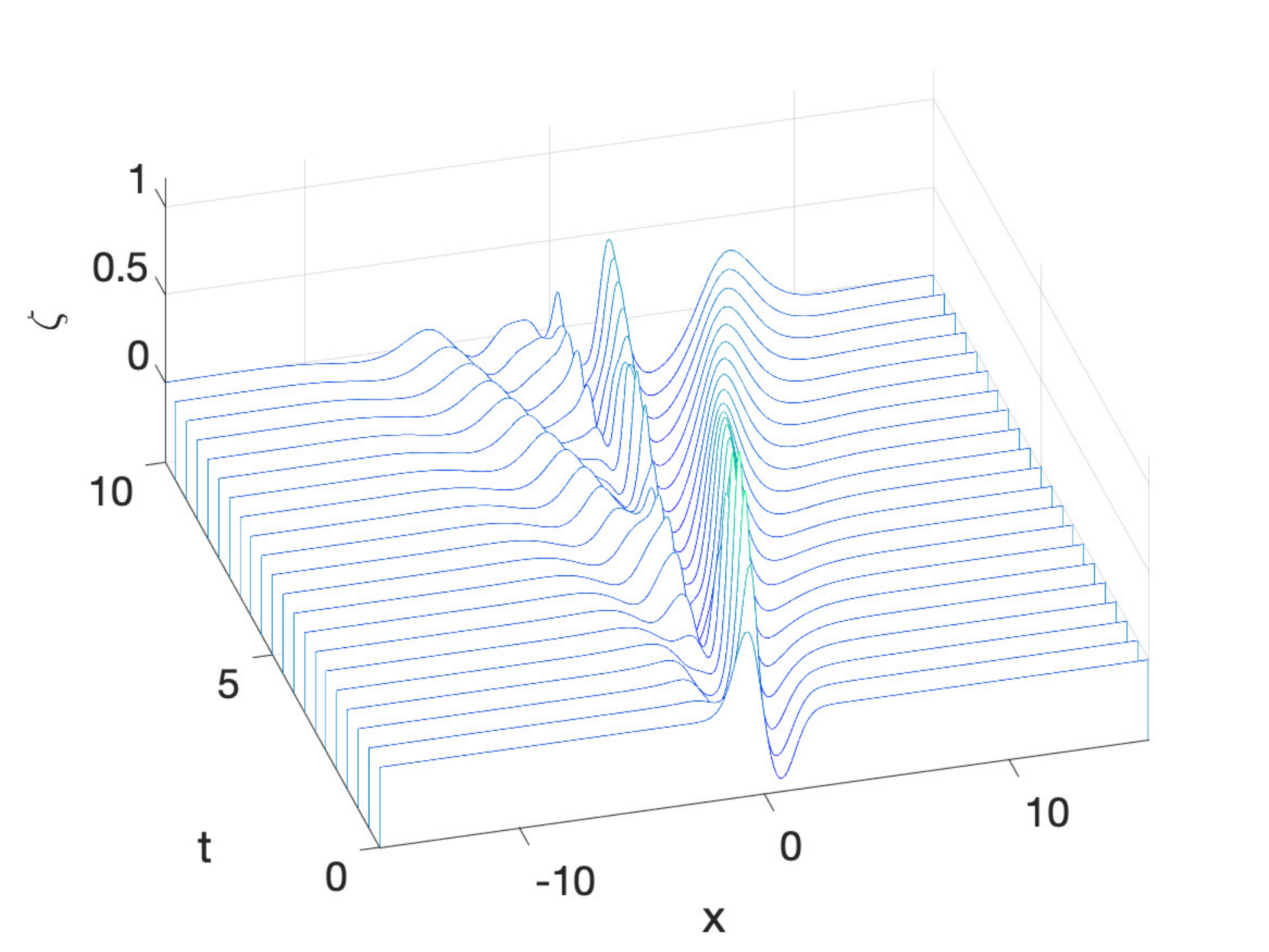}
 \caption{Solution to the SGN equations for the initial data 
\eqref{CHtoSGN2} with  $\delta = 0.1$.}
 \label{SGNCH}
\end{figure}

We consider the example $\delta^2 = 0.1$.
 For the computation we use $N=2^{11}$ Fourier modes 
for $x\in5[-\pi,\pi]$ and $N_{t}=10^{4}$ time steps for $t\leq 10$. 
The SGN solution for the initial data~\eqref{CHtoSGN2} can 
be seen in Fig.~\ref{SGNCH}. Obviously strong gradients 
appear, but these do not lead to a shock formation.

\begin{figure}[htb!]
  \includegraphics[width=0.49\textwidth]{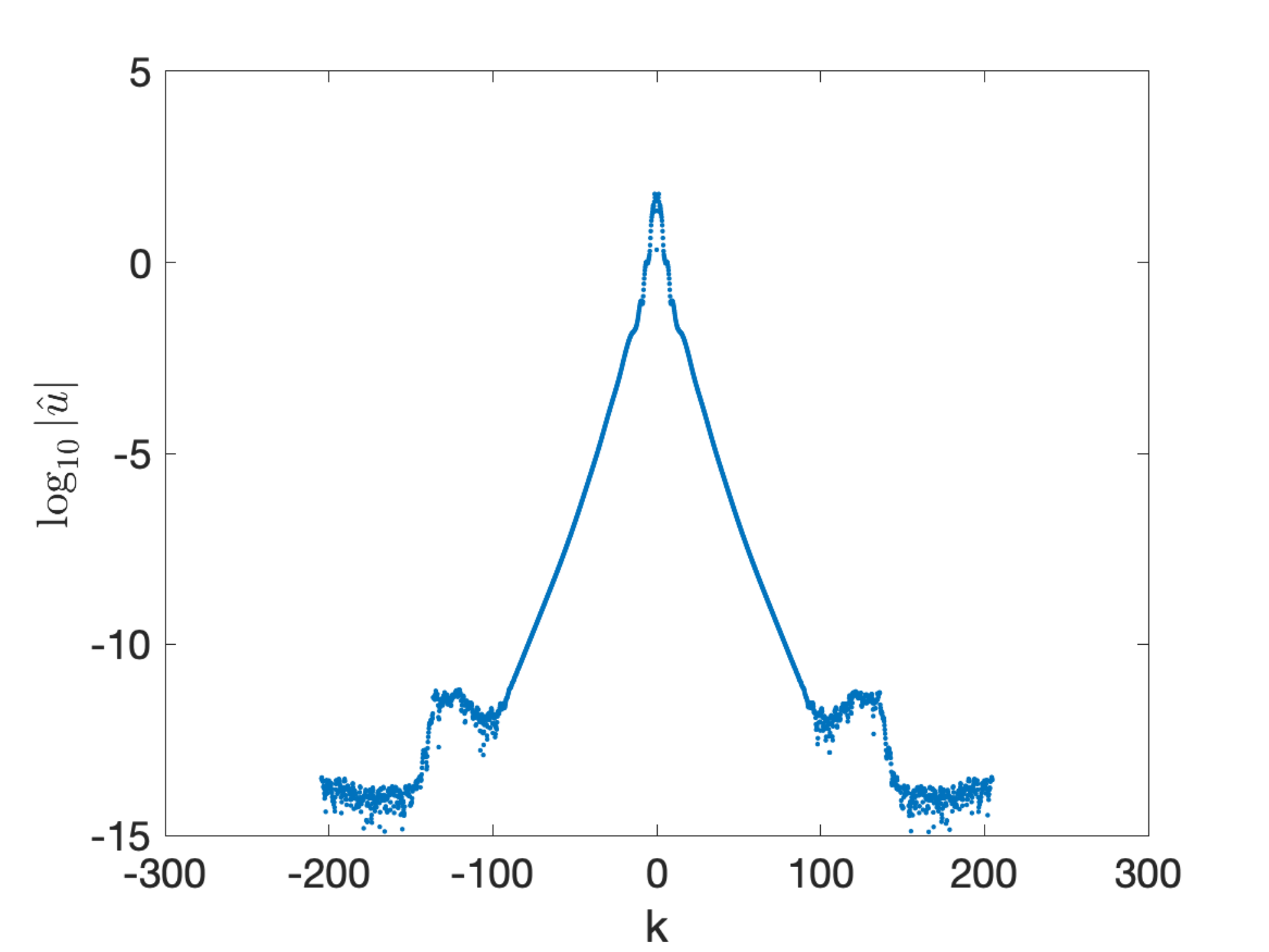}
  \includegraphics[width=0.49\textwidth]{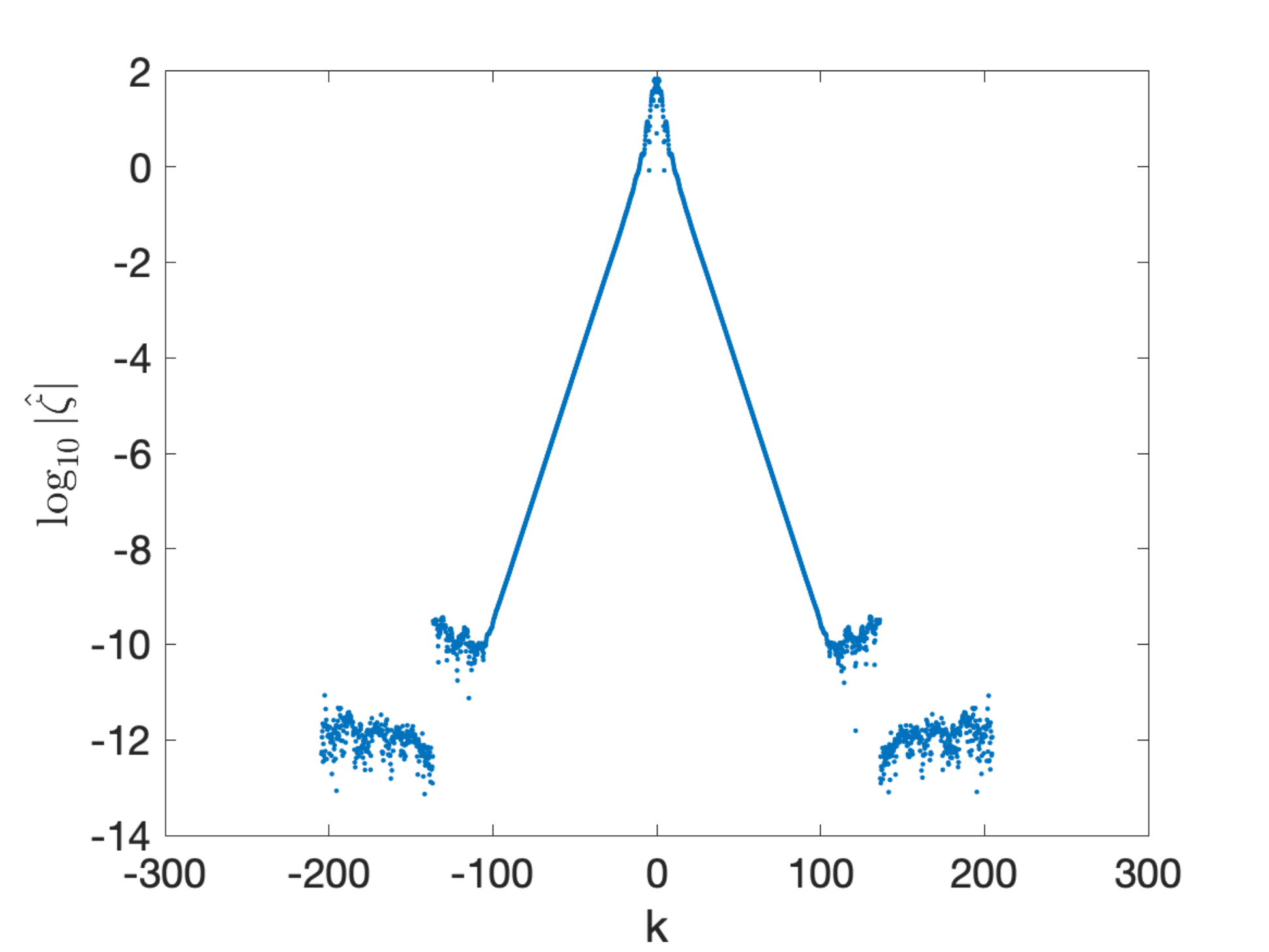}
 \caption{Fourier coefficients  of the solutions shown in 
 Fig.~\ref{SGNCH} for $t=10$.}
 \label{SGNCHfourier}
\end{figure}

This example is numerically challenging because of convergence problems 
for the GMRES algorithm which imply 
a pollution of the Fourier coefficients at high wave numbers. 
Therefore we use a dealiasing according to the 2/3-rule which means 
that the Fourier coefficients corresponding to the highest one third 
of the wave numbers are put equal to 0.  If this is 
done, the example is well resolved in space as can 
be seen  in Fig.~\ref{SGNCHfourier} where the Fourier 
coefficients of the solution in Fig.~\ref{SGNCH} are shown for $t=10$, and in time as inferred by the relative conservation of the third quantity in~\eqref{SGNcons} to the order of $10^{-12}$. Once again, using the WGN equations~\eqref{WGN-delta} instead of~\eqref{SGN-delta} does not modify substantially the behavior of the solution, and we do not show the corresponding pictures.

\section{Near-cavitation initial data}\label{S.cavitation}

In Sections~\ref{S.DSW} and~\ref{S.CH}, we used initial data leading to shocks for simplified models, and observed for each scenario that the solutions to both the \ref{SGN} and \ref{WGN} equations remained smooth. This leaves open the important question of global-in-time well-posedness of the equations. 
Let us recall that the local well-posedness of the \ref{SGN} equations has been proved in~\cite{Li06}\footnote{very recently improved in \cite{Inci} to allow less regular initial data.} ---and later on in more general frameworks in~\cite{Alvarez-SamaniegoLannes08a,Israwi11,FujiwaraIguchi15,DucheneIsrawi18}--- but that the technique in these works (essentially energy methods) does not provide any global-in-time result, 
in particular due to the fact that the functional setting is not controlled by conserved functionals listed in~\eqref{SGNcons}.
 It should be mentioned that the Boussinesq system obtained when neglecting the nonlinear dispersive terms in~\eqref{SGN} is known to be globally well-posed provided that the non-cavitation assumption $\inf_\RR (1+\zeta) >0$ is initially satisfied, by~\cite{Schonbek81,Amick84}; see also~\cite{MolinetTalhoukZaiter20}. This result does not generalize easily to the \ref{SGN} or the \ref{WGN} equations. Very recently, Bae and Granero-Belinch\'{o}n showed in~\cite{BaeGranero-Belinchon} that if the non-cavitation assumption initially fails to hold at one single point and some symmetry assumptions are enforced, then solutions to~\eqref{SGN} (or rather an equivalent reformulation when the non-cavitation assumption holds) preserve these assumptions for positive times and cannot remain smooth globally in time.

Motivated by this result, we consider the initial data
\begin{equation}
\zeta(x,t=0)=-0.9 \exp(-x^2), \qquad   u(x,t=0)=-x\exp(-x^2)	
	\label{cavinitial}.
\end{equation}
The non-cavitation assumption is valid initially and hence will remain valid as long as the solution remains smooth. Indeed we can define, given any $x\in\RR$ and $t>0$, $h_{x,t}(\tau)=1+\zeta(X_{x,t}(\tau),\tau)$ where $X_{x,t}(\tau)$ is the backward characteristic defined by $X_{x,t}'(\tau)=u(X_{x,t}(\tau),\tau)$ and $X_{x,t}(t)=x$. One then observes that, by the conservation of mass equation, one has for any $\tau\in[0,t]$, $h_{x,t}'(\tau)=- h_{x,t}(\tau)(\partial_xu)(X_{x,t}(\tau),\tau) $, and hence $1+\zeta(x,t)=\big(1+\zeta(X_{x,t}(0),t=0)\big)\exp\big(\int_0^t- (\partial_xu)(X_{x,t}(\tau),\tau)\, \dd\tau\big)>0$.
As we see in the numerical results below, such initial data produce very steep gradients,
and a possible blowup scenario which deserves to be investigated in more details.

\begin{figure}[htb!]
  \includegraphics[width=0.49\textwidth]{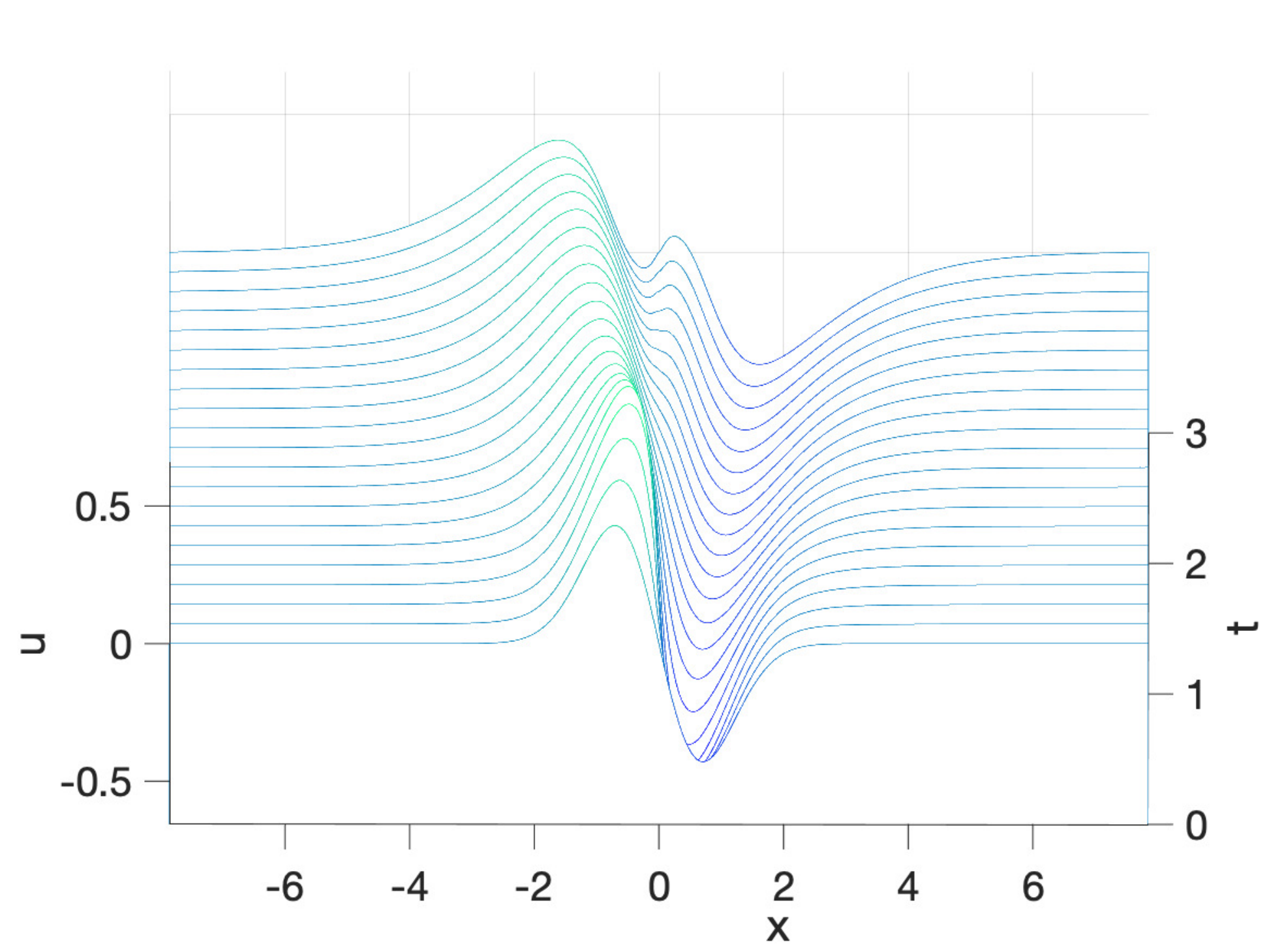}
  \includegraphics[width=0.49\textwidth]{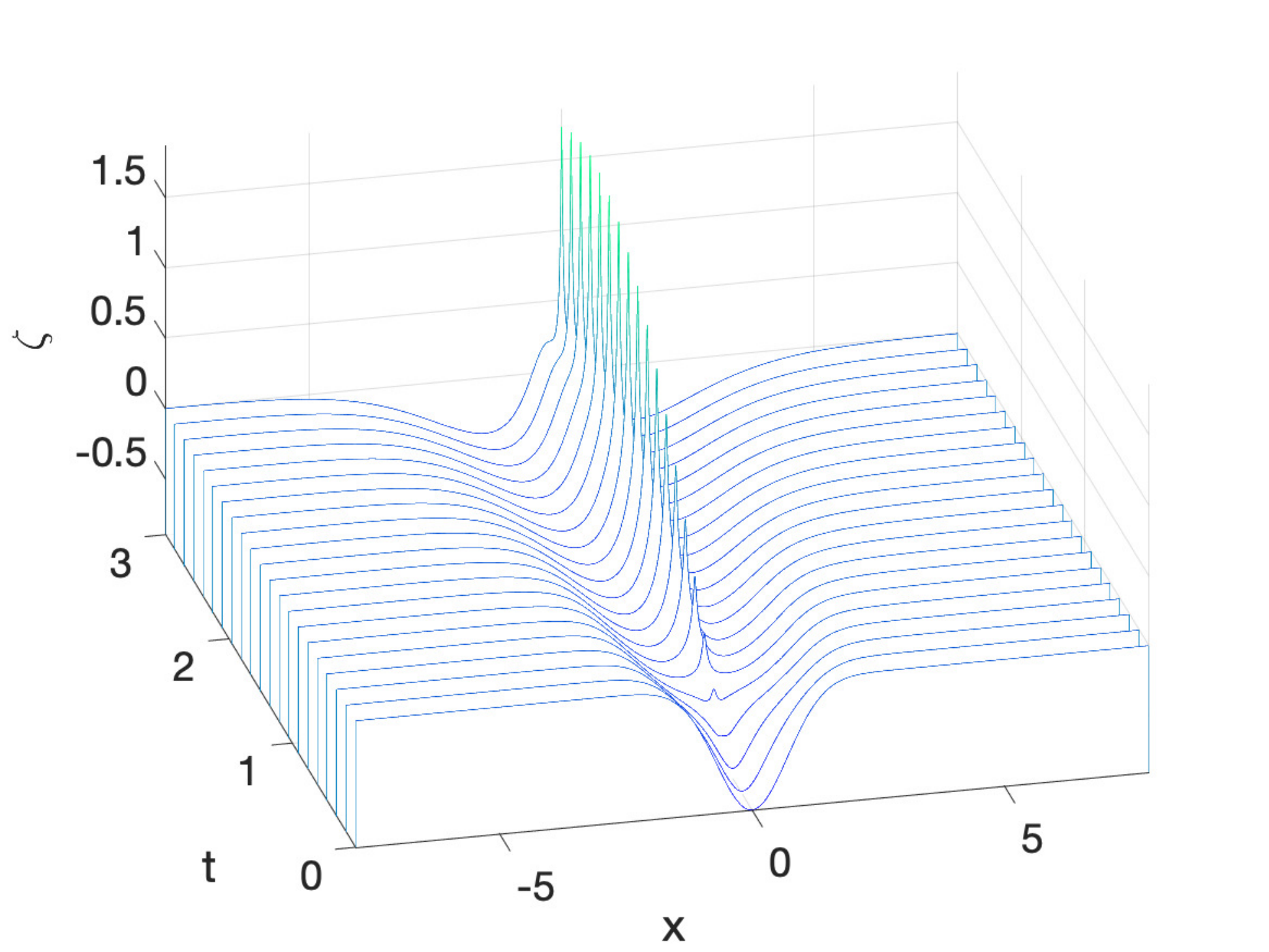}
 \caption{Solution to the \ref{SGN} equations for the initial data 
 (\ref{cavinitial}), on the left $u$, on the right $\zeta$.}
 \label{SGNcavwater}
\end{figure}

To address this question we use $N=2^{12}$ Fourier modes for $x\in 
2.5[-\pi,\pi]$ with dealiasing and 
$N_{t}=10^{4}$ time steps for $t\leq 3$ to solve the \ref{SGN} equations for the initial 
data (\ref{cavinitial}). The solution can be seen in 
Fig.~\ref{SGNcavwater}. 

\begin{figure}[htb!]
  \includegraphics[width=0.49\textwidth]{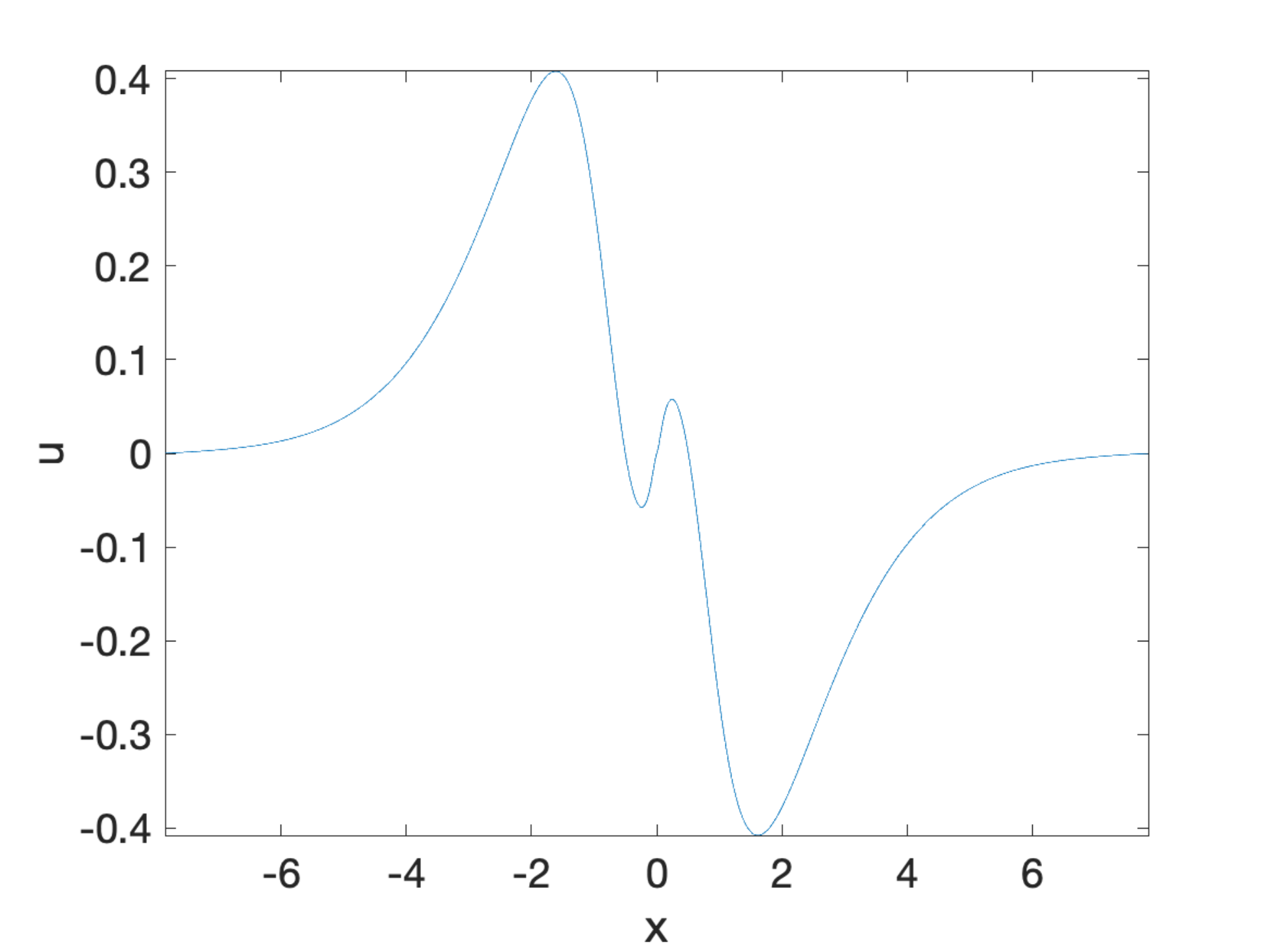}
  \includegraphics[width=0.49\textwidth]{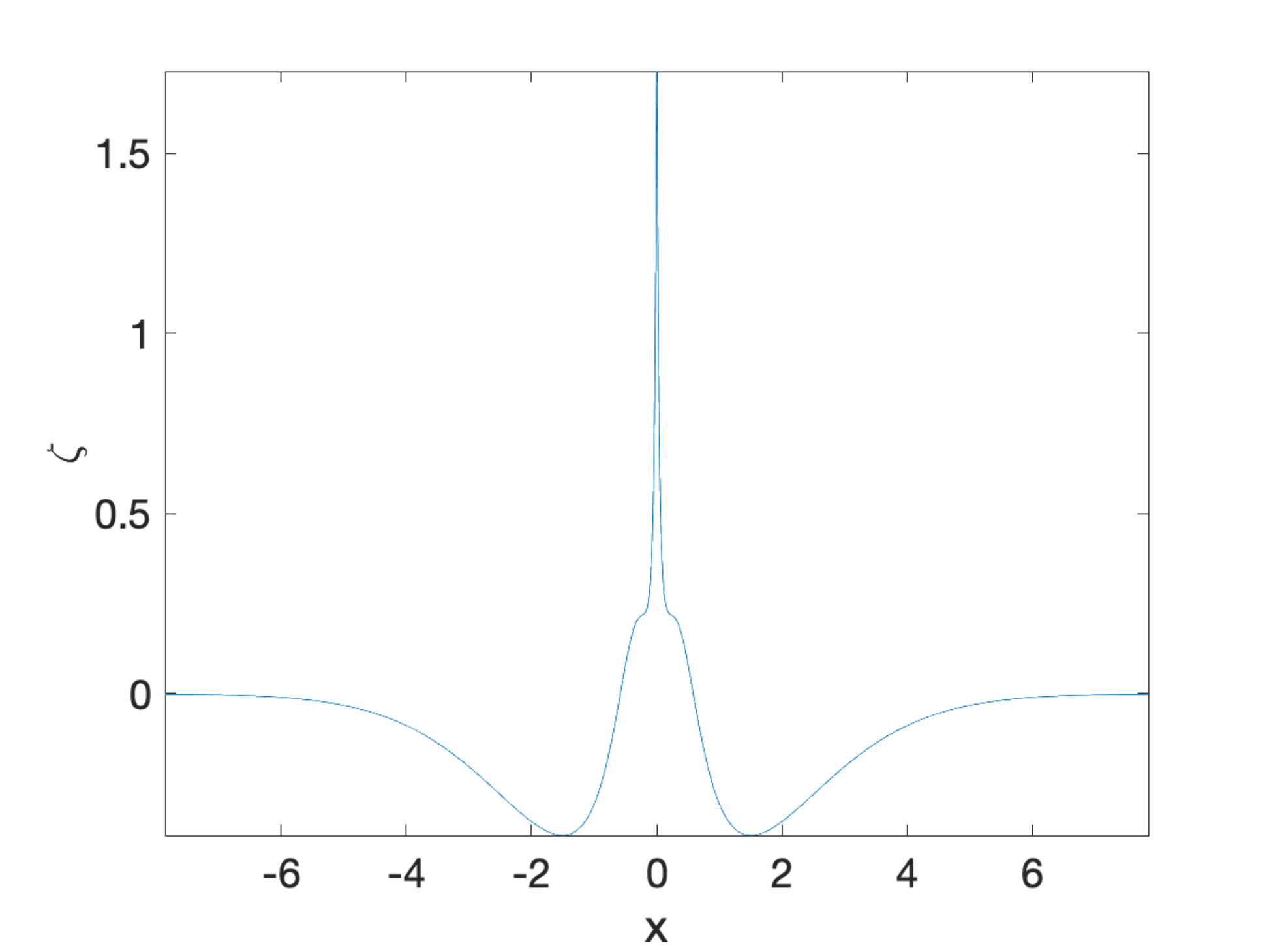}
 \caption{Solution to the \ref{SGN} equations for the initial data 
~\eqref{cavinitial} for $t=3$, on the left $u$, on the right $\zeta$.}
 \label{SGNcavt3}
\end{figure}

The function $\zeta$ 
develops some cusp-like structure which is strongly peaked. The 
solution at the final time is shown in Fig.~\ref{SGNcavt3}. The 
function $u$ appears to stay smooth. 

\begin{figure}[htb!]
  \includegraphics[width=0.49\textwidth]{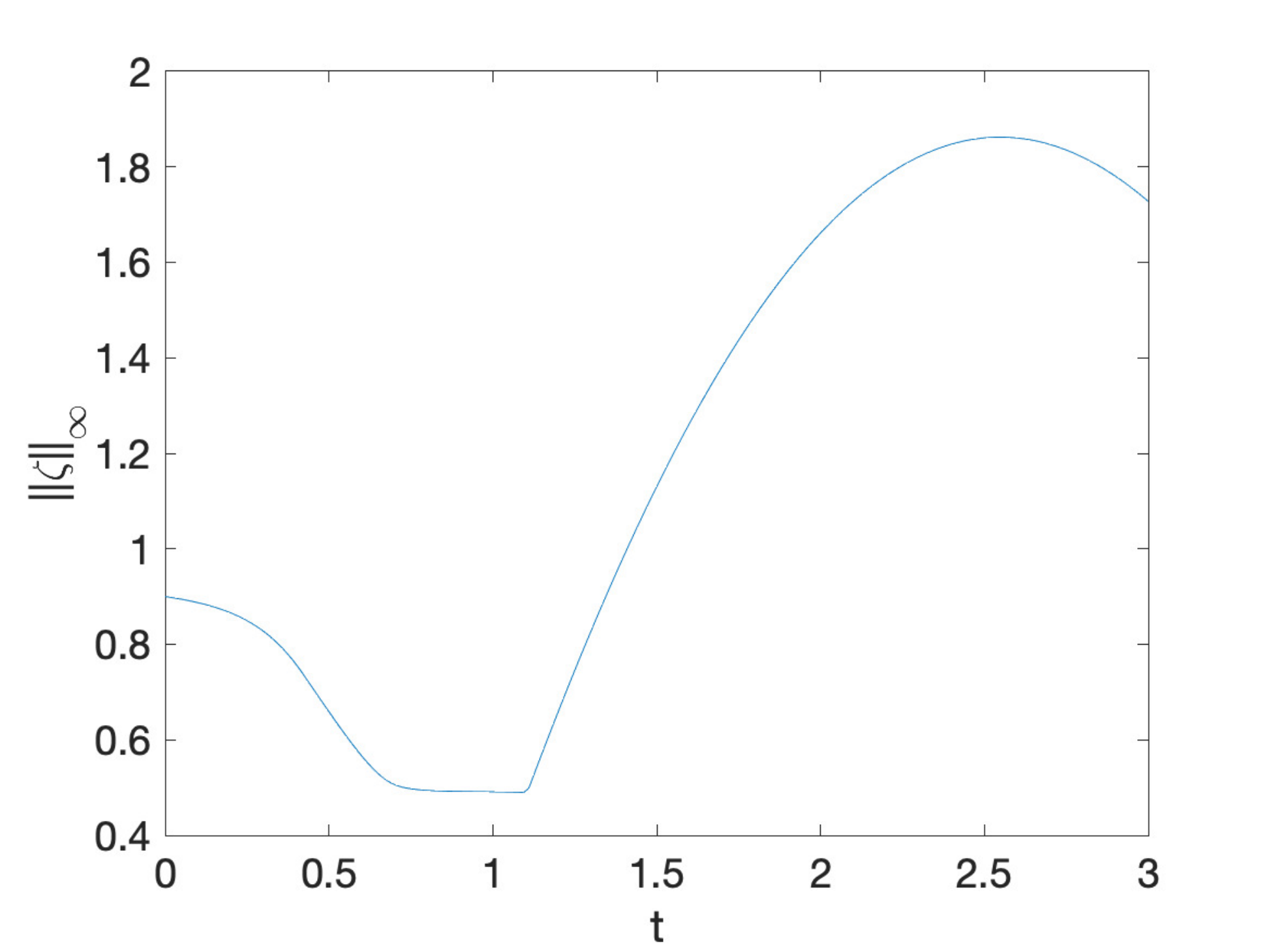}
  \includegraphics[width=0.49\textwidth]{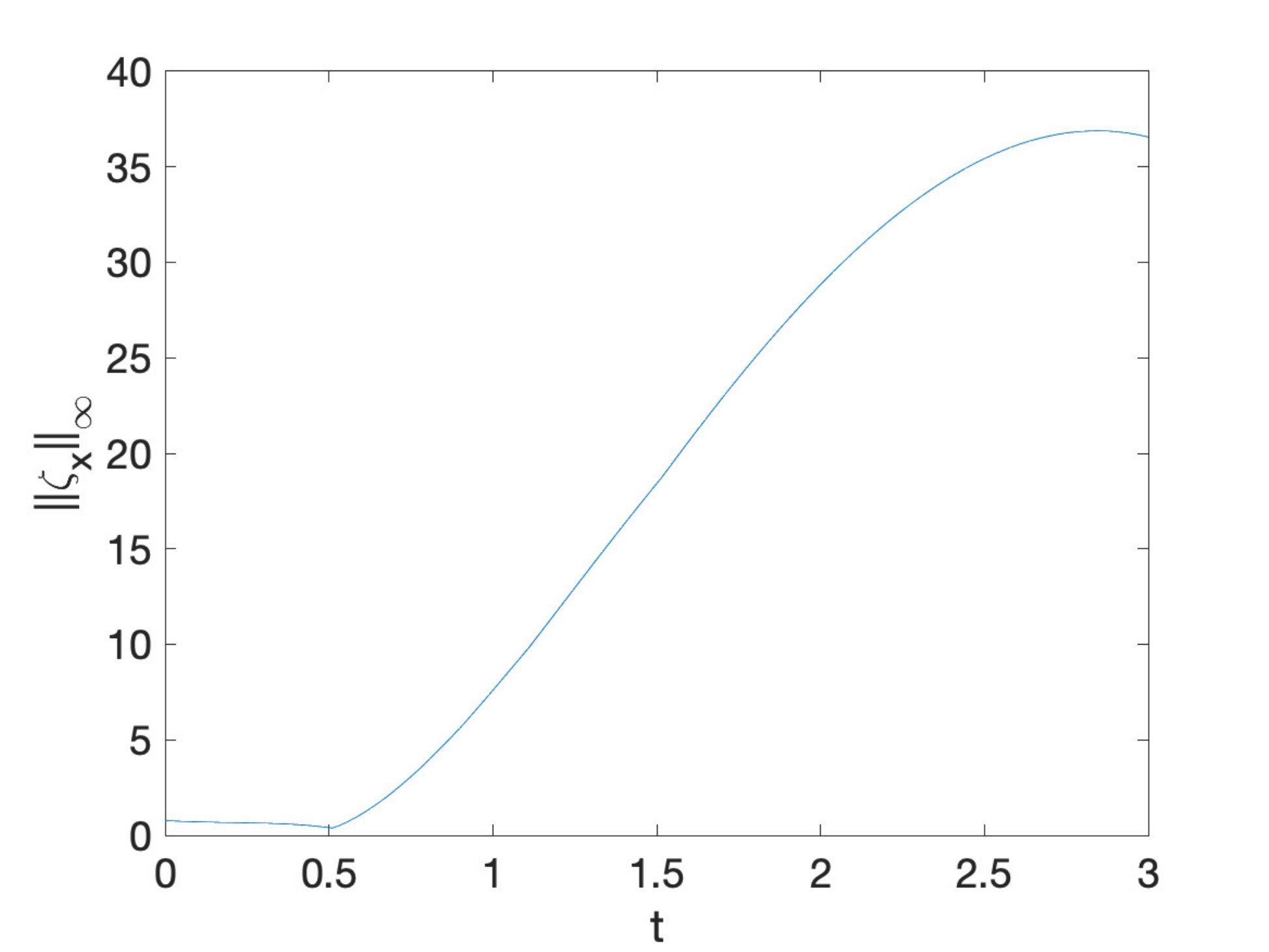}
 \caption{$L^{\infty}$ norms of the solution $\zeta$ to the \ref{SGN} equations for the initial data 
~\eqref{cavinitial} on the left and for its gradient  on the right.}
 \label{SGNcavnorms}
\end{figure}

To decide whether a blow-up is possible in this case, we show the 
$L^{\infty}$ norms of both $\zeta$ and $\zeta_{x}$ in Fig.~\ref{SGNcavnorms} in 
dependence of time. Whereas the $L^{\infty}$ norm of $\zeta$ grows 
for some time, it reaches a maximum for $t\sim 2.5$ and decreases 
then. Thus there is no $L^{\infty}$ blow-up, but the strong gradient 
can be seen in the same figure on the right. But also the 
$L^{\infty}$ norm of the gradient appears to reach a finite maximum. 
This would indicate that one is close to a cusp-like situation, but 
that the solution stays smooth in this example. 

\begin{figure}[htb!]
  \includegraphics[width=0.49\textwidth]{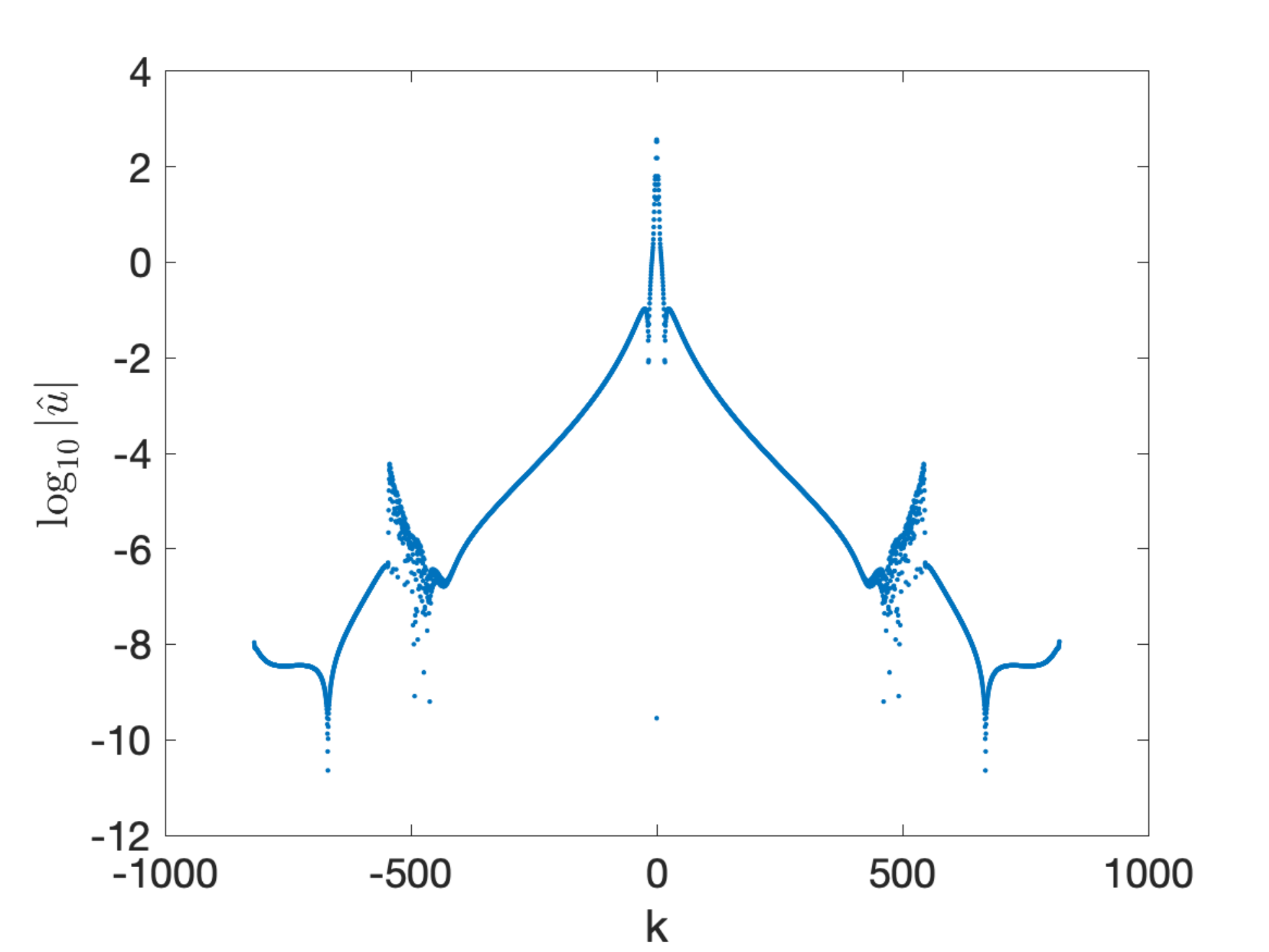}
  \includegraphics[width=0.49\textwidth]{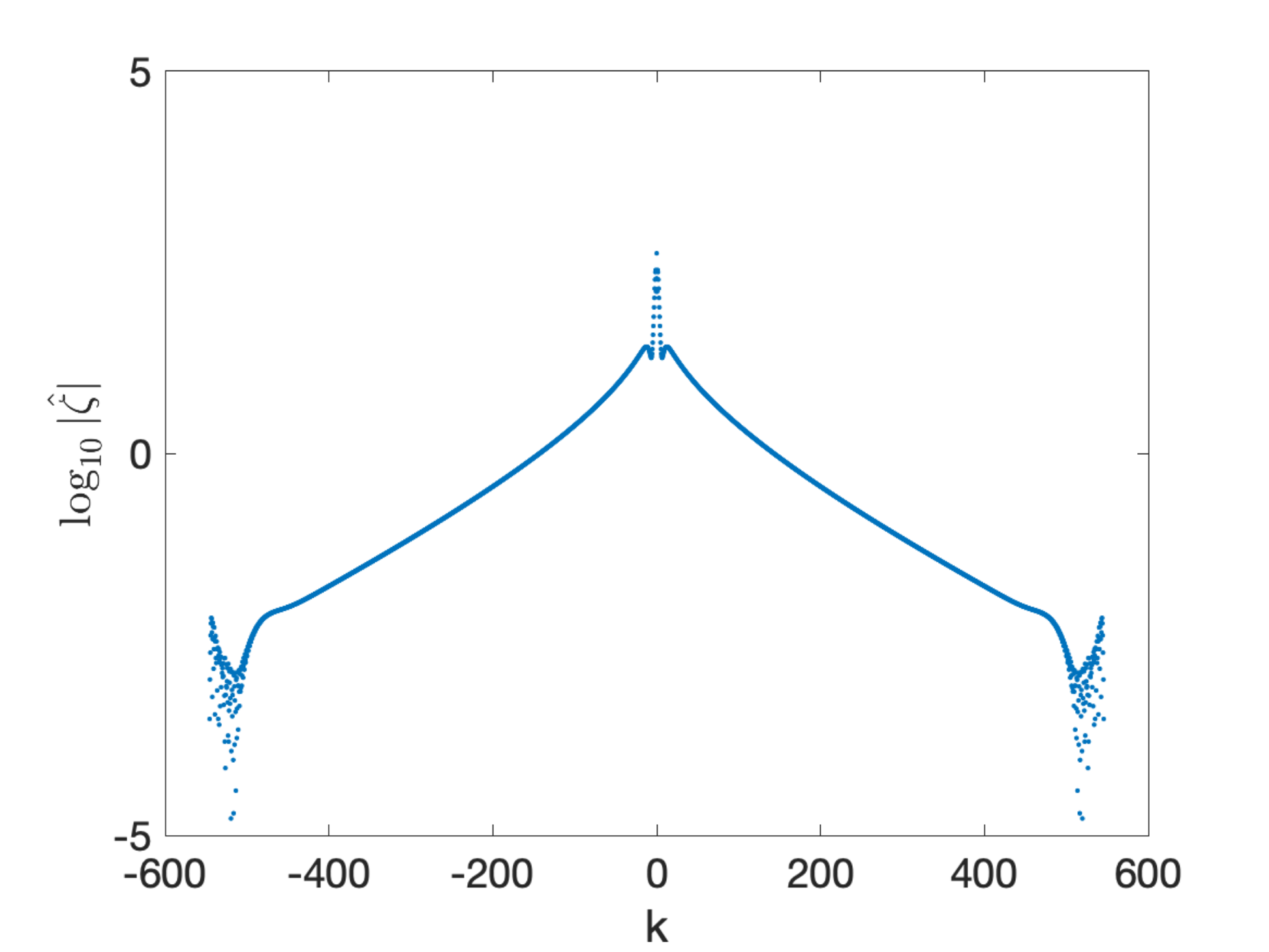}
 \caption{The Fourier coefficients of the solution  to the \ref{SGN} equations for the initial data 
 (\ref{cavinitial}), on the left $u$,   on the right for $\zeta$.}
 \label{SGNcavfourier}
\end{figure}

Note that the solution is well resolved: the resolution in time as 
indicated by the relative conservation of the conserved quantities is 
of the order of $10^{-9}$, and the Fourier coefficients of the 
solution are shown in Fig.~\ref{SGNcavfourier}. 
The above 
results do not change within numerical precision if the computation 
is repeated with $N=2^{14}$ Fourier modes.

\section{Outlook}
\label{S.conc}

In this paper we have investigated numerically several aspects 
of solutions to Serre-Green-Naghdi type equations. First we have 
obtained a family of supercritical solitary waves of the fully 
dispersive system with no upper bound on 
the admissible velocities, alike the explicit family of solitary waves of the
original Serre-Green-Naghdi system. Investigating the dynamic stability of 
these solitary waves, we have found no sign of instability, even for
large velocities. We have also set up several experiments for which
solutions to either the original or  fully dispersive Serre-Green-Naghdi system 
develop zones of rapid modulated oscillations and/or steep gradients, 
but we have never monitored finite-time singularity formations.

On the numerical side, we have shown that an approach based on 
a Fourier spectral method combined with the Krylov subspace iterative 
method GMRES is very efficient and allows to study with high accuracy 
 computationally demanding problems.
%The same approach could be applied to these equations 
%in 2D which will be considered in the future. 
From a numerical point of view, there are two main 
directions of research worth exploring to further improve the code: first 
better preconditioners for GMRES adapted to the situations to be 
studied could increase the efficiency (which would be helpful in 
higher dimensions) and allow to study even more 
extreme situations which is mainly interesting from a theoretical 
point of view. Secondly one could improve the time integration by 
studying stiff integrators for PDEs with stiffness in the nonlinear 
part, see for instance~\cite{leja} and references therein.

\bibliographystyle{abbrv}
\bibliography{./Biblio}

\end{document}